%% file: main.tex
\documentclass[11pt,a4paper]{article}

\usepackage[margin=1in]{geometry}
\usepackage[T1]{fontenc}
\usepackage[utf8]{inputenc}
\usepackage{lmodern}
\usepackage{microtype}
\input{macros}
\usepackage{array}
\usepackage{booktabs}
\usepackage{makecell}
\usepackage{listings}
\usepackage{tikz}
\usetikzlibrary{arrows.meta,positioning,fit,calc}
\usepackage{graphicx}
\usepackage{xcolor}
\usepackage{enumitem}
\usepackage{csquotes}
\usepackage{authblk}

\usepackage[numbers,sort&compress]{natbib}
\usepackage{hyperref}
\usepackage{cleveref}
\crefname{lstlisting}{Listing}{Listings}
\Crefname{lstlisting}{Listing}{Listings}
\crefname{assumption}{Assumption}{Assumptions}
\Crefname{assumption}{Assumption}{Assumptions}
\crefname{proposition}{Proposition}{Propositions}
\Crefname{proposition}{Proposition}{Propositions}
\crefname{definition}{Definition}{Definitions}
\Crefname{definition}{Definition}{Definitions}
\crefname{remark}{Remark}{Remarks}
\Crefname{remark}{Remark}{Remarks}
\crefname{theorem}{Theorem}{Theorems}
\Crefname{theorem}{Theorem}{Theorems}

\hypersetup{
  colorlinks=true,
  linkcolor=blue!50!black,
  citecolor=blue!50!black,
  urlcolor=blue!50!black,
  pdftitle={Sequential Portfolio Optimizer (SPO-CPU and SPO-GPU): C++ Interior-Point Solvers and Differentiable Layers for Quadratic Cone Programs},
  pdfauthor={Fan Zhao, Yan Zaretskiy, Yuwen Chen, Jerry Shi}
}

\title{Sequential Portfolio Optimizer (\cpu{} and \gpu{}):\\
C++ Interior-Point Solvers and\\
Differentiable Layers for Quadratic Cone Programs}
\author[1]{Fan Zhao\thanks{Corresponding author. \gpu{} and its differentiable optimization layer were developed primarily while the author was at NVIDIA; present address: independent researcher. Email: \texttt{zhaofan816@gmail.com}}}
\author[1]{Yan Zaretskiy}
\author[1]{Yuwen Chen}
\author[1]{Jerry Shi}
\affil[1]{NVIDIA}
\date{}

\begin{document}
\maketitle

\begin{abstract}
Quadratic cone programs (QCPs) are a practical canonical form for
convex optimization, including parametric portfolio problems that
must be re-solved to high accuracy as market data change.
We present Sequential Portfolio Optimizer (\spo{}), with
\cpu{} and \gpu{} the C++ and CUDA realizations
of the Clarabel homogeneous-embedding interior-point method.
Clarabel is a registered trademark; \spo{} is our name for these
implementations.
Each exposes a forward optimizer and an in-solver differentiable
layer.
The Mehrotra algorithm is unchanged.
The contribution is an implementation stack: a reduced two-block
vector-Jacobian product with gradients in the original coordinates;
a documented variant of the Chen--Goulart--Jones smoothing warm
start, applied in equilibrated coordinates; and a GIL-free GPU
batch architecture whose isolated throughput is at parity with a
warm solver loop up to $B{=}16$.
On an RTX~A6000, subsequent solves of a $5000$-asset factor-model
QP average $24$\,ms on \gpu{}, compared with $98$\,ms on
Clarabel.rs and $122$\,ms on single-thread Gurobi; at $10{,}000$
assets the corresponding GPU time is $59$\,ms against $439$\,ms
for one-thread Mosek.
Clarabel.rs, \cpu{}, and \gpu{} agree pairwise to
$2.2{\times}10^{-7}$ in relative objective on this family.
Against a CVXPYLayers baseline driven by Clarabel, the GPU layer
reproduces the solution and the return-gradient on a $5000$-asset
factor QP/SOCP, while a subsequent forward-plus-backward step is
$101$--$137\times$ faster at $5000$ assets (and $17$--$30\times$
at $1000$ assets).
The smoothing warm start reduces QP rebalance iterations but can
increase them on a large GPU second-order cone.
Wall-clock rankings apply to the factor-model family only; an
$18$-problem Maros--M\'esz\'aros subset is reported as a
cross-check.
We do not claim a speedup over CuClarabel, and the batch
architecture is not a measured throughput improvement.
\end{abstract}



\tableofcontents
\newpage

\input{sections/01_intro}
\input{sections/02_background}

\input{sections/03_related}

\input{sections/04_overview}
\input{sections/05_cpu}

\input{sections/06_gpu}

\input{sections/07_warmstart}
\input{sections/08_diff}

\input{sections/09_setup}

\input{sections/10_evaluation}

\input{sections/11_discussion}

\input{sections/12_conclusion}

\section*{Acknowledgments}
We thank Jianchi Chen (Kendall Square Capital) for substantial
support on feature completion and bug reporting.
Private deployment reports are outside the measurement protocol
of this paper.

\section*{Use of AI tools}
The authors used the Kimi-3 model in Cursor to assist with drafting,
editing, and polishing the manuscript. The authors reviewed all
content and take full responsibility for it.

\bibliographystyle{plainnat}
\bibliography{refs}

\appendix
\input{sections/appendix}

\end{document}

%% file: macros.tex
\usepackage{amsmath,amssymb,amsthm,mathtools}
\usepackage{bm}

\newcommand{\R}{\mathbb{R}}
\newcommand{\Splus}{\mathbb{S}_{+}}
\newcommand{\K}{\mathcal{K}}
\newcommand{\Kstar}{\mathcal{K}^{*}}
\newcommand{\Kcirc}{\mathcal{K}^{\circ}}

\newcommand{\dataP}{P}
\newcommand{\dataA}{A}
\newcommand{\dataq}{q}
\newcommand{\datab}{b}

\newcommand{\proj}{\Pi}
\newcommand{\Dpi}{D\pi}
\newcommand{\Hs}{H_{s}}
\newcommand{\Ft}{F^{\top}}
\newcommand{\transp}{^{\top}}

\newcommand{\nncone}{\R_{+}}
\newcommand{\soc}{\mathcal{K}_{\mathrm{soc}}}
\newcommand{\expcone}{\mathcal{K}_{\mathrm{exp}}}
\newcommand{\powcone}{\mathcal{K}_{\mathrm{pow}}}
\newcommand{\genpow}{\mathcal{K}_{\mathrm{gp}}}

\newcommand{\tauvar}{\tau}
\newcommand{\kappavar}{\kappa}
\newcommand{\mubar}{\mu}

\newcommand{\spo}{SPO}
\newcommand{\cpu}{SPO-CPU}
\newcommand{\gpu}{SPO-GPU}
\newcommand{\diffsolver}{DiffSolver}
\newcommand{\batchdiff}{BatchDiffSolver}

\theoremstyle{plain}
\newtheorem{theorem}{Theorem}
\newtheorem{proposition}{Proposition}

\theoremstyle{definition}
\newtheorem{definition}{Definition}
\newtheorem{assumption}{Assumption}
\newtheorem{remark}{Remark}

\DeclareMathOperator{\diag}{diag}

\DeclareMathOperator{\nnz}{nnz}
\newcommand{\cdim}{\mathrm{dim}}

\DeclareMathOperator*{\minimize}{minimize}
\DeclareMathOperator*{\maximize}{maximize}
\DeclareMathOperator*{\subjectto}{subject~to}

%% file: sections/01_intro.tex
\section{Introduction}
\label{sec:intro}

Quadratic cone programs (QCPs) minimize a convex quadratic over the
intersection of an affine subspace and a product of convex cones
\cite{boyd2004convex,goulart2024clarabel}.
The form is rich enough to express linear and quadratic programs,
second-order cone programs (SOCPs), exponential- and power-cone models,
and semidefinite programs (SDPs), while remaining sparse enough that a
single interior-point code can serve as a backend for modelling
languages \cite{diamond2016cvxpy,agrawal2018rewriting}.
In quantitative investment the same template covers Markowitz
mean--variance allocation \cite{markowitz1952}, factor-model risk with
a quadratic or SOC constraint \cite{boyd2017portfolio,lobo2007portfolio},
turnover and leverage limits, and occasional entropy or SDP regularizers.
Those problems are \emph{parametric}: expected returns, factor loadings,
and risk aversion change every rebalance, while the sparsity pattern of
$(P,A)$ and the cone product $\K$ stay fixed.
A production solver must therefore (i)~reach tight feasibility and
duality-gap tolerances on each instance, (ii)~reuse symbolic
factorization across a stream of nearby instances, and (iii)~when the
QCP sits inside a learning loop, return vector-Jacobian products (VJPs)
of the solution map with respect to $(P,q,A,b)$ without rebuilding the
KKT pattern.
These three requirements organize the paper: parametric re-solve
latency (\Cref{sec:cpu,sec:gpu,sec:warmstart}), an in-solver
adjoint (\Cref{sec:diff}), and a GIL-free multi-problem GPU
architecture whose measured isolated throughput does not yet beat a
Python loop (\Cref{sec:batch,sec:exp-batch}).

First-order operator-splitting methods \cite{odonoghue2016scs,stellato2020osqp,garstka2021cosmo}
warm-start easily but are a poor match for high-accuracy daily or
intraday rebalancing: the tail convergence of ADMM-type schemes is slow
once residuals fall below $10^{-6}$--$10^{-8}$.
Primal-dual interior-point methods (IPMs) \cite{mehrotra1992,wright1997primaldual,nesterov1998primaldual}
remain the reliable choice at those tolerances, at the cost of a sparse
$LDL^\top$ factorization that grows increasingly ill-conditioned as
the iterate approaches the boundary of $\K$.
Clarabel \cite{goulart2024clarabel} specialized the homogeneous embedding
of a monotone complementarity problem to QCPs, keeping the quadratic
term in the objective rather than lifting it into an extra SOC, and
demonstrated a clear advantage on quadratic objectives relative to
linear-objective HSDE codes.
CuClarabel \cite{chen2025cuclarabel} ported the same algorithm to the GPU
with a mixed-family cone parallelization and NVIDIA cuDSS
\cite{cudss2024} as the direct KKT solver, and listed device-side
equilibration and KKT assembly as future work.
Neither ships an in-solver differentiable layer, nor a
multi-problem batch interface for training.

The computational cost of an IPM is concentrated in one sparse
symmetric KKT factorization per iteration (\Cref{sec:background}).
When the QCP is embedded in a learning loop, training requires not
only $x^\star$ but its sensitivity to the problem data.
Differentiating the KKT conditions, rather than unrolling the IPM
iterations, reduces that gradient to one additional linear solve
of the same size class, independent of the forward iteration count.

In quantitative research the prevailing workflow remains
two-stage: a statistical or neural model predicts expected returns
(and sometimes covariances), after which a separate portfolio QCP
computes the allocation \cite{markowitz1952,boyd2017portfolio}.
End-to-end training composes the predictor with that optimizer and
backpropagates a \emph{decision}-level loss---realized PnL, a
Sharpe-type objective, or turnover---through the solution map
\cite{amos2017optnet,agrawal2019cvxpylayers}.
The argument, which this paper does not re-evaluate on a trading
dataset, is that a forecast loss (for example the MSE on $\mu$) is
only loosely related to the allocated portfolio, whereas a loss on
$x^\star$ already incorporates the transaction-cost and tax
penalties encoded in the QCP.
The component that makes the loop possible is a differentiable
optimization layer: a network block whose forward pass is a
mathematical solver and whose backward pass differentiates the
interior-point KKT system at the computed optimum.
Unlike a feed-forward allocation head, the layer enforces hard
constraints---budget equalities, position and factor-exposure
bounds, and SOC risk caps---to solver tolerance, rather than as
soft penalties in the training loss.

Existing differentiable cone layers address this requirement only
in part.
CVXPYLayers \cite{agrawal2019cvxpylayers} differentiates linear cone
programs on the CPU through CVXPY canonicalization and a homogeneous
embedding \cite{agrawal2019differentiating}; every backward pass
crosses the host, and repeated solves typically rebuild the linear
system.
On a universe of a few thousand names that CPU path is often too
slow for training or fine-tuning (\Cref{sec:exp-cvxlayer}).
Interest in solver-in-the-loop research has therefore outpaced the
wall-clock performance of the layers that enable it.
\texttt{diffqcp} \cite{healey2025diffqcp} extends the implicit-function
argument to QCPs, including the power cone, and evaluates Jacobian
and vector-Jacobian products as an abstract linear operator on CPU
or GPU via JAX.
It is solver-agnostic by design: the forward QCP is an input, not a
fused IPM.
A single solver stack that (a)~solves the QCP with Clarabel's
embedding, (b)~differentiates the computed
$(x^\star,z^\star,s^\star)$ with a reduced two-block adjoint that
reuses sparsity, and (c)~on the GPU, batches many independent
same-pattern problems through native cuDSS uniform batches without
the Python GIL, is not available in those systems.
The GPU \diffsolver{} of this paper is designed for that setting:
\Cref{sec:exp-cvxlayer} reports a one-to-two-order subsequent-step
reduction versus CVXPYLayers on 5000-asset factor QP/SOCPs.
We do not claim an end-to-end learning-curve improvement; the
measurements concern the layer.

\paragraph{Contributions.}
This paper describes Sequential Portfolio Optimizer (\spo{}):
\cpu{} and \gpu{}, the CPU and CUDA realizations of
that stack, aimed at continuous portfolio optimization and at
differentiable QCP layers more generally.
Clarabel is a registered trademark.
We use \spo{} as the product name of our solvers;
they implement the published Clarabel homogeneous-embedding
interior-point method \cite{goulart2024clarabel} and are not the
Clarabel product.
The contributions are as follows.

\begin{enumerate}[leftmargin=1.4em,itemsep=0.35em]
\item \textbf{\cpu{}.}
A host interior-point solver that follows Clarabel.rs in algorithm,
cones, and defaults, with a SIMD-aware QDLDL factorization
\cite{stellato2020qdldl} extended by a multi-right-hand-side solve
that fuses the per-iteration affine and constant systems into one
$\mathrm{nrhs}{=}2$ pass over the factors, a SuiteSparse KLU
\cite{davis2011klu} adjoint for the unsymmetric two-block map
$F^\top$, and pooled
iteration workspaces so that a parametric re-solve does not
reallocate. A forward optimizer and a \diffsolver{} share the same
IPM; they are alternative interfaces, not two algorithms.

\item \textbf{\gpu{}.}
A CUDA solver in the lineage of CuClarabel that implements the
setup path the Julia paper deferred (device-side Ruiz equilibration
and two-phase CSR KKT assembly), and adds mixed-dimension PSD cones
with chordal decomposition \cite{vandenberghe2015chordal},
dimension-tiered SOC and generalized-power kernels (thread / warp /
block / cooperative multi-block, split at $d{=}8,64,1024$) with five
CUDA streams overlapping heterogeneous cone families, and a cuDSS
reuse protocol: cached elimination trees,
two-right-hand-side solves, hybrid host/device execute, and
iterative refinement gated by an Off--Watch--On state machine that
latches on only near tolerance (\Cref{sec:gpu-cudss}).
The KKT is kept in IEEE-754 binary64: CuClarabel's mixed-precision
experiments \cite[Table~6]{chen2025cuclarabel} are a deliberate
non-goal, because IPM systems become increasingly ill-conditioned.

\item \textbf{Smoothing-operator warm start.}
Both solvers implement Algorithm~1 of Chen, Goulart, and Jones
\cite{chen2025warmstart} \emph{with documented modifications}
(\Cref{sec:smoothing-impl}): the operator runs in equilibrated
coordinates; $\mubar_0$ uses feasibility residuals only; SOC uses
Theorem~4.1 ($\lambda=s_0^\star/z_0^\star$) rather than the Julia
prototype's dual form at $\lambda=1$; and the GPU additionally caps
$\mubar_{\mathrm{soc}}$.
Nonnegative, PSD, and three-dimensional cones use the dual proximal
form at smoothing parameter $\lambda=1$; the Newton damping
$\lambda^\star=2-\sqrt{3}$ in \Cref{sec:smoothing-impl} is a
different scalar.
The control is distinct from CVXPY's instance-cache warm start.
On the SOC risk-constrained family it is not uniformly helpful
(\Cref{sec:exp-warm}).

\item \textbf{In-solver \diffsolver{}.}
Given upstream $\partial L/\partial x$ (and optionally
$\partial L/\partial z$, $\partial L/\partial s$), a reduced two-block
adjoint $F^\top w=-g$ of size $n+m$ yields
$(\partial L/\partial q,\partial L/\partial b,\partial L/\partial P,\partial L/\partial A)$.
Large SOC blocks add at most three Woodbury columns to the adjoint
system rather than densifying it (the projection Jacobian itself is
rank at most two; \Cref{sec:dpi}).
Gradients are returned in the original (unequilibrated) coordinates.

\item \textbf{GIL-free GPU batch architecture (\batchdiff{}).}
Independent same-pattern problems are bound to stable
\diffsolver{} slots; compatible KKT phases are gathered into
cuDSS uniform batches whose shard width is independent of the
problem count.
This is a deployment architecture, not a measured speedup: on the
isolated portfolio QP a Python loop of warm \diffsolver{} objects
remains slightly faster up to $B=16$ (\Cref{sec:exp-batch}).
We do not claim a GIL-contention win; that experiment is not in
the protocol.
\cpu{} has no multi-problem QDLDL batch.
\end{enumerate}
\paragraph{Scope of claims.}
The Mehrotra homogeneous embedding, Nesterov--Todd scalings, and
the existence theory of the Chen--Goulart--Jones smoother are
inherited; we do not re-prove them.
The original mathematical content is
(i)~the reduced two-block adjoint
\eqref{eq:F}--\eqref{eq:grads} and its derivation
(\Cref{prop:vjp,app:reduced-jacobian});
(ii)~the documented variants of Algorithm~1, in particular the
GPU SOC cap, which is a different operator from
\cite[Theorem~4.1]{chen2025warmstart}
(\Cref{def:smoothing-impl});
and (iii)~the problem distance $d_{\mathrm{KKT}}$ as a
parameter-agnostic predictor of warm-start cost
(\Cref{def:dkkt,prop:dkkt-invar}).
Wall-clock comparisons, iteration ratios $R$, and
vector-Jacobian products checked against finite differences are
measured facts on the protocol of \Cref{sec:experiments}.
They are not claimed outside the factor-model family and the
$18$-problem Maros--M\'esz\'aros subset.
We do not claim a speedup over CuClarabel, an end-to-end
learning-curve improvement, or a throughput gain for
\batchdiff{}.

The paper is organized as follows.
\Cref{sec:background} recalls the QCP, the homogeneous embedding, and
the Mehrotra loop, and \Cref{sec:related} places the work against
Clarabel.rs, CuClarabel, the warm-start paper, and \texttt{diffqcp}.
\Cref{sec:overview} gives the system overview---two interfaces over
one IPM core (\Cref{fig:system})---and the design rules that recur
on both backends.
\Cref{sec:cpu,sec:gpu} describe the numerical and systems techniques
on each backend, with a measured effect where the protocol of
\Cref{sec:experiments} isolates one.
\Cref{sec:warmstart,sec:diff} treat warm start, the adjoint, and
batch training.
\Cref{sec:experiments} states the experimental protocol, and
\Cref{sec:evaluation} reports first-versus-subsequent timings and
accuracy for Clarabel.rs, \cpu{}, \gpu{}, Mosek, and Gurobi
\cite{gurobi} on
factor-model portfolios, together with warm start, \diffsolver{}, a
CVXPYLayers-as-benchmark layer comparison (\Cref{sec:exp-cvxlayer}),
and \batchdiff{} measurements on the same machine.
A Maros--M\'esz\'aros QP subset, including CuClarabel timings, is
reported in \Cref{sec:exp-maros}.
\Cref{sec:discussion,sec:conclusion} discuss CPU versus GPU regimes and
limitations.

%% file: sections/02_background.tex
\section{Background}
\label{sec:background}

We use the notation of Goulart and Chen \cite{goulart2024clarabel} and
do not re-derive the algorithm; this section records the objects used
in the rest of the paper.

\paragraph{Notation.}
Vectors are columns.
The Euclidean inner product is $x^\top y$.
For a closed convex cone $\K\subset\R^m$, the dual is
$\Kstar=\{z:z^\top s\ge 0\text{ for all }s\in\K\}$ and the polar is
$\Kcirc=-\Kstar$.
Euclidean projection onto a closed convex set $C$ is written
$\proj_C$.
A star denotes a computed primal-dual solution.
Other symbols are defined at first use.

\subsection{Quadratic cone programs}
\label{sec:qcp}

\begin{definition}[Quadratic cone program]
\label{def:qcp}
A QCP in standard form is
\begin{equation}
\label{eq:qcp}
\begin{aligned}
\minimize_{x,s}\quad
&\tfrac12 x\transp \dataP x + \dataq\transp x \\
\subjectto\quad
& \dataA x + s = \datab,\qquad s\in\K,
\end{aligned}
\end{equation}
with $x\in\R^n$, $s\in\R^m$, $\dataP\in\Splus^n$ (possibly zero),
$\dataA\in\R^{m\times n}$, $\dataq\in\R^n$, $\datab\in\R^m$, and $\K$ a
closed convex cone.
The dual is
\begin{equation}
\label{eq:qcp-dual}
\begin{aligned}
\maximize_{x,z}\quad
& -\tfrac12 x\transp \dataP x - \datab\transp z \\
\subjectto\quad
& \dataP x + \dataA\transp z = -\dataq,\qquad z\in\Kstar,
\end{aligned}
\end{equation}
where $\Kstar$ is the dual cone.
\end{definition}
Under strong duality the KKT system is
\begin{equation}
\label{eq:kkt}
\dataA x+s=\datab,\quad
\dataP x+\dataA\transp z=-\dataq,\quad
s\transp z=0,\quad
(s,z)\in\K\times\Kstar.
\end{equation}
At a primal-dual feasible point the duality gap equals
$s\transp z$; the solver terminates when that gap and the
residuals fall below $10^{-8}$.

\begin{assumption}[Standing hypotheses]
\label{ass:standing}
Throughout, $\dataP\in\Splus^n$, $\K\subset\R^m$ is a closed convex
cone of the product form in \Cref{sec:cones}, and \eqref{eq:qcp}
admits a primal-dual solution (strong duality).
Differentiability of the solution map is invoked only under the
additional hypotheses of \Cref{prop:vjp}.
\end{assumption}

\subsection{Homogeneous embedding}
\label{sec:hsde}

When $\dataP=0$, \eqref{eq:qcp} admits the classical homogeneous
self-dual embedding (HSDE) \cite{ye1994hsde,andersen2003hsde}.
Keeping a nonzero quadratic term yields a homogeneous embedding that is
\emph{not} self-dual \cite{goulart2024clarabel,odonoghue2021splitting}.
With slack variables $\tauvar,\kappavar\ge 0$ and the change of
variables $x\leftarrow x/\tauvar$, $z\leftarrow z/\tauvar$,
$s\leftarrow s/\tauvar$, one solves
\begin{equation}
\label{eq:embedding}
G(x,z,s,\tauvar,\kappavar)
=
\begin{bmatrix}
0\\ s\\ \kappavar
\end{bmatrix}
-
\begin{bmatrix}
\dataP & \dataA\transp & \dataq \\
-\dataA & 0 & \datab \\
-\dataq\transp & -\datab\transp & 0
\end{bmatrix}
\begin{bmatrix} x\\ z\\ \tauvar \end{bmatrix}
+
\begin{bmatrix}
0\\ 0\\ \tfrac1\tauvar x\transp \dataP x
\end{bmatrix}
=0
\end{equation}
over $(x,z,s,\tauvar,\kappavar)\in\R^n\times\Kstar\times\K\times\R_+\times\R_+$.
If $\tauvar^\star>0$, then $(x^\star/\tauvar^\star,s^\star/\tauvar^\star)$
(resp.\ $(x^\star/\tauvar^\star,z^\star/\tauvar^\star)$) solves
\eqref{eq:qcp} (resp.\ \eqref{eq:qcp-dual}).
If $\kappavar^\star>0$, the embedding produces a primal or dual
infeasibility certificate.
The pathological case $\tauvar^\star=\kappavar^\star=0$ is discussed in
\cite{goulart2024clarabel}.
Lifting the quadratic into an epigraphical SOC
\cite{goulart2024clarabel} is always possible, but it enlarges the cone
and typically slows the IPM on genuine QPs; both \cpu{} and \gpu{}
keep $\dataP$ in the KKT block as Clarabel does.
The same loop searches for an optimum
($\tauvar^\star>0$) and, when $\kappavar^\star>0$, returns a
Farkas-type infeasibility certificate.

\subsection{Mehrotra interior-point method}
\label{sec:ipm}

After optional chordal expansion of sparse PSD blocks
(\Cref{sec:cones}) and Ruiz equilibration \cite{ruiz2001scaling}
---row and column scaling of $(\dataP,\dataA)$ together with a cost
scaling $c>0$---the solver runs a Mehrotra predictor-corrector
\cite{mehrotra1992} on \eqref{eq:embedding}.
Each iteration (i)~forms residuals and a centrality measure
$\mubar$, (ii)~updates a Nesterov--Todd or dual scaling matrix on
each cone \cite{nesterov1998primaldual,goulart2024clarabel},
(iii)~assembles and factors a quasi-definite KKT matrix, (iv)~computes
an affine (predictor) step, a centering parameter $\sigma$, and a
combined (corrector) step, and (v)~takes a fraction-to-boundary step
with a barrier line search on nonsymmetric cones.
The Newton linear system has the symmetric saddle-point form
\begin{equation}
\label{eq:kkt-matrix}
\begin{bmatrix}
\dataP+\sigma I & \dataA\transp \\
\dataA & -\Hs
\end{bmatrix}
\begin{bmatrix} \Delta x \\ \Delta z \end{bmatrix}
=
\text{rhs},
\end{equation}
stored as an upper triangle for $LDL^\top$.
Here $\Hs$ is the (block-diagonal) scaled Hessian contribution of the
cones.
For a second-order cone with Nesterov--Todd scaling point $w$,
$\Hs=\eta^2(2ww^\top-J)$ with $J=\diag(1,-1,\ldots,-1)$, and this
factors further into the diagonal-plus-rank-two form
$\Hs=\eta^2\diag(d,1,\ldots,1)+\eta^2(uu^\top-vv^\top)$
\cite{goulart2024clarabel}, where the
scalar $d$ and the vectors $u,v$ are derived in closed form from $w$;
for a large generalized power cone the dual-scaled Hessian admits the
analogous diagonal-plus-rank-three form
$\Hs=\nu(D+pp^\top-qq^\top-rr^\top)$, where $\nu$ is the
barrier-parameter scale (distinct from the centrality $\mubar$ of
\Cref{sec:ipm}).
Rather than inserting such a dense $d\times d$ block into
\eqref{eq:kkt-matrix} once $d$ exceeds a threshold ($d>4$ for SOC,
$d>32$ for generalized power; the 3D primal-dual slice stays dense),
both solvers \emph{expand} the cone, bordering the KKT system with
the low-rank factors,
\begin{equation*}
\begin{bmatrix}
\dataP+\sigma I & \dataA\transp & 0 \\
\dataA & -\Hs^{d} & V \\
0 & V\transp & D_{\mathrm{aux}}
\end{bmatrix},
\end{equation*}
so that eliminating the auxiliary block recovers exactly $-\Hs$,
while the factorization sees only two (SOC) or three (generalized
power) dense columns per cone.
The expanded system has size $(n+m+p)\times(n+m+p)$ with
$p=2\,n_{\mathrm{spSOC}}+3\,n_{\mathrm{spGP}}$.
The auxiliary diagonal signs ($-,+$ for SOC; $+,-,-$ for generalized
power) are prescribed to the $LDL^\top$ inertia, and static and
dynamic regularization keep the factorization
quasi-definite \cite{goulart2024clarabel,stellato2020qdldl}.

The expansion is what makes large second-order cones affordable:
a dense $d\times d$ block would grow the fill-in like $d^2$ and
the factorization like $d^3$, whereas the bordered form adds two
rows and columns per cone and removes that dependence on $d$.
This is the regime of the risk cone in the portfolio SOCP of
\Cref{sec:portfolio-layer}.

Iterative refinement is optional; both backends provide an AUTO mode
whose Off--Watch--On gate switches refinement on only once the
iterate approaches tolerance and then latches it
(\Cref{sec:gpu-cudss}), so the refinement solves are not paid on every
iteration.

\subsection{Supported cones}
\label{sec:cones}

The product $\K=\K_1\times\cdots\times\K_p$ may contain, in caller
order,
\begin{itemize}[leftmargin=1.2em,itemsep=0.15em]
\item the zero cone $\{0\}^d$ (equalities);
\item the nonnegative orthant $\nncone^d$;
\item second-order cones
      $\soc^d=\{(t,u):\|u\|_2\le t\}$ \cite{alizadeh2003socp};
\item exponential cones
      $\expcone=\{(x,y,z): y\,e^{x/y}\le z,\; y>0\}$ \cite{serrano2015exp};
\item three-dimensional power cones
      $\powcone^\alpha=\{(x,y,z): x^\alpha y^{1-\alpha}\ge |z|\}$
      \cite{chares2009power};
\item generalized power cones
      $\genpow=\{x\in\R^{\cdim_1},\,u\in\R^{\cdim_2}:
      \prod_{i=1}^{\cdim_1} x_i^{\alpha_i}\ge\|u\|_2,\ x\ge0\}$,
      $\alpha_i>0$, $\sum_i\alpha_i=1$;
\item PSD cones in scaled lower-triangular (\texttt{svec}) coordinates.
\end{itemize}
For a generalized power cone the solver either splits it into
three-dimensional power cones (the default) or, on the
$(\cdim_1,\cdim_2)=(2,1)$ slice, keeps it and uses the native
high-dimensional dual-scaled IPM.
Nonsymmetric cones use dual scaling, with a primal-dual BFGS Hessian
on three-dimensional power and on the $(\cdim_1,\cdim_2)=(2,1)$
slice of $\genpow$, matching Clarabel.jl.
Sparse PSD constraints may be replaced by a chordal clique
decomposition with parent--child or clique-graph merge
\cite{fukuda2001exploiting,vandenberghe2015chordal}; the original
$(x,s,z)$ is recovered by a Vandenberghe completion after the IPM.
Barriers used by the warm-start operator are listed in
\Cref{app:barriers}.
Modelling languages reduce entropy and $p$-norm constraints to many
three-dimensional cones, so per-copy latency matters; a PSD cone is
the most expensive primitive per dimension
(\Cref{sec:gpu-setup}).

%% file: sections/03_related.tex
\section{Related work and lineage}
\label{sec:related}

\Cref{tab:lineage} summarizes the six systems against which we compare.
The rest of this section states, for each reference, the
property we inherit and the difference realized here.
The paper is not a new interior-point algorithm: the Mehrotra
homogeneous embedding is that of Clarabel.
Originality, in the sense of \Cref{sec:intro}, lies in the reduced
adjoint, the documented smoothing variants, the problem distance,
and the fused CPU/GPU implementation stack, not in a new embedding.

\input{tables/tab_lineage}

\paragraph{Clarabel.}
Goulart and Chen \cite{goulart2024clarabel} give the homogeneous
embedding \eqref{eq:embedding}, the Mehrotra loop, NT and dual
scalings, chordal PSD support, and a numerical comparison against
Mosek, ECOS, SCS, and others, with a particular advantage on quadratic
objectives.
The reference implementation is Clarabel.rs \cite{clarabelrs}.
Clarabel is a registered trademark; \cpu{} and \gpu{} are the
\spo{} implementations of that published method, not the Clarabel
product.
\cpu{} follows the same mathematics on the host, so
the performance difference is almost entirely a matter of
implementation (SIMD QDLDL, pooled workspaces, and OpenMP residuals), plus two
features that Clarabel.rs does not provide: the KLU adjoint of
\Cref{sec:diff} and the smoothing warm start of \Cref{sec:warmstart}.
There is no multi-problem batch on either CPU solver; parametric work is
a sequential update of $q$ or $b$.

\paragraph{CuClarabel.}
Chen et al.\ \cite{chen2025cuclarabel} introduce a mixed parallel
computing strategy: linear (zero/nonnegative) cones are processed
first, then each remaining family (SOC, exp, power, PSD) runs on its
own GPU stream, with cones reordered so that a family occupies a
contiguous memory range.
The KKT solve uses cuDSS; mixed-precision factorization with iterative
refinement can reduce time on some dense-constraint QPs at some cost
in robustness.
The paper's conclusion defers GPU equilibration and KKT population,
PSD batching is restricted to same-dimension blocks, and
there is no differentiable layer.
\gpu{} keeps the mixed-family strategy and cuDSS, implements the
deferred setup path on device, generalizes PSD to mixed dimensions via
chordal merge, and adds the cuDSS settings of \Cref{sec:gpu-cudss}
(\texttt{nrhs=2}, hybrid execute, pooled device memory, latched AUTO IR,
unsplit large SOCs) together with \diffsolver{} and
\batchdiff{}.
It does \emph{not} adopt mixed-precision KKT: the IPM Jacobian is
factored in binary64.

\paragraph{Warm start.}
Classical interior-point warm starts include Gondzio's multiple
centrality corrections \cite{gondzio1996mcc}, the
Y{\i}ld{\i}r{\i}m--Wright neighbourhood analysis
\cite{yildirim2002warm}, and the Gondzio--Grothey reoptimization
framework \cite{gondzio2003reopt}.
Those methods keep a previous interior iterate and restore
centrality by extra correctors or a controlled residual.
Chen, Goulart, and Jones \cite{chen2025warmstart} (arXiv:2512.00693,
posted 30 November 2025) replace Euclidean
projection by a smoothing operator---the proximal map of a
logarithmically homogeneous self-concordant barrier---and prove that
Algorithm~1 with $\lambda=1$ produces a point on the central path
whose residuals remain comparable to those of the previous optimum for
nonnegative, second-order, and PSD cones.
The Julia prototype lives on a Clarabel.jl branch
\cite{cuclarabeljl} and applies the
dual form at $\lambda=1$ for every cone, including SOC.
Both of our solvers implement the operator in equilibrated
coordinates, with Theorem~4.1 for SOC, the Moreau rule for the zero
cone (\Cref{sec:warmstart}), and GPU kernels for every supported
family, including PSD.

\paragraph{diffqcp and CVXPYLayers.}
Agrawal et al.\ \cite{agrawal2019differentiating,agrawal2019cvxpylayers}
differentiate linear cone programs through a homogeneous embedding and
Minty's parameterization of the complementarity set.
Healey, Nobel, and Boyd \cite{healey2025diffqcp} extend the argument to
QCPs, give JVP/VJP algorithms via cone projection derivatives, and
implement them as a JAX operator that can run on a GPU without a host
round-trip---addressing the long-standing data-movement overhead of
CVXPYLayers.
\texttt{diffqcp} remains the clean solver-agnostic theory reference.
DiffOpt.jl \cite{besancon2022diffopt} is the modelling-language
counterpart: it differentiates through JuMP models and is not a
fused Clarabel adjoint.
CVXPYLayers driven by Clarabel is the
production-stack baseline of \Cref{sec:exp-cvxlayer}: the comparison
is stack-versus-stack (canonicalization plus a CPU cone solve plus
the \texttt{diffcp} adjoint), not IPM-versus-IPM.
\diffsolver{} solves a \emph{reduced} two-block adjoint of size
$n+m$ at a Clarabel solution, caches cuDSS or KLU factors across
training steps, and, on the GPU, batches many QCPs.
It is not a drop-in replacement for the homogeneous Minty map: we drop
the $\tauvar$ row and column for conditioning, and we never form an
abstract Jacobian operator independent of this IPM.

\paragraph{Other solvers.}
Mosek \cite{mosek2024} is the commercial baseline for high-accuracy
conic solves in \Cref{sec:experiments}.
ECOS \cite{domahidi2013ecos} and SCS \cite{odonoghue2016scs} appear as
optional SOCP baselines.
None of these expose a fused QCP VJP with multi-problem GPU batching.

%% file: tables/tab_lineage.tex
\begin{table}[t]
\centering
\caption{Lineage of QCP solvers and differentiable layers. ``---'' means the capability does not arise; ``n/a'' means the column's design does not expose it.
``3D split or native'' is the choice described in \Cref{sec:cones}: split a large generalized-power cone into 3D cones, or keep it and use the native high-dimensional IPM.
\cpu{} and \gpu{} share the Clarabel homogeneous embedding; mixed-precision KKT is a
deliberate non-goal.}
\label{tab:lineage}
\scriptsize
\setlength{\tabcolsep}{2.4pt}
\resizebox{\textwidth}{!}{%
\begin{tabular}{@{}lcccccc@{}}
\toprule
Capability &
\makecell{Clarabel.rs\\{\scriptsize\cite{goulart2024clarabel}}} &
\makecell{CuClarabel\\{\scriptsize\cite{chen2025cuclarabel}}} &
\cpu{} &
\gpu{} &
\makecell{diffqcp\\{\scriptsize\cite{healey2025diffqcp}}} &
\makecell{CVXPYLayers\\{\scriptsize\cite{agrawal2019cvxpylayers}}} \\
\midrule
Language & Rust & Julia & C++ & CUDA C++ & JAX/Python & Python \\
KKT backend & QDLDL & cuDSS & QDLDL+KLU & cuDSS & host/GPU LA & SCS (diffcp) \\
GPU equilibration / KKT fill & --- & host (deferred) & --- & device & --- & --- \\
Mixed-precision KKT & no & yes & no & no (FP64) & n/a & no \\
Mixed-dim PSD / chordal & yes & same-dim PSD & yes & yes & n/a & via CVXPY \\
GenPower cone & GenPower & no & 3D split or native & 3D split or native & n/a & via DCP \\
Smoothing warm start & no & no & yes & yes & --- & --- \\
In-solver \diffsolver{} & no & no & yes & yes & post-hoc op. & yes (CPU) \\
Multi-problem cuDSS batch & no & no & no & \batchdiff{} & no & no \\
CVXPY DPP instance reuse & yes & yes & yes & yes & no & compile once \\
Zero-copy Torch/CuPy & no & limited & no & $q,b$ and grads & JAX device & CPU arrays \\
\bottomrule
\end{tabular}%
}
\end{table}

%% file: sections/04_overview.tex
\section{System overview}
\label{sec:overview}

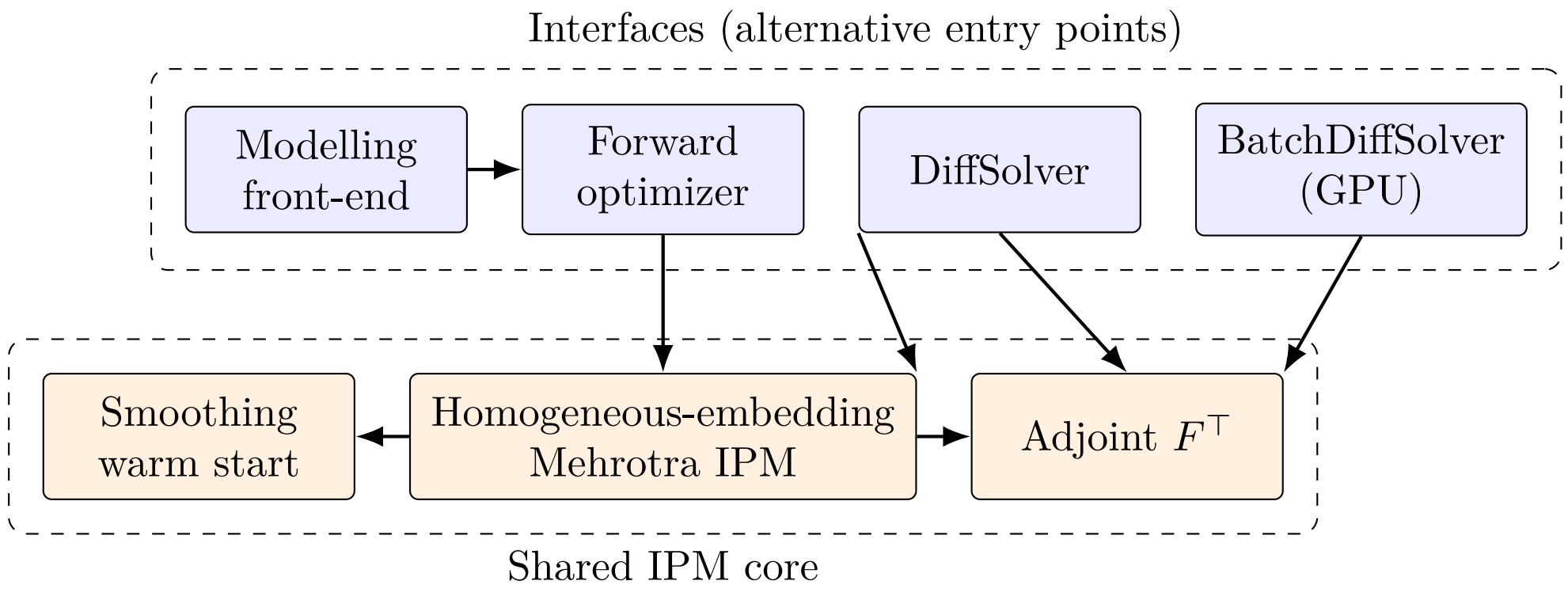
\begin{figure}[t]
\centering
\input{figs/fig_system}
\caption{Two interfaces share one interior-point core.
They are alternative entry points, not stages of one pipeline:
a modelling front-end drives the forward optimizer (IPM only);
\diffsolver{} calls the IPM and the adjoint $F^\top$;
\batchdiff{} (GPU only) owns a pool of \diffsolver{} slots.}
\label{fig:system}
\end{figure}

\cpu{} and \gpu{} are one algorithm on two backends: the Clarabel
homogeneous-embedding Mehrotra IPM of \Cref{sec:background},
instantiated on the host (\Cref{sec:cpu}) and on the device
(\Cref{sec:gpu}).
Around that single numerical core the stack provides three
production capabilities, each developed in its own section:
parametric re-solves with a reused symbolic factorization
(\Cref{sec:cpu,sec:gpu}), a smoothing-operator warm start
(\Cref{sec:warmstart}), and an in-solver differentiable layer with
a GIL-free GPU batch architecture (\Cref{sec:diff}).
\Cref{fig:system} shows the layout; the rest of this section
summarizes the interfaces before the per-backend detail.

\subsection{Two interfaces, one IPM}
\label{sec:dual-compile}

A modelling rebalance requires a conic optimizer, not an adjoint; a
training job requires a vector-Jacobian product and must not
recanonicalize on the inner loop.
Each backend therefore exposes two interfaces over the same IPM,
a forward solver and a \diffsolver{}.
The split is operational rather than algorithmic---one KKT assembly, one
cone set, one equilibration path.
Solve and backward drop the Python GIL; \batchdiff{} materializes
inputs once and then runs the batch without returning to the
interpreter.
GPU adjoint contexts are created on first backward, so a
forward-only process does not allocate them.

\subsection{Design rules}
\label{sec:overview-rules}

Three rules recur in both backends and organize the design
sections that follow.
\begin{enumerate}[leftmargin=1.4em,itemsep=0.35em]
\item \textbf{Structure is computed once.}
Everything that depends only on the sparsity pattern---cone layout,
equilibration maps, fill-reducing orderings, elimination trees---is
built at construction and reused across the parametric stream;
the per-iteration loop touches numeric values only.
\item \textbf{The device owns the hot path.}
On the GPU, equilibration, KKT assembly, cone projections, and the
adjoint all run on-device; the host orchestrates but never
round-trips a matrix.
\item \textbf{Binary64 everywhere.}
The KKT is kept in IEEE-754 binary64 on both backends;
mixed-precision factorizations are a deliberate non-goal because
IPM systems become increasingly ill-conditioned near the boundary
(\Cref{sec:gpu}).
\end{enumerate}

%% file: figs/fig_system.tex
\begin{tikzpicture}[
  font=\small,
  box/.style={draw, rounded corners=2pt, align=center, inner sep=5pt, minimum height=1.05cm},
  py/.style={box, fill=blue!8, minimum width=2.35cm},
  core/.style={box, fill=orange!12, minimum width=2.6cm},
  arr/.style={-{Latex}, thick}
]
  \node[py] (cvxpy) {Modelling\\front-end};
  \node[py, right=0.45cm of cvxpy] (fwd) {Forward\\optimizer};
  \node[py, right=0.45cm of fwd] (ds) {\diffsolver{}};
  \node[py, right=0.45cm of ds] (bds) {\batchdiff{}\\(GPU)};

  \node[core, below=1.15cm of fwd] (ipm) {Homogeneous-embedding\\Mehrotra IPM};
  \node[core, left=0.45cm of ipm] (ws) {Smoothing\\warm start};
  \node[core, right=0.45cm of ipm] (adj) {Adjoint $\Ft$};

  \draw[arr] (cvxpy) -- (fwd);
  \draw[arr] (fwd) -- (ipm);
  \draw[arr] (ds.south) -- (adj.north);
  \draw[arr] (ds.south west) -- (ipm.north east);
  \draw[arr] (bds.south) -- (adj.north east);
  \draw[arr] (ipm) -- (ws);
  \draw[arr] (ipm) -- (adj);

  \node[draw, dashed, rounded corners=4pt, inner sep=8pt,
        fit=(cvxpy)(fwd)(ds)(bds), label=above:{Interfaces (alternative entry points)}] {};
  \node[draw, dashed, rounded corners=4pt, inner sep=8pt,
        fit=(ws)(ipm)(adj), label=below:Shared IPM core] {};
\end{tikzpicture}

%% file: sections/05_cpu.tex
\section{\cpu{}: numerical techniques}
\label{sec:cpu}

\cpu{} is a host interior-point solver for \eqref{eq:qcp}.
The algorithm, cones, and default tolerances follow Clarabel.rs
\cite{goulart2024clarabel}; the Mehrotra homogeneous embedding is
unchanged.
This section records the numerical and systems choices that
matter on parametric and differentiable workloads, each with a
measured effect where \Cref{sec:experiments} isolates one.

\subsection{Solver structure}
\label{sec:cpu-arch}

Construction optionally chordally decomposes sparse PSD cones,
collapses adjacent nonnegative cones (caller order is otherwise
preserved), applies Ruiz equilibration, and runs a single AMD
ordering \cite{amestoy2004amd} plus QDLDL analysis
\cite{stellato2020qdldl}.
Unlike \gpu{}, there is no type-sort permutation of
$(\dataA,\datab)$: no batched kernels require contiguous family
runs on the host.

Each IPM iteration updates residuals and cone scalings, scatters the
numeric values of \eqref{eq:kkt-matrix} through a precomputed
permutation of size $O(\nnz)$, then refactorizes and solves.
Nonsymmetric cones start from a unit point and may fall back from
primal-dual to dual-only scaling; a pure centering recovery is used
when a primal-dual step becomes too short, matching Clarabel.jl.
Updates of $q$ or $b$ rewrite the original-space vectors, re-apply
the cached Ruiz diagonals, and invalidate numeric factors while
keeping the AMD permutation.

The organizing principle is a strict separation between
\emph{one-time structural analysis} and the \emph{per-iteration
numeric loop}.
Everything that depends only on the sparsity pattern---cone layout,
equilibration maps, the fill-reducing AMD ordering, and QDLDL's
symbolic elimination tree---is computed once at construction,
because a parametric portfolio loop re-solves thousands of
instances that differ only in numerical values, never in sparsity
pattern.
An iteration then updates only floating-point values: they are
gathered into the KKT, the numeric factor is recomputed, and the
system is solved.
This separation is the principal design choice for re-solve
latency; the remainder of the section refines it.

\subsection{Forward hot path}
\label{sec:cpu-hot}

Per iteration the numeric path is a value gather, a numeric
factorization, and a triangular solve.
Symbolic analysis is not repeated.
Each iteration needs three solves against the same $LDL^\top$
factors: the constant system $K[x_2;z_2]=[-\dataq;\datab]$, the
affine (predictor) right-hand side, and the combined (corrector)
right-hand side.
The constant system is invariant across iterations and independent
of the affine one, so the two are fused into a single
two-right-hand-side solve; the combined step then reuses the
$(x_2,z_2)$ computed in that call, and three logical systems cost
two triangular-solve invocations per iteration.
Upstream QDLDL exposes only a single-vector solve; we added a
multi-right-hand-side entry point that, at $\mathrm{nrhs}{=}2$,
streams the stored $L$ once and applies each saxpy update to both
columns, falling back to per-column solves beyond two right-hand
sides.
QDLDL uses Clarabel-style dynamic pivot regularization and the
inertia signs required by the sparse-SOC expansion columns.
When a supernodal column has enough nonzeros, the triangular solve
uses a gather/FMA kernel; short columns stay scalar, since the
gather setup is not recovered there.
The instruction set is selected at runtime (scalar through AVX-512,
and NEON), so a single binary remains correct on hosts without the
widest vector unit.
We did not isolate the CPU $\mathrm{nrhs}{=}2$ saving in
\Cref{sec:experiments}: the portfolio KKT is factorization-bound,
and the GPU ablation of the same fusion is $11\%$
(\Cref{tab:ablate}).

\subsection{Techniques for parametric loops}
\label{sec:cpu-opts}

None of the following changes the algorithm; each removes an
overhead that a direct port of Clarabel.rs would carry.

\paragraph{Vectorization.}
Dense length-$n$ kernels (residual dots, SOC inner products, diagonal
scaling) follow an explicit SIMD path rather than the compiler's
cost model; equilibration norms, residual reductions, and CSR SpMV
use OpenMP, while cone scaling stays serial because per-block
Hessian caches are not thread-safe.
BLAS threads are pinned before PSD eigendecompositions so that a
host BLAS does not oversubscribe against the residual OpenMP region.
The measured effect on the portfolio QPs is a small
\emph{regression}: forcing scalar code on the $n{=}5000$ instance of
\Cref{tab:ablate} was faster than the default SIMD dispatch
(122\,ms versus 132\,ms), because the dense kernels the path
vectorizes are a small fraction of QDLDL on this sparsity and the
gather/dispatch prologue costs about 8\%.
Models dominated by a large dense cone block shift the work into
exactly those kernels, and \Cref{tab:simd-dense} measures that
regime directly: on a max-eigenvalue SDP with one PSD block of
dimension $100$ (svec length $5050$), runtime dispatch beats
forced scalar by $1.29\times$ per iteration ($7.51$\,s versus
$9.70$\,s, six iterations each).
On the portfolio SOCPs, whose largest block is a SOC of dimension
$5051$ or $20051$, the same comparison stays within run-to-run
scatter ($\pm 3\%$), because QDLDL still dominates those models.
We therefore claim vectorization as a measured contribution for
PSD-dominated workloads only; it is retained elsewhere because
runtime dispatch keeps the binary portable.

\input{tables/tab_simd_dense}

\paragraph{Pooled workspaces.}
A parametric loop calls the IPM thousands of times on the same
sparsity pattern, so per-iteration allocation is pure overhead: it
adds latency variance and allocator contention under solver threads.
Host arrays owned by a solver are therefore taken from a pool, and
same-size iteration workspaces are reused rather than allocated.
There is no user-facing setting.
We do not isolate this effect in \Cref{sec:experiments}; its weight
is concentrated at small $n$, where allocator noise is a visible
fraction of a sub-millisecond solve.

\paragraph{Correctly rounded reductions.}
IPM reductions must remain reproducible at $10^{-8}$ tolerances, so
the build does not enable unsafe floating-point contraction, and
vectorized reductions are pairwise within a register plus a scalar
tail.
Across the QP rows of \Cref{tab:portfolio-acc}, the pairwise
relative objective difference between \cpu{} and Clarabel.rs stays
at most $5.3\times10^{-8}$ (at $n{=}10000$); including SOCPs the
three Clarabel solvers stay within $2.2\times10^{-7}$
(\Cref{sec:experiments}).
AUTO iterative refinement likewise stays off until the iterate is
close to tolerance (the gating state machine is specified in
\Cref{sec:gpu-cudss}): refinement steps are wasted while the
residuals are far from converged, and the well-scaled portfolio
QPs of \Cref{sec:experiments} never need them.

\subsection{Delta versus Clarabel.rs}
\label{sec:cpu-delta}

\input{tables/tab_cpu_vs_rs}

\Cref{tab:cpu-vs-rs} summarizes the differences.
The algorithm is Clarabel's, and the numerical equivalence of
\Cref{tab:portfolio-acc} is the evidence that the port is faithful.
SIMD triangular solves and pooled workspaces target per-iteration
overheads that matter in a re-solve loop but are invisible in a
single solve; on the factorization-bound portfolio QP their measured
effect is negligible (\Cref{sec:cpu-opts}), and we report that
rather than claim an isolated speedup.
Struct-of-arrays layout for three-dimensional power cones serves a
modelling need: CVXPY reduces high-dimensional power constraints
to many tiny cones, and a batched layout runs their scaling and
projection as dense vector operations instead of per-cone dispatch.
The KLU adjoint exists because the differentiable layer of
\Cref{sec:diff} needs an unsymmetric factorization that QDLDL cannot
provide; placing it in the same solver lets forward and backward
share cones and equilibration maps.

The deliberate limitation is the absence of a multi-problem QDLDL
batch.
CPU training is either a sequential parameter-update loop on one
\diffsolver{}, or process-level parallelism outside the solver.
That split is intentional: QDLDL is a single-matrix sparse direct
solver, so a batched facade would serialize on one AMD permutation
without a GPU to hide it.

%% file: tables/tab_simd_dense.tex
\begin{table}[t]
\centering
\caption{Dense-cone SIMD benchmark (CPU, one thread). Runtime SIMD
dispatch versus a forced scalar code path
on three dense-block instances: the portfolio SOCPs, whose largest
block is a second-order cone of dimension $n{+}k{+}1$, and a
max-eigenvalue SDP with one PSD block of dimension $100$ (svec length
$5050$), where the PSD scaling application dominates the iteration.
Times are medians of subsequent cold solves (ten, five, and three
repeats respectively, seed~42, tolerances $10^{-8}$); iterations are
identical across modes.}
\label{tab:simd-dense}
\small
\setlength{\tabcolsep}{5pt}
\begin{tabular}{@{}llrrrrr@{}}
\toprule
 & Dominant & \multicolumn{2}{c}{ms / iteration} & \multicolumn{2}{c}{total solve} & \\
\cmidrule(lr){3-4}\cmidrule(lr){5-6}
Instance & block & SIMD & scalar & SIMD & scalar & speedup \\
\midrule
SOCP $n{=}5000$  & SOC $d{=}5051$  & 6.3 & 6.5 & 139.0\,ms & 143.6\,ms & $1.03\times$ \\
SOCP $n{=}20000$ & SOC $d{=}20051$ & 28.9 & 28.1 & 723.2\,ms & 701.4\,ms & $0.97\times$ \\
Max-eig SDP      & PSD dim $100$   & 7512 & 9702 & 45.1\,s & 58.2\,s & $\mathbf{1.29\times}$ \\
\bottomrule
\end{tabular}
\end{table}

%% file: tables/tab_cpu_vs_rs.tex
\begin{table}[t]
\centering
\caption{Clarabel.rs versus \cpu{}. The Mehrotra homogeneous embedding is unchanged;
differences are numerical techniques, the in-solver adjoint, and warm start.}
\label{tab:cpu-vs-rs}
\small
\begin{tabular}{@{}llp{5.6cm}@{}}
\toprule
Topic & Clarabel.rs & \cpu{} \\
\midrule
Language & Rust & C++ \\
Forward $LDL^\top$ & QDLDL & QDLDL with SIMD triangular solves \\
Adjoint $F^\top$ & none & KLU (default) or quasi-definite lift \\
Workspaces & per-solve allocator & pooled, reused across re-solves \\
Many 3D power cones & enum dispatch & struct-of-arrays batch \\
$A,b$ permutation & none required & none (caller cone order) \\
Multi-problem batch & no & no (reuse one solver) \\
\bottomrule
\end{tabular}
\end{table}

%% file: sections/06_gpu.tex
\section{\gpu{}: numerical techniques}
\label{sec:gpu}

\gpu{} is a CUDA interior-point solver for the same QCP
\eqref{eq:qcp}.
The IPM is CuClarabel's Algorithm~1 \cite{chen2025cuclarabel}; the
systems contribution is a device-resident setup path, an extended
cuDSS reuse protocol, mixed-dimension PSD cones, and the
differentiable stack of \Cref{sec:diff}.

\subsection{Mapping onto CuClarabel}
\label{sec:gpu-map}

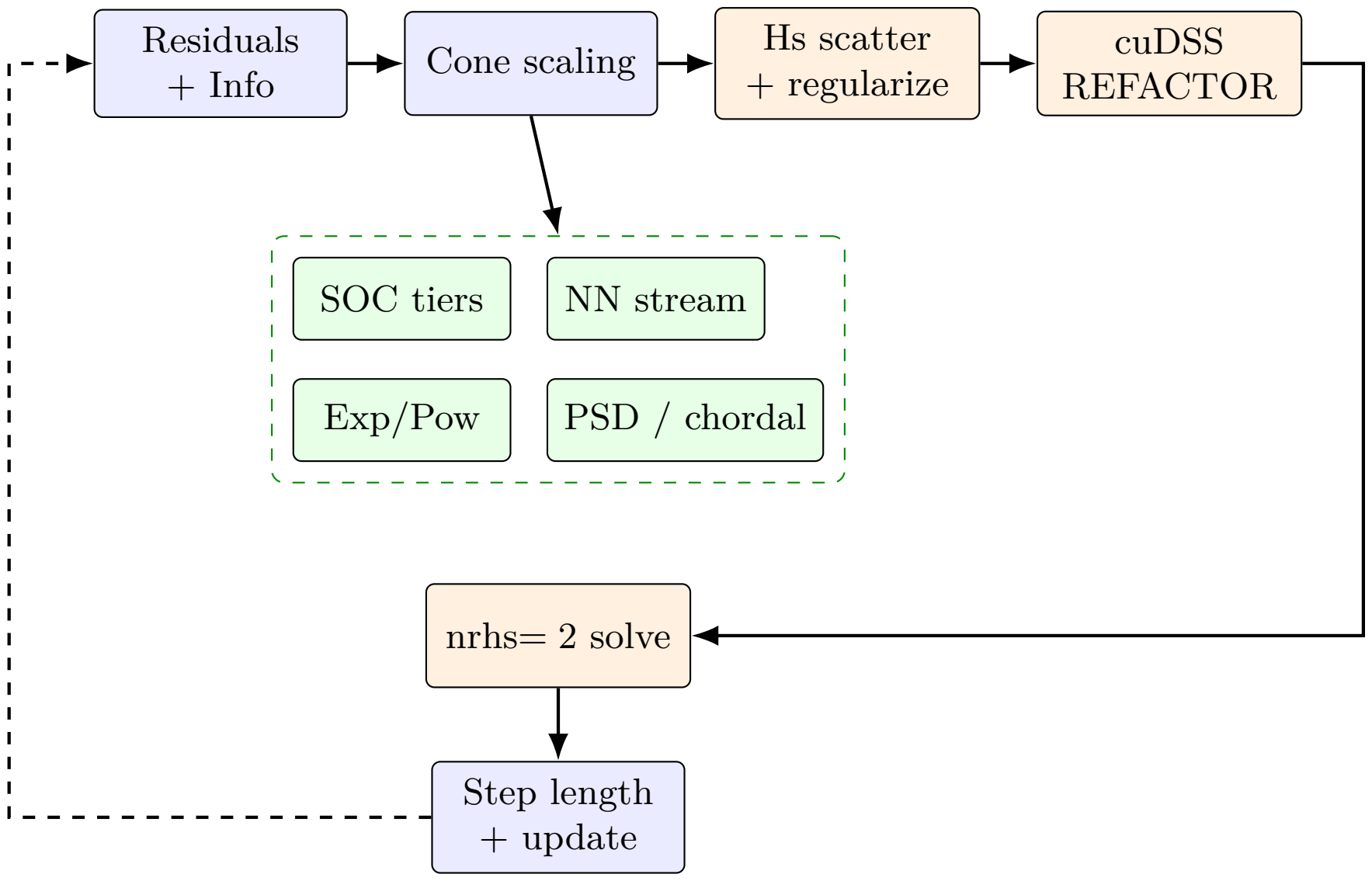
\begin{figure}[t]
\centering
\input{figs/fig_gpu_ipm}
\caption{One \gpu{} iteration.
Cone families run on separate streams after a type-sort of $(A,b)$.
cuDSS stays in refactorization across IPM iterations once analysis
has run; the affine residual and the constant block share a
two-right-hand-side solve.
The dashed edge is the iterate update.}
\label{fig:gpu-ipm}
\end{figure}

Setup still comprises data upload, equilibration, and solver
initialization; the solve phase is the Mehrotra loop of
\Cref{sec:ipm}.
Cones are \emph{reordered by family} so that zero, nonnegative, SOC,
exponential, power, generalized power, and PSD blocks occupy
contiguous device ranges---the mixed parallel computing strategy of
\cite{chen2025cuclarabel}.
Five dedicated CUDA streams are assigned to the nonlinear families
(SOC, exponential, power, generalized power, PSD); the linear cones
stay on the solver stream.
Multi-stream execution is enabled when at least two nonlinear
families are present and together contribute at least eight cones;
an event-based fan-out/fan-in against the main stream then lets
different families run their scaling, $\Hs$ assembly, shift, and
projection work concurrently.
A synchronization is required only where a later IPM step consumes a
global reduction (step length, barrier).
Matrices are stored in CSR; the KKT matrix is the full (not purely
triangular) pattern expected by cuDSS, with an upper-triangular
\emph{view} for the forward symmetric factorization.
\Cref{fig:gpu-ipm} sketches one iteration.

\subsection{CUDA math libraries}
\label{sec:gpu-stack}

The KKT solve, sparse matvecs, dense updates, and PSD
eigendecompositions are CUDA math-library kernels (cuDSS
\cite{cudss2024}, cuSPARSE \cite{cusparse2024}, cuBLAS
\cite{cublas2024}, cuSOLVER \cite{cusolver2024}).
In this paper, ``\gpu{} optimization'' refers primarily to
eliminating redundant work \emph{between} those kernels---fusion,
cached permutations, phase reuse, and a single device memory
pool---rather than to handwritten replacements for them.
Accordingly, \Cref{sec:gpu-cudss} measures most items through
configuration ablations rather than kernel rewrites.

\subsection{Closing CuClarabel's deferred setup path}
\label{sec:gpu-setup}

CuClarabel's conclusion notes that KKT population and matrix
equilibration ``take non-negligible time'' and leaves moving them to
the GPU as future work \cite{chen2025cuclarabel}.
\gpu{} implements both.
A host-side setup forces a device--host round-trip between
construction and the first solve, which dominates first-execution
latency on small and medium instances; the same maps are re-applied
on every parametric update (\Cref{sec:gpu-dpp}), so leaving them on
the host would add that cost to every re-solve.

\paragraph{Device Ruiz equilibration.}
Column norms, clamps, and the scalings of
$\dataP$, $\dataA$, $\dataq$, $\datab$ run entirely on device, with
the inner replace-zero / inverse-square-root / clip / accumulate
fused (eight launches collapsed to two) and a single synchronize at
the end of the loop.
The ablation of \Cref{tab:ablate} measures equilibration itself
rather than its placement: disabling it on the $n{=}5000$ QP leaves
the iteration count unchanged ($13$) and the steady solve within
run-to-run noise ($24.6$\,ms default versus $24.1$\,ms with Ruiz
off).
We did not retain a host-side switch, so the device-placement benefit
is architectural---no round-trip in the iteration loop---rather than
an isolated measurement.

\paragraph{Two-phase CSR KKT assembly.}
A first kernel counts row nonzeros; a second fills values.
Scatter indices from $(P,A)$ into the KKT are stored once, so
subsequent IPM iterations never rebuild structure: they scatter
$O(\nnz)$ values in place and mark the numeric factor dirty.
This is the device analogue of the CPU value gather
(\Cref{sec:cpu-hot}), and it makes per-iteration KKT updates
proportional to the data rather than to the pattern.

\paragraph{Cross-solve symbolic reuse.}
Nested dissection on the CPU is a one-time cost per sparsity pattern.
After the first analysis, the fill-reducing permutation and
nested-dissection tree are cached and re-injected, skipping the host
nested dissection on later rebuilds.
A values-only update leaves the next solve as a refactorization; a
fused analysis-plus-factor at construction (default on) moves lazy
kernel loading and workspace allocation out of the first solve.
The measured effect is on first execution: with prefactoring the
setup phase absorbs $96$\,ms and the first solve returns in
$\sim$\,$25$\,ms, whereas disabling it defers the cost to a
$102$\,ms first solve and leaves setup at $12$\,ms
(\Cref{tab:ablate}); subsequent solves are unaffected.

\paragraph{Cones beyond the Julia paper.}
SOC kernels are dispatched in four dimension tiers: one thread per
cone for $d\le 8$, one warp for $8<d\le 64$, one block for
$64<d\le 1024$, and for $d>1024$ a multi-block cooperative-grid
kernel whose block count follows the device's SM
count---orthogonal to CuClarabel's ``split into small cones''
advice.
Dense SOCs apply $\Hs$ in the implicit rank-one form
$\eta^2(2ww^\top-J)$ and store it as a triangular block; sparse
SOCs ($d>4$) carry only the diagonal $\eta^2\diag(d,1,\ldots,1)$ in
the KKT, while the rank-two part $uu^\top-vv^\top$ is written as two
expansion columns by the same fused kernel that scatters $\Hs$.
Generalized power cones use the same $8/64/1024$ cutoffs, but their
large tier is a multi-block kernel without cooperative launch; a
dense Hessian is used up to $d\le 32$ and a sparse rank-three
expansion $\mu(D+pp^\top-qq^\top-rr^\top)$ above it (except the 3D
primal-dual slice, which stays dense).
PSD cones follow a different model: they are grouped by matrix
order, blocks of order $n\le 32$ use batched cuSOLVER Jacobi
eigensolvers, larger cones use a per-cone divide-and-conquer
eigensolver, and independent dimension groups may overlap on up to
four worker streams.
If chordal decomposition is enabled, a host pipeline builds the
sparsity graph, supernodes, and a merge tree.
Repeated merges of the smallest clique into its parent produce equal
block sizes for the batched eigensolver; if that fails, chordal is
abandoned rather than launching a mixed-dimension batch incorrectly.

Optional CUDA graphs capture the $\tauvar$-independent residual
prefix; they are off by default, because cuDSS stays outside the
graph and capture has a setup cost that is recovered only on long
DPP loops.

A portfolio risk constraint of dimension $\sim 20{,}000$ falls in
the cooperative multi-block tier above.

\subsection{Techniques versus CuClarabel}
\label{sec:gpu-cudss}

The factorization dominates wall time on portfolio-sized KKTs.
Kernel fusion outside cuDSS is therefore a second-order effect, and
we do not claim it as the primary source of speedup.
Relative to CuClarabel---whose description is mixed-family cone
launches plus a generic cuDSS call \cite{chen2025cuclarabel}---the
items below are either absent from that paper, listed there as
future work, or present only as a one-line library invocation.

\paragraph{Kernel fusion.}
Every launch costs a few microseconds, and an unfused equilibration
or residual pass would idle the device between short kernels.
Ruiz therefore fuses the inner scaling stages; Hs scatter into the
KKT is fused with the sparse-SOC rank-two expand; the Mehrotra
combined right-hand side writes the cone shift once; residual dots
$[b^\top z,\,s^\top z,\,q^\top x,\,x^\top Px]$ are one reduction.
The aggregate effect is a few percent of an iteration; the value of
fusion is that it removes fixed per-launch costs that would otherwise
grow relative to the solve as $n$ shrinks.

\paragraph{Hybrid execute (default on).}
The root supernodes of the elimination tree are small and dense,
and the host handles them faster than the device can; letting
cuDSS run those nodes on the CPU while the device processes the
remainder overlaps the two.
Together with two right-hand sides this is the default.
Hybrid on/off is numerically equivalent, not bitwise identical;
\Cref{tab:batch} reports both modes, and the ``Sequential RHS''
row of \Cref{tab:ablate} isolates the fused two-right-hand-side
solve ($24.6$\,ms versus $27.2$\,ms on the $n{=}5000$ QP).

\paragraph{A single device memory pool.}
cuDSS allocates internal LU and workspace buffers inside each
execute; served from the raw allocator, those allocations sit on
the critical path of the first solves and can fragment device
memory in a long-lived process.
Registering the application's pool as cuDSS's device allocator
puts solver memory in the same pool as training tensors and makes
the footprint predictable.
The effect is on latency variance rather than mean throughput and
is not isolated in \Cref{sec:experiments}.

\paragraph{Nested-dissection depth.}
The reordering depth is controlled by the cuDSS
\texttt{ND\_NLEVELS} parameter (default $10$): the minimum number of
levels the nested-dissection partition tree must have
\cite{cudss2024}.
The parameter trades the cost of matrix analysis against the
parallelism exposed in the factorization.
Setting it too large inflates the analysis (symbolic) phase and its
partition-tree bookkeeping; setting it too small leaves fewer
independent subtrees at the top of the tree, lowering the parallelism
the numeric factorization can exploit.
For our QP KKT systems the library default $10$ is essentially
sufficient, but the factor-model SOCP produces one large unsplit SOC
whose dense root supernodes defeat a shallow dissection: deepening
the tree to $12$ keeps those root supernodes small enough for the
hybrid path.
The protocol of \Cref{sec:experiments} therefore uses $10$ on QPs and
$12$ on SOCPs, and does not re-sweep the depth; it is a
task-dependent knob rather than a universal optimum.

\paragraph{Latched AUTO iterative refinement.}
Refinement is only useful near convergence: while the iterate is far
from the solution the factorization error is not the limiting term,
and every refined solve incurs extra triangular passes.
Toggling it on and off between iterations, however, changes the step
and can push the convergence state machine into insufficient-progress
restorations.
Both backends therefore gate refinement with the \emph{same}
three-state Off--Watch--On machine, evaluated once per IPM iteration
from the solver info: the homogeneous ratio $\kappa/\tau$, the
absolute and relative duality gap, the primal and dual residuals, the
previous iteration's residuals, and the last step length.
Writing $\mathrm{tol}_{\mathrm{feas}}$ and $\mathrm{tol}_{\mathrm{gap}}$
for the feasibility and gap tolerances and $\mathrm{tol}^{\mathrm{red}}_{\mathrm{gap}}$
for the looser reduced gap, the transitions are:

\begin{itemize}[leftmargin=1.2em,itemsep=0.2em]
\item \textbf{Off $\to$ Watch} when $\kappa/\tau\le 1$, both residuals
fall below $10^{6}\,\mathrm{tol}_{\mathrm{feas}}$, and the gap satisfies
$\mathrm{gap}<10^{4}\,\mathrm{tol}_{\mathrm{gap}}$ or already meets
$\mathrm{tol}^{\mathrm{red}}_{\mathrm{gap}}$.
Watch applies no refinement; it only arms the latch.
\item \textbf{Watch $\to$ On} (latched) under any of three rules:
(a)~the high-accuracy band is reached---$\kappa/\tau\le 1$, the gap
below $10^{3}\,\mathrm{tol}_{\mathrm{gap}}$, and both residuals below
$10^{4}\,\mathrm{tol}_{\mathrm{feas}}$; (b)~one further iteration is
spent inside the watch band; or (c)~progress stalls near the reduced
gap---a primal or dual residual regresses by more than $2\times$ from
an already in-band value, or the last step falls below $5\times$ the
minimum switching step length.
\end{itemize}

Once On, refinement is \emph{latched} for the rest of the solve and is
never re-tested. Independently, a single insufficient-progress
termination restores the previous iterate, forces the latch once per
solve, and retries with the full budget.
On every iteration the machine returns either $0$ (Off/Watch) or the
configured step budget; the affine (predictor) KKT solve then receives
half that budget rounded up and the combined (corrector) solve the full
budget.
The $5\times$ step threshold in rule~(c) differs by backend, because
the two use different minimum switching step lengths: it trips at
$\alpha<0.5$ on \cpu{} and $\alpha<0.05$ on \gpu{}.

The gate is shared; the AUTO modes differ only in \emph{which}
refinement consumes the budget.
On \cpu{}, AUTO is a host residual loop against the QDLDL factors.
On \gpu{}, the same machine drives either that host residual loop
(AUTO1) or cuDSS's built-in refinement (AUTO2): AUTO1 applies the
solver's own iteration but drops the two-right-hand-side batch to
sequential columns, whereas AUTO2 keeps the batched
two-right-hand-side solve by handing the step count to cuDSS's
internal refinement.
On the protocol of \Cref{sec:experiments} refinement is off on QPs
and AUTO on SOCPs (\cpu{} AUTO, \gpu{} AUTO2); the SOCP rows of
\Cref{tab:portfolio-fwd} terminate as solved at $10^{-8}$ under that
policy.

\paragraph{Unsplit large SOCs.}
CuClarabel's SOC advice is to split a large cone into small cones
and launch one thread per cone.
Splitting is a modelling intervention: it changes the cone product,
enlarges the KKT, and must be undone before the solution is
reported.
\gpu{} keeps the cone unsplit and matches kernel granularity to cone
size (thread / warp / block / cooperative multi-block at the
$d{=}8/64/1024$ cutoffs of \Cref{sec:gpu-setup}).
Devices without cooperative launch fall back automatically.
The factor-model SOCP of \Cref{sec:portfolio-layer} therefore reaches
the solver as one cone of dimension $k+n+1$ as written; \gpu{} leads
the SOCP rows of \Cref{tab:portfolio-fwd} from $n{=}5000$ onward on
instances whose KKT is dominated by that cone.

\paragraph{Workspace preallocation.}
Dense GEMM workspaces and reduction temporaries are sized once at
construction and reused across IPM iterations; the motivation is
the same as for the memory pool---allocation belongs at
construction, not inside the loop---and the effect is on latency
stability rather than mean time.

\paragraph{Device-side warm-start smoothing.}
After a solved status, the equilibrated $(x,s,z)$ is cached before
unscaling.
The next solve applies Algorithm~1 of \cite{chen2025warmstart} on
device for every cone family (\Cref{sec:warmstart}), so no host
round-trip is inserted between the previous solution and the next
warm-started solve.
This is not CVXPY's instance-cache flag.
The measured effect on the $n{=}5000$ QP is a reduction of the
steady $q$-update solve from $26.1$\,ms to $22.7$\,ms
(\Cref{tab:portfolio-dpp}).
QP re-solves drop iterations on both backends; GPU sequential SOCP
at relative step $10^{-2}$ does not
(\Cref{tab:warmstart,tab:warmstart-seq}).

\paragraph{Asynchronous parallel primitives.}
Default Thrust \cite{cccl2025} policies synchronize the stream after
each algorithm,
which serializes sequences of small transforms that have no host
data dependence.
The solver stream uses a no-sync policy, with a single synchronize
deferred to residual copies, cuDSS boundaries, and the Python
hand-off.
The change is invisible to users and its effect is folded into the
few-percent fusion figure above.

\paragraph{Two-right-hand-side KKT solve.}
Each iteration solves three systems against the same numeric
factor---the affine (predictor) right-hand side, the combined
(corrector) right-hand side, and the constant system
$K[x_2;z_2]=[-\dataq;\datab]$.
The affine residual changes at every iteration, but the constant
system does not, and the two are independent; packing them into one
execute with a two-column dense right-hand side lets a single
numeric factorization serve both, streaming the stored factors once
for both columns instead of twice.
The combined Mehrotra step then reuses the $(x_2,z_2)$ already
computed in the affine step, so the three systems cost two solve
invocations.
This is the one fusion with an isolated measurement: on the
$n{=}5000$ QP, disabling the batch raised the subsequent-solve time
from $24.6$\,ms to $27.2$\,ms, about $11\%$ of the steady solve
(\Cref{tab:ablate}).

\subsection{Parametric updates}
\label{sec:gpu-dpp}

Updates rewrite $\dataP,\dataq,\dataA,\datab$ in place while applying
the cached Ruiz maps
($c\,D\dataP D$, $E\dataA D$, $c\,D\dataq$, $E\datab$).
Updates of $\dataq$ or $\datab$ skip KKT rebuild.
Updates of $\dataA$ use a precomputed $A\to A^\top$ index map and a
fused scatter into $A$, $A^\top$, and the KKT, replacing a transpose
with an $O(\nnz)$ permutation.
Pattern-changing writes are rejected.
CVXPY instance reuse is keyed on $(n,m,\text{cone layout})$.
Equilibration must be turned off when a single solver object is fed
unrelated problems: the first instance's Ruiz diagonals would
otherwise scale later data incorrectly.
Same-structure DPP loops (only $\dataq$ changing) may keep
equilibration on.

\subsection{Delta versus CuClarabel}
\label{sec:gpu-delta}

\input{tables/tab_gpu_vs_cu}

\Cref{tab:gpu-vs-cu} summarizes the differences.
Device setup removes a host bottleneck from cold and first solves.
Elimination-tree reuse, two right-hand sides, hybrid execute, and
values-only updates reduce a high-frequency portfolio re-solve to
an upload of values, a refactorization, and a solve, with no host
participation between iterations.
Forward cuDSS uses an upper-triangular view and no pivoting, which
is safe because the forward KKT is quasi-definite; the adjoint uses
the full matrix, automatic pivoting, and library iterative
refinement because $\Ft$ is unsymmetric and therefore \emph{not}
quasi-definite, so the $LDL^\top$ pivot-none contract of the
forward factor no longer holds (\Cref{sec:diff}).
Mixed-dimension PSD cones, generalized power cones, and unsplit
large SOCs broaden the accepted models relative to the Julia paper,
so models reach the solver as written rather than after a
user-side rewriting pass.
Omitting mixed precision trades some peak performance on dense-$A$
QPs \cite[Table~6]{chen2025cuclarabel} for robustness at $10^{-8}$
gap and feasibility tolerances; \Cref{tab:portfolio-acc} shows the
payoff, with \gpu{} agreeing with Clarabel.rs to within
$2.2\times10^{-7}$ in relative objective across the sweep.
This comparison with CuClarabel is feature-level;
\Cref{sec:exp-maros} times it on an 18-QP subset of this machine,
and we claim no speedup over it.
Differentiation and multi-problem batching are absent from
CuClarabel; they occupy \Cref{sec:diff}.

%% file: figs/fig_gpu_ipm.tex
\begin{tikzpicture}[
  font=\footnotesize,
  box/.style={draw, rounded corners=2pt, align=center, inner sep=4pt},
  ipmstep/.style={box, fill=blue!8, minimum width=2.15cm, minimum height=0.88cm},
  cone/.style={box, fill=green!10, minimum width=1.85cm, minimum height=0.70cm},
  kkt/.style={box, fill=orange!12, minimum width=2.25cm, minimum height=0.88cm},
  arr/.style={-{Latex}, thick}
]
  \node[ipmstep] (res) {Residuals\\+ Info};
  \node[ipmstep, right=0.48cm of res] (scale) {Cone scaling};
  \node[kkt, right=0.48cm of scale] (upd) {Hs scatter\\+ regularize};
  \node[kkt, right=0.48cm of upd] (fac) {cuDSS\\REFACTOR};

  \node[cone, below=1.20cm of scale, xshift=-1.10cm] (soc) {SOC tiers};
  \node[cone, right=0.30cm of soc] (nn) {NN stream};
  \node[cone, below=0.32cm of soc] (nonsym) {Exp/Pow};
  \node[cone, right=0.30cm of nonsym] (psd) {PSD / chordal};
  \node[draw=green!55!black, dashed, rounded corners=3pt, inner sep=5pt,
        fit=(soc)(nn)(nonsym)(psd)] (conefam) {};

  \node[kkt, below=0.85cm of conefam] (sol) {nrhs$=2$ solve};
  \node[ipmstep, below=0.62cm of sol] (stepn) {Step length\\+ update};

  \draw[arr] (res) -- (scale);
  \draw[arr] (scale) -- (upd);
  \draw[arr] (upd) -- (fac);
  \draw[arr] (scale.south) -- (conefam.north);

  \draw[arr] (fac.east) -- ++(0.52,0) |- (sol.east);
  \draw[arr] (sol) -- (stepn);

  \draw[arr, dashed]
    (stepn.west) -- ([xshift=-0.72cm]res.west |- stepn.west)
    -- ([xshift=-0.72cm]res.west) -- (res.west);
\end{tikzpicture}

%% file: tables/tab_gpu_vs_cu.tex
\begin{table}[t]
\centering
\caption{CuClarabel versus \gpu{}. Mixed-family cone parallelism and cuDSS are shared;
device setup, PSD generality, factorization reuse, parametric updates, and
differentiation are not.}
\label{tab:gpu-vs-cu}
\small
\setlength{\tabcolsep}{3.2pt}
\begin{tabular}{@{}llp{5.6cm}@{}}
\toprule
Topic & CuClarabel & \gpu{} \\
\midrule
Language & Julia + CUDA.jl & CUDA C++ \\
Equilibration / KKT fill & host; paper future work & fully on device, fused Ruiz \\
PSD cones & same dimension, batched & mixed dim.\ + chordal \\
SOC parallelism & one thread per cone & four tiers + cooperative groups \\
cuDSS reuse & per solve / basic & etree cache, two RHS, hybrid execute \\
Iterative refinement & $\le 10$ library steps & latched AUTO near tolerance \\
ND depth & library default ($10$) & $12$ (large unsplit SOC) \\
Warm start & none in the paper & device smoothing, all families \\
Mixed precision & optional & binary64 only \\
Parametric updates & limited & values-only path, cached maps \\
Differentiation & none & \diffsolver{} + \batchdiff{} \\
\bottomrule
\end{tabular}
\end{table}

%% file: sections/07_warmstart.tex
\section{Smoothing-operator warm start}
\label{sec:warmstart}

Interior-point methods resist naive warm starts: a previous optimum
lies on the boundary of $\K$, where the barrier Hessian
diverges and the next Newton step collapses \cite{john2008warm,skajaa2014warm}.
First-order solvers have no such restriction
\cite{odonoghue2016scs,stellato2020osqp}.
Chen, Goulart, and Jones \cite{chen2025warmstart} propose a
\emph{smoothing operator} that maps a previous $(s,z)$ to a point on
the central path of the new instance, prioritizing centrality over
feasibility: it perturbs the boundary point just enough to
restore centrality. The user-facing control is a single boolean.

\subsection{Smoothing operator}
\label{sec:smoothing-theory}

Let $f$ be a $\nu$-logarithmically homogeneous self-concordant barrier
for $\K$ and $f^*$ its conjugate barrier on $\Kstar$
\cite{nesterov1997selfscaled}.
Replacing the indicator in the Euclidean projection by $\eta f$ gives
\begin{equation}
\label{eq:smoothing}
S_{\K,\eta}(c)
\;=\;
\arg\min_{s}\;
\tfrac12\|s-c\|_2^2 + \eta f(s).
\end{equation}
This is a proximal map \cite{parikh2014proximal,moreau1962} and
decomposes over the product $\K=\K_1\times\cdots\times\K_p$.
Algorithm~1 of \cite{chen2025warmstart} applies the
\emph{dual} form at $\lambda=1$
\begin{equation}
\label{eq:smoothing-dual}
\min_z\;
\tfrac12\|z-w\|_2^2 + \mubar_0 f^*(z),
\qquad
w=z_{\mathrm{old}}-s_{\mathrm{old}},
\qquad
s_{\mathrm{new}}=z_{\mathrm{new}}-w,
\end{equation}
which is Moreau-equivalent to \eqref{eq:smoothing} at $\lambda=1$ and
places $(s,z)$ on the central path
$z=-\mubar\nabla f(s)$ (\cite[Theorem~3.1]{chen2025warmstart}).
For a second-order cone the paper instead takes
$\lambda=s_0^\star/z_0^\star$, so that $c=s^\star-\lambda z^\star$ has
vanishing $0$-component; the proximal map then lifts only the Lorentz
scalar (\cite[Theorem~4.1]{chen2025warmstart}).
For nonnegative, second-order, and PSD cones the additional residual
introduced by smoothing is on the order of $\|\dataA\|\mubar_0$.
Closed forms exist for the zero, nonnegative, and second-order cones;
three-dimensional dual cones use a damped Newton iteration with
$\lambda^\star=2-\sqrt{3}$ (\cite[Theorem~5.1]{chen2025warmstart}).

\begin{theorem}[{\cite[Theorems~3.1 and~4.1]{chen2025warmstart}}]
\label{thm:cgj}
Let $f$ be a $\nu$-logarithmically homogeneous self-concordant
barrier for $\K$.
The dual proximal map \eqref{eq:smoothing-dual} at $\lambda=1$
produces a point on the central path $z=-\mubar\nabla f(s)$.
For a second-order cone with $\lambda=s_0^\star/z_0^\star$, writing
$R$ for the KKT residual in the coordinates of the operator,
\begin{equation}
\label{eq:socbound}
\|R(x^{0},s^{0},z^{0})\|_{\infty}
\;\le\;
\Big[\,1+\frac{\|A\|_{\infty}+\lambda}{2\,s_{0}^{\star}}\Big]\,
\|R(x^{\star},s^{\star},z^{\star})\|_{\infty},
\qquad \lambda=\tfrac{s_{0}^{\star}}{z_{0}^{\star}}.
\end{equation}
\end{theorem}
\begin{remark}
\label{rem:cgj-scope}
The bound is stated in a fixed coordinate system.
It does not address a preceding Ruiz equilibration, nor a
subsequent truncation of $\mubar_{\mathrm{soc}}$.
\end{remark}

\subsection{Implementation}
\label{sec:smoothing-impl}

Both solvers apply the smoother in \emph{equilibrated} coordinates,
because the operator is only meaningful on the well-scaled central
path that the IPM actually traverses.
The Clarabel.jl prototype instead unscales, smooths in user space,
and re-equilibrates; that round-trip does not commute with Ruiz
diagonal scaling, so a point that is central in user coordinates
need not be central in the coordinates of the next Newton step.

After a solved or almost-solved termination, the
equilibrated iterate is cached \emph{before} post-process unscaling,
then divided by the homogeneous $\tau^\star$ so that the seed lies
on the $\tau=1$ slice of Algorithm~1
(\cite[\S6.1]{chen2025warmstart}).
The embedding typically terminates at $\tau^\star\neq 1$; restoring
$\tau=1$ without that division makes $Ax+s-b$ of order
$|1-1/\tau^\star|$ and forces $\mubar_0\sim 1$, which is an
oversmoothed cold start rather than a warm one.
On the next solve, if smoothing is enabled
and the cache is valid, the solver restores that triple, recomputes
$\mubar_0$ from the new instance's feasibility residuals, and calls
the per-cone smoother; otherwise it uses Clarabel's default start.
Failed solves invalidate the cache.
The duality gap is not folded into $\mubar_0$.
The Julia prototype uses
$\max(\mathrm{pres},\mathrm{dres},\min(\mathrm{gap}_{\mathrm{abs}},\mathrm{gap}_{\mathrm{rel}}),10^{-6})$;
after cost equilibration, $|\Delta q^\top x|/c$ can be $O(1)$ for a
tiny $q$-jump and that extra term again forces $\mubar_0\sim 1$.
Feasibility residuals track the perturbation size.
Homogeneous scalars after smoothing follow \cite[\S6.1]{chen2025warmstart}:
$\tau=1$, $\kappa=\mubar_0$.
The user-facing control is a single boolean; the smoothing
parameter is not exposed.
The control is distinct from a modelling-language instance cache,
which only skips constructing a new solver object.

\begin{definition}[Implemented smoothing map]
\label{def:smoothing-impl}
Given a cached complementary pair $(s^\star,z^\star)$ on the
equilibrated $\tau=1$ slice and a new instance, set
\begin{equation}
\label{eq:mu0}
\mubar_0
=
\max\bigl(\mathrm{res}_{\mathrm{primal}},\,
\mathrm{res}_{\mathrm{dual}},\,10^{-6}\bigr),
\end{equation}
the homogeneous feasibility residuals of the restored point on the
new data (primal and dual only; the duality gap is excluded).
The per-cone map is that of \cite{chen2025warmstart} except as
noted:
the zero cone uses the Moreau rule $z\leftarrow z-s$, $s\leftarrow 0$;
a second-order cone uses
\cite[Theorem~4.1]{chen2025warmstart} with
$\lambda=s_0^\star/z_0^\star$ and
$\mubar_{\mathrm{soc}}=\lambda\mubar_0$, subject to the relative
floor below and an optional relative cap
(production default: off on \cpu{}, on on \gpu{});
nonnegative, PSD, and three-dimensional cones use the dual proximal
form at $\lambda=1$;
exponential, power, and generalized-power cones use a damped Newton
iteration with $\lambda^\star=2-\sqrt{3}$.
\end{definition}
The GPU SOC cap is an algorithmic variant: it defines a different
operator from the CPU map (floor only) and from the closed form
$\mubar_{\mathrm{soc}}=\lambda\mubar_0$ of
\cite[Theorem~4.1]{chen2025warmstart}.
The iteration bound of that theorem does not cover the capped map.

Per-cone maps, matching \cite{chen2025warmstart} except where noted:

\begin{itemize}[leftmargin=1.2em,itemsep=0.2em]
\item \textbf{Zero cone.}
      Clarabel.jl Moreau: $z\leftarrow z-s$, $s\leftarrow 0$, rather
      than keeping $(s^\star,z^\star)$.
\item \textbf{SOC.}
      Theorem~4.1 of \cite{chen2025warmstart}:
      $\lambda=s_0^\star/z_0^\star$ (so $c_0=0$), smoother parameter
      $\mubar_{\mathrm{soc}}=\lambda\mubar_0$, closed form that
      preserves the Lorentz generator $s_1=c_1/2$.
      Dual Appendix~B at $\lambda=1$---the Julia prototype for every
      SOC---is the fallback when $s_0^\star$ or $z_0^\star$ is
      non-positive.
      A relative floor
      $\mubar_{\mathrm{soc}}\ge 10^{-4}\,\|c_1\|_2^2/4$ keeps the
      Lorentz residual representable on a large cone.
      On the GPU we also \emph{cap} $\mubar_{\mathrm{soc}}$ at the
      same relative gap against $\max(\|c_1\|_2^2/4,s_0^2)$:
      problem $\mubar_0$ is often a linear dual residual from
      $\Delta q$, and feeding $\lambda\mubar_0$ into a complementary
      SOC of dimension $\sim 10^3$ can lift $s_0$ by $\sim 10\%$ and
      increase the Mehrotra iteration count in the sparse four-tier
      SOC scaling.
      The cap is an \emph{algorithmic variant} of Theorem~4.1, not
      an implementation detail: the GPU map is a different operator
      from the CPU map (floor only) and from the paper's closed
      form $\mubar_{\mathrm{soc}}=\lambda\mubar_0$.
      The theorem's iteration bound does not cover the capped GPU
      map.
      Isolating the cap on the CPU (\Cref{tab:soc-cap}) leaves
      sequential $R=0.85<1$; the production floor-only CPU map
      has $R=0.77$.
      The cap therefore trims the CPU saving but does
      \emph{not} reproduce the GPU sequential failure
      ($R=1.20$).
      The GPU floor cell of \Cref{tab:soc-cap} is worse still
      (sequential $R=1.34$), so the cap is neither sufficient nor
      necessary for $R>1$; the failure is a GPU-path effect of
      Theorem~4.1's lift of $s_0$ with four-tier sparse SOC
      scaling.
\item \textbf{Nonnegative.}
      Closed form of \cite[eq.~(32)]{chen2025warmstart}, $\lambda=1$.
\item \textbf{Exp / power / GenPower.}
      Damped Newton on the dual Hessian, at most one hundred
      iterations, with $\lambda^\star=2-\sqrt{3}$.
\item \textbf{PSD.}
      Eigen-decomposition (host LAPACK on CPU; batched Jacobi on
      GPU), smooth the eigenvalues, reconstruct, $\lambda=1$.
      The Julia prototype does not implement this cone.
\end{itemize}

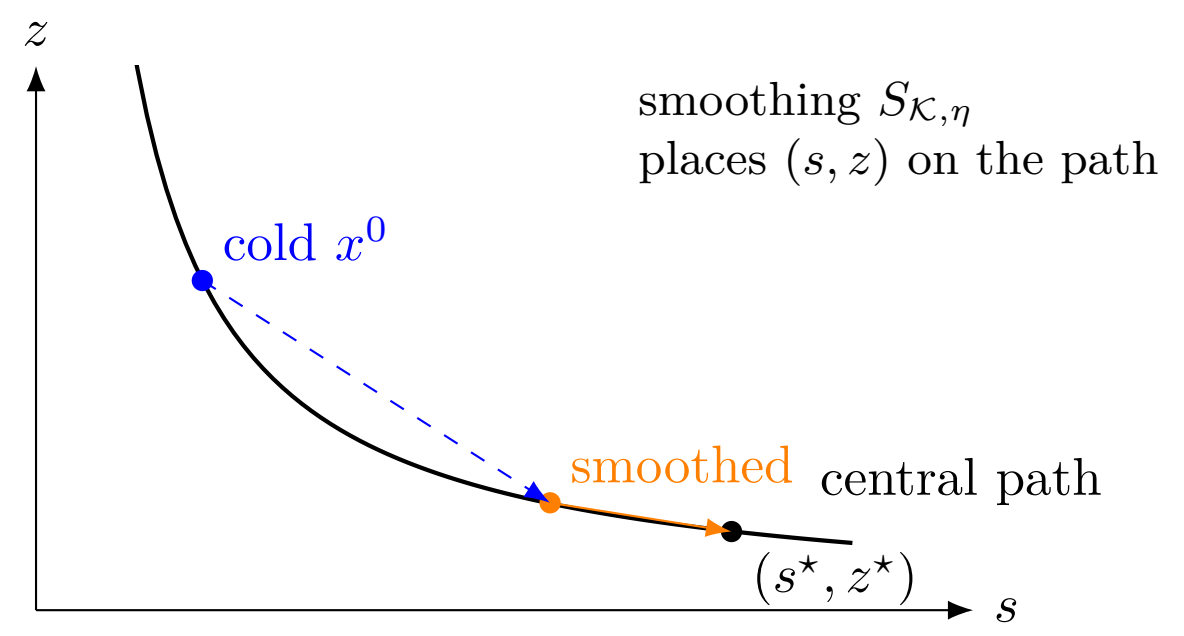
\begin{figure}[t]
\centering
\input{figs/fig_warmstart}
\caption{Cold start versus smoothed warm start.
The previous optimum lies near the boundary.
Smoothing $S_{\K,\eta}$ returns a point on the central path from which
the next Mehrotra iteration can take a long step.}
\label{fig:warmstart}
\end{figure}

\Cref{fig:warmstart} sketches the geometry.
On QPs (zero and nonnegative cones) the prediction of
\cite{chen2025warmstart} is sharp: nearby re-solves drop iterations,
and the gain shrinks as the jump grows.
On a large complementary SOC the same prediction holds only while
$\mubar_0$ is small enough that Theorem~4.1 does not over-lift
$s_0$.
\Cref{sec:exp-warm} tests both regimes on the factor-model family,
as a one-shot perturbation of a cached solve and as a sequential
rebalancing walk in $\mu$.

\subsection{A parameter-agnostic problem distance}
\label{sec:dkkt-def}

The one-shot and sequential tests of \Cref{sec:exp-warm} perturb a
single coordinate of the data.
In applications several inputs may change at once.
The following residual aggregates every perturbation into one
scalar, independently of which parameter moved.

\begin{definition}[Problem distance]
\label{def:dkkt}
Let $(x^\star,s^\star,z^\star)$ be a cached solution of a previous
instance and $(P',q',A',b')$ the data of a new instance.
The problem distance is the normalised residual of that old point
in the new KKT system,
\begin{equation}
\label{eq:dkkt}
d_{\mathrm{KKT}}
=\frac{\|P'x^\star+q'+A'^{\top}z^\star\|_2}{1+\|q'\|_2}
+\frac{\|A'x^\star+s^\star-b'\|_2}{1+\|b'\|_2}.
\end{equation}
\end{definition}
This is the size of the first corrector an interior-point restart
must take, so it predicts the warm-start cost directly.

\begin{remark}[Equilibration invariance]
\label{prop:dkkt-invar}
Equation \eqref{eq:dkkt} is written in the unequilibrated data
$(P',q',A',b')$ and the unequilibrated cached point
$(x^\star,s^\star,z^\star)$.
No Ruiz diagonal enters the formula, so the scalar is invariant
under a subsequent equilibration of the IPM.
It is a definitional invariance, not an independent theorem.
\end{remark}
As a geometric complement, \Cref{sec:exp-warm} also records the
fraction $d_{\mathrm{act}}$ of holdings whose active bound flips.

%% file: figs/fig_warmstart.tex
\begin{tikzpicture}[
  font=\small,
  >=Latex
]
  \draw[->] (0,0) -- (6.2,0) node[right] {$s$};
  \draw[->] (0,0) -- (0,3.6) node[above] {$z$};
  \begin{scope}
    \clip (0,0) rectangle (6.2,3.6);
    \draw[thick, domain=0.35:5.4, samples=80] plot (\x, {2.4/\x});
  \end{scope}
  \node[anchor=west] at (5.05,0.85) {central path};
  \fill (4.6,0.52) circle (2pt) node[below right] {$(s^\star,z^\star)$};
  \fill[blue] (1.1,2.18) circle (2pt) node[above right] {cold $x^0$};
  \fill[orange] (3.4,0.71) circle (2pt) node[above right] {smoothed};
  \draw[blue, dashed, ->] (1.1,2.18) -- (3.4,0.71);
  \draw[orange, ->] (3.4,0.71) -- (4.6,0.52);
  \node[align=left, font=\footnotesize, anchor=west] at (3.85,3.15)
    {smoothing $S_{\K,\eta}$\\places $(s,z)$ on the path};
\end{tikzpicture}

%% file: sections/08_diff.tex
\section{Differentiable layer and batch training}
\label{sec:diff}

This section describes the contribution that none of the four
reference systems contains as an integrated design: a reduced
two-block VJP fused to this IPM, and, on the GPU, a multi-problem
cuDSS batch.
The two-interface layout this section builds on is summarized in
\Cref{sec:overview}.

\subsection{Implicit differentiation of a QCP layer}
\label{sec:ift}

A differentiable optimization layer makes a solver behave like a
neural-network layer \cite{amos2017optnet,agrawal2019cvxpylayers}:
parameters $\theta$ (here the QCP data $(\dataP,\dataq,\dataA,\datab)$)
are the input, the optimizer $x^\star(\theta)$ is the output, and
backpropagation needs $\partial L/\partial\theta$ given upstream
$\partial L/\partial x^\star$.
Writing the KKT conditions at a locally unique primal-dual solution
as $G(x^\star,\theta)=0$ and differentiating the identity
$G(x^\star(\theta),\theta)\equiv 0$ gives
$\partial x^\star/\partial\theta
= -(\partial G/\partial x)^{-1}\,\partial G/\partial\theta$
\cite{amos2017optnet,agrawal2019differentiating,healey2025diffqcp}.
Backpropagation never forms that Jacobian: one linear solve
\begin{equation}
\label{eq:vjp}
\Bigl(\frac{\partial G}{\partial x}\Bigr)^{\!\top}\lambda
=
\frac{\partial L}{\partial x^\star}
\end{equation}
followed by the closed-form
$\partial L/\partial\theta=-(\partial G/\partial\theta)^\top\lambda$
constitutes the backward pass: one factorization-class solve
instead of an unrolled IPM.
Two standard realizations of $G$ exist: the \emph{direct} KKT map,
which requires $\theta$ to keep the problem feasible, and the
\emph{homogeneous embedding} of
\cite{agrawal2019differentiating,healey2025diffqcp}, which also
covers infeasible $\theta$ through Farkas certificates.
We use the former (\Cref{sec:two-block}).
The implicit-function theorem requires $\partial G/\partial x$
invertible, which fails at degenerate points; on a boundary SOC we
project out the singular direction before the adjoint solve
(\Cref{sec:dpi}).
The cost of a backward pass is independent of the forward
iteration count: one linear solve of the same size class, whether
the IPM took nine iterations or ninety.

\subsection{Two-block adjoint versus the homogeneous Minty map}
\label{sec:two-block}

\texttt{diffqcp} \cite{healey2025diffqcp} treats the QCP solution map
$S:\theta\mapsto(x,y,s)$ as $\phi\circ s$, where $s$ sends problem
data to a complementary solution of the homogeneous embedding
\eqref{eq:embedding} and $\phi$ recovers $(x,y,s)$.
That is the homogeneous realization of $G$ in \Cref{sec:ift}.
Differentiability follows from the implicit-function theorem applied
to Minty's parameterization $M(z)=(\proj_K z,-\proj_{K^\circ}z)$ of
the complementarity set, as in Agrawal et al.\
\cite{agrawal2019differentiating}.
The resulting operator is solver-agnostic: any QCP solver may supply
$s(\theta)$.

\diffsolver{} instead takes the \emph{direct} KKT route: it linearizes
\eqref{eq:kkt} at a Clarabel solution $(x^\star,z^\star,s^\star)$
already computed by \emph{this} IPM.
Writing complementarity in Minty form $z=\proj_{\Kstar}(z-s)$, the
conditions \eqref{eq:kkt} read $G(x,z,s;\theta)=0$ with
\begin{equation*}
G(x,z,s;\theta)
=
\begin{bmatrix}
\dataP x + \dataA\transp z + \dataq\\
\dataA x + s - \datab\\
z - \proj_{\Kstar}(z-s)
\end{bmatrix}.
\end{equation*}
Let $v=z^\star-s^\star$ and let $\Dpi=\Dpi_{\Kstar}(v)$ denote the
derivative of the projection onto $\Kstar$ at $v$ (block diagonal
over cones).
Linearizing $G$ at the computed solution, write
$u=dz-ds$ for the differential of the Minty parameter $v=z-s$
(so $dz=\Dpi u$ and $ds=-(I-\Dpi)u$, using that $\Dpi$ is a
Euclidean-projection Jacobian and hence symmetric).
Eliminating $ds$ through the linearized primal feasibility block
then yields a reduced Jacobian
$F\in\R^{(n+m)\times(n+m)}$
in the coordinates $(dx,u)$, not
$(dx,dz)$:
\begin{equation}
\label{eq:F}
F
=
\begin{bmatrix}
\dataP & \dataA\transp\Dpi \\
-\dataA & I-\Dpi
\end{bmatrix}.
\end{equation}
The adjoint system \eqref{eq:vjp} is therefore
\begin{equation}
\label{eq:Ft}
\Ft w
=
\begin{bmatrix}
\dataP & -\dataA\transp \\
\Dpi\transp \dataA & I-\Dpi\transp
\end{bmatrix}
\begin{bmatrix} w_1 \\ w_2 \end{bmatrix}
=
-g .
\end{equation}
The $(1,2)$ block of $\Ft$ is $F_{21}^\top=(-A)^\top=-A^\top$,
not the transpose of $F_{12}=A^\top\Dpi$: $F$ is unsymmetric, so
those two off-diagonal blocks are unrelated.
The linearization, the sign of the second row, the assembly of
$g$, and the extraction of $(w_1,w_2)$ are written out in
\Cref{app:reduced-jacobian}.
We use this direct $(n{+}m)$ Jacobian rather than the
$(n{+}m{+}1)$ homogeneous Jacobian of \eqref{eq:embedding}: at a
computed solution with $\tauvar^\star>0$ the embedding carries an
extra nearly singular direction in $\tauvar$, which is exactly the
conditioning the adjoint solve must avoid.
The right-hand side $g$ assembles the upstream gradients: with
$g_x=\partial L/\partial x$ and optional
$g_z=\partial L/\partial z$, $g_s=\partial L/\partial s$,
\begin{equation}
\label{eq:g}
g
=
\begin{bmatrix}
g_x\\
\Dpi\,(g_z+g_s)-g_s
\end{bmatrix},
\end{equation}
where the second block is zero when only $g_x$ is given, and the
singular direction of any boundary SOC is projected out of $g$
before the solve (\Cref{sec:dpi}).
Solving \eqref{eq:Ft} once and extracting
\begin{equation}
\label{eq:grads}
\frac{\partial L}{\partial \dataq}=w_1,
\quad
\frac{\partial L}{\partial \datab}=w_2,
\quad
\frac{\partial L}{\partial \dataP}=\tfrac12\bigl(w_1 x\transp + x w_1\transp\bigr),
\quad
\frac{\partial L}{\partial \dataA}
=
\proj_{\Kstar}(z-s)\,w_1\transp - w_2 x\transp,
\end{equation}
with $\partial L/\partial \dataP$ and $\partial L/\partial \dataA$
returned on the sparsity pattern of $\dataP$ and $\dataA$.

\begin{proposition}[Reduced two-block VJP]
\label{prop:vjp}
Let $(x^\star,z^\star,s^\star)$ be a locally unique solution of
\eqref{eq:kkt} at data $\theta=(\dataP,\dataq,\dataA,\datab)$, and
assume that $\Dpi=\Dpi_{\Kstar}(z^\star-s^\star)$ exists and that
the reduced Jacobian $F$ of \eqref{eq:F} is invertible.
Let $L$ be continuously differentiable in a neighbourhood of
$(x^\star,z^\star,s^\star)$.
Then there exists a unique $w=(w_1,w_2)$ solving
$\Ft w=-g$ with $g$ as in \eqref{eq:g}, and the derivatives of $L$
with respect to $\theta$ are given by \eqref{eq:grads}.
The hypothesis that $\Dpi$ exists excludes a boundary SOC, where
the Euclidean projection is not differentiable; the
implementation then returns one Clarke element after projecting
the face direction out of $g$ (\Cref{sec:dpi}), which is outside
the scope of this statement.
\end{proposition}
\begin{proof}
See \Cref{app:reduced-jacobian}.
\end{proof}
Two implementation facts follow.

\paragraph{$F$ is not the IPM KKT.}
The last factored $\Hs$-system \eqref{eq:kkt-matrix} is a
regularized, scaled Newton matrix on the embedding.
$\Ft$ is an unsymmetric Jacobian of the projection form of
\eqref{eq:kkt}.
Forward $LDL^\top$ factors cannot be backsolved for the VJP
\cite{goulart2024clarabel,healey2025diffqcp}.
\cpu{} therefore factors $\Ft$ with KLU by default
\cite{davis2011klu}; an optional QDLDL path lifts
$Q=\bigl[\begin{smallmatrix}\varepsilon I & F^{\top}\\F & -\varepsilon I\end{smallmatrix}\bigr]$.
If a PSD cone is present or $n+m\le 800$, a dense LU is used instead.
The KLU (or QDLDL) pattern is analyzed once; later
updates of $q$ or $b$ only scatter values.

\paragraph{Original-space gradients.}
When Ruiz equilibration is on, the IPM works with scaled
$(\tilde P,\tilde q,\tilde A,\tilde b)$.
\gpu{} keeps device copies of the original $\dataP,\dataA,\dataq,\datab$
and of $A^\top$, and applies $\Dpi$ and \eqref{eq:grads} in that
space, so the VJP is with respect to the user's parameters.
Device pointers on the solver stream update $q$ and $b$ without a
host round-trip.

\subsection{Cone projection derivatives}
\label{sec:dpi}

$\Dpi_{\Kstar}(v)$ is block diagonal.
Zero and nonnegative cones contribute diagonal $0/1$ masks.
A second-order cone is the identity (interior of the dual), zero
(interior of the polar), or a rank-at-most-two operator on the
boundary.
If $d>16$ that operator is stored as the Woodbury form
$\Dpi=D_d+V D_a^{-1}V^\top$ with $V\in\R^{d\times 3}$,
$D_a=\mathrm{diag}(+1,-1,-\mathrm{sign}\,\tau)$
(\Cref{app:soc-dpi}): three signed outer products whose span is
two-dimensional, because the third generator is parallel to
$p-q$.
Three columns, not two, keep $D_a$ diagonal and the expanded
KKT pattern independent of $\tau$; $\Ft$ is expanded by those
columns rather than densified.
This adjoint expansion is distinct from the rank-two forward
expansion of \Cref{sec:ipm} (a factorization of $\Dpi$, not of
$\Hs$).
Correctness of the adjoint expansion is checked only through the
end-to-end VJP-versus-FD test of \Cref{sec:exp-diff}, not by a
separate residual on $\Dpi$ itself.
A boundary SOC also introduces a singular face direction, which is
projected out of the right-hand side $g$ before the KKT solve and
is \emph{not} a fourth auxiliary column.
Exponential and 3D power cones use analytic $3\times 3$ Jacobians.
A generalized power cone uses Hien's 3D Jacobian
\cite{hien2024power} when $(\cdim_1,\cdim_2)=(2,1)$, and otherwise a
central-difference Jacobian of a Newton projection---the same scheme
on CPU and GPU.
PSD cones use the spectral formula for projection onto $\Splus$
\cite{agrawal2019differentiating,healey2025diffqcp}; a
dense-diagonal kernel reduces the matvec from $O(n^3)$ to a gather
of the $O(n^2)$ svec.

On the GPU, a reduction over the $\Dpi$ diagonal chooses among full
analysis, refactorization, or a cached backsolve, without a host-side
pattern walk.
The common training path (fixed cones, updated $\dataq$) is
refactorization of $\Ft$ after each forward solve, then a cached
backsolve on subsequent backwards with the same projection type.

The projection type is not, however, guaranteed to be stable across
backward passes.
When the forward optimum lands on the boundary of a second-order
cone, $\Dpi$ changes from the identity (interior) or zero (polar) map
to the rank-at-most-two boundary operator of \Cref{sec:dpi}, so the
sparsity pattern of $\Ft$ can change from one backward to the next
even though the QCP dimensions and cone layout are unchanged.
A naive implementation would respond by re-running the full symbolic
analysis on every such transition.
\gpu{} instead separates the \emph{permutation} from the
\emph{pattern}: the fill-reducing permutation and its elimination
tree depend only on the fixed sparsity of $\Ft$---determined by the
cone dimensions and layout---so they are computed once at
construction and cached.
Only when the reduction reports a genuine change in the $\Dpi$
structure (an interior/boundary flip on some SOC block) is the
symbolic factorization refreshed, and even then the cached
elimination tree is re-injected rather than recomputed from scratch.
This makes the matrix-analysis reuse robust precisely in the regime
where a training trajectory visits the SOC boundary repeatedly.

\subsection{\batchdiff{}: batching for differentiable training}
\label{sec:batch}

\begin{figure}[t]
\centering
\input{figs/fig_batch}
\caption{\batchdiff{}.
Each problem owns a stable \diffsolver{} slot so that
backward $i$ differentiates forward $i$ even when
iteration counts diverge.
Workers enqueue cuDSS phases; the scheduler dispatches full
shards as one uniform batch and isolates latency-expired partial
shards.
Forward and backward use separate role channels because symmetry and
pivoting differ (\Cref{sec:gpu-cudss}).}
\label{fig:batch}
\end{figure}
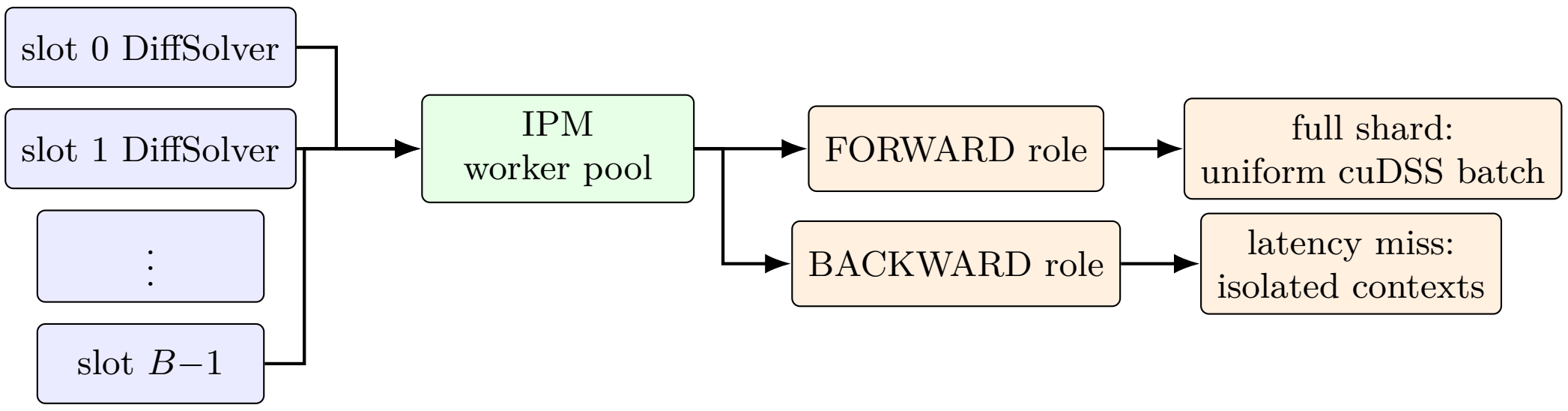

A Python loop over $B$ independent \diffsolver{} objects already
drops the GIL inside each solve and backward
(\Cref{sec:dual-compile}), so an isolated training step is not
serialized on the interpreter: the measured Python-loop baseline in
\Cref{sec:exp-batch} scales almost linearly to $B=16$.
What the loop cannot do is (i)~issue $B$ cuDSS factorizations as
one uniform batch, or (ii)~tolerate divergent IPM iteration
counts without a lock-step kernel stalling on the slowest cone.
\batchdiff{} provides those two
contracts: stable slots, C++ workers, and a cuDSS gather whose
shard width is independent of $B$.
\Cref{sec:exp-batch} times isolated throughput on an idle host,
where a Python loop of warm \diffsolver{} objects remains
2--13\% faster at every $B\le 16$; we do not claim a
GIL-contention speedup, since that experiment is not in the
protocol.
\Cref{fig:batch} shows the data flow.

\paragraph{Structural contract.}
All $B$ items share $n$, $m$, the cone list, and the sparsity
pattern of $\dataP$ and $\dataA$.
Numerical values may differ.
Pattern-changing updates are rejected.

\paragraph{Stable slots and workers.}
Slot $i$ owns a complete \diffsolver{}, including the nested
forward solver, the backward KKT, and the last solution.
A reusable pool of $W$ CPU worker threads ($1\le W\le B$) runs
solve (resp.\ backward) with the GIL released; each task is a full
independent IPM, so iteration counts may diverge across slots.
The public call is synchronous: it returns after the slowest slot;
a max-iteration status does not cancel siblings.

\paragraph{Dynamic cuDSS batching.}
The $B$ slots are partitioned into fixed shards of width $S$, a knob
independent of both $B$ and $W$ ($1\le S\le W$,
$W\ge\lceil B/S\rceil$; default $S{=}1$).
Workers enqueue factor/solve requests tagged by a (phase,
right-hand-side count, refinement-steps) signature, kept FIFO within
a signature; a full signature executes as one uniform cuDSS batch,
while a partial shard waits until it fills or its oldest request
exceeds a latency bound (default $1.5\,\mathrm{ms}$), in which case
that shard migrates to independent device-only contexts.
The scheme is the solver analogue of dynamic batching in
inference-serving systems: requests accumulate per signature until a
full batch or the latency deadline, whichever comes first.
Decoupling $S$ from $B$ is deliberate: an oversized uniform batch
can saturate the device and lose to smaller shards, so the shard
width is exposed separately from the problem count.
All shards of one role share the first successful analysis
permutation and nested-dissection tree.
Hybrid execute is offered for shard sizes $1$ through $8$;
larger shards run device-only.
Hybrid on/off is numerically equivalent, not bitwise identical.

\paragraph{Updates and device I/O.}
A parameter update writes every slot's $\dataq,\datab$; a matrix
update writes new $\dataP,\dataA$ values on the same pattern,
optionally deferring backward-cache invalidation until the next
backward.
A device-side backward accepts $\partial L/\partial x$ on a producer
stream without a host round-trip.

\paragraph{Time-dependent portfolios.}
Independent trajectories should be batched at the \emph{same} time
index.
Consecutive times of one trajectory must not be placed in
parallel: those steps form a sequential warm-start chain.
A fixed trajectory-to-slot assignment keeps smoothing
semantically meaningful.

\paragraph{CPU counterpart.}
\cpu{} has no analogue, for the reason given in
\Cref{sec:cpu-delta}: QDLDL exposes no uniform-batch interface, so a
batched facade would serialize on a single AMD permutation without a
GPU to hide it.
A training loop instead updates one \diffsolver{} or farms
independent problems to processes.

%% file: figs/fig_batch.tex
\begin{tikzpicture}[
  font=\footnotesize,
  box/.style={draw, rounded corners=2pt, align=center, inner sep=4pt},
  slot/.style={box, fill=blue!8, minimum width=2.0cm, minimum height=0.7cm},
  work/.style={box, fill=green!10, minimum width=2.4cm, minimum height=0.85cm},
  sch/.style={box, fill=orange!12, minimum width=2.6cm, minimum height=0.75cm},
  arr/.style={-{Latex}, thick}
]
  \node[slot] (s0) {slot 0 \diffsolver{}};
  \node[slot, below=0.18cm of s0] (s1) {slot 1 \diffsolver{}};
  \node[slot, below=0.18cm of s1] (sd) {$\vdots$};
  \node[slot, below=0.18cm of sd] (sb) {slot $B{-}1$};

  \node[work, right=1.1cm of s1] (w) {IPM\\worker pool};

  \node[sch, right=1.0cm of w] (fwd) {FORWARD role};
  \node[sch, below=0.25cm of fwd] (bwd) {BACKWARD role};
  \node[sch, right=0.7cm of fwd] (full) {full shard:\\uniform cuDSS batch};
  \node[sch, right=0.7cm of bwd] (part) {latency miss:\\isolated contexts};

  \draw[arr] (s0.east) -- ++(0.35,0) |- (w.west);
  \draw[arr] (s1) -- (w);
  \draw[arr] (sb.east) -- ++(0.35,0) |- (w.west);
  \draw[arr] (w) -- (fwd);
  \draw[arr] (w.east) -- ++(0.25,0) |- (bwd.west);
  \draw[arr] (fwd) -- (full);
  \draw[arr] (bwd) -- (part);
\end{tikzpicture}

%% file: sections/09_setup.tex
\section{Experimental setup}
\label{sec:experiments}

This section fixes the measurement protocol; the results themselves
are reported in \Cref{sec:evaluation}, and the Maros--M\'esz\'aros
subset in \Cref{sec:exp-maros}.
Numbers in \Cref{sec:evaluation,sec:exp-maros} come from one
measurement campaign on the
hardware below, except
\Cref{tab:diff,tab:gpu-path}, which were produced by a second
pass after the finite-difference checker was extended (relative error, kink
filter, an SOCP row, and a native-GPU path),
\Cref{tab:cvxlayer-acc,tab:cvxlayer-perf}, which are a later
isolated CVXPYLayers-as-benchmark pass (\Cref{sec:exp-cvxlayer}),
and
\Cref{tab:maros-sub}, which is a later Maros--M\'esz\'aros pass.
The portfolio, DPP, warm-start, batch, and ablation suites were
not rerun.
Each factor-model instance is canonicalized once with CVXPY,
then timed through native solver interfaces (Clarabel.rs, \cpu{}, \gpu{}
\diffsolver{}) and Gurobi through CVXPY.
QP instances are solved by Mosek as quadratic programs;
SOCP instances still use CVXPY's conic Mosek path.
Solver time excludes modelling-stack canonicalization except in the
CVXPYLayers row.
The portfolio, DPP, warm-start, batch, ablation, finite-difference,
and CVXPYLayers layer-baseline suites all use the multi-period
factor-model family of \Cref{sec:portfolio-layer}, with the
single-instance suites taking its $T{=}1$ specialization.
The paper is a performance evaluation of the stack on parametric
portfolio QCPs, plus an 18-problem Maros--M\'esz\'aros QP subset
(\Cref{sec:exp-maros}).
We do not claim corpus-wide rankings; CBLIB \cite{cbf2016}, Netlib LP,
entropy-dominated, and SDP instances remain outside the protocol.

\subsection{Benchmark family: factor-model portfolio QCPs}
\label{sec:portfolio-layer}

We benchmark on a purpose-built factor-model family rather than on a
standard corpus such as Maros--M\'esz\'aros for a structural reason.
A large fraction of the intended financial users model risk through
style factors, so the resulting KKT systems are highly sparse: the
risk enters through a low-rank factor block plus a diagonal, not
through a dense covariance matrix.
That sparsity is exactly the regime the stack is engineered for, and
a generic QP corpus---dominated by denser or less structured
matrices---would not isolate it.
The factor-model family below is therefore the primary workload; an
18-problem Maros--M\'esz\'aros QP subset is reported separately in
\Cref{sec:exp-maros} as a cross-check, not as the main evidence.

A $T$-period factor model with $n$ assets and $k$ factors uses
$\Sigma=F\Omega F^\top+D$, $D=\diag(d)$.
Expected returns are written $r_t$ in this section; the tables keep
the solver-variable name $\mu$ for the corresponding objective block.
Decision variables are weights $x_t\ge 0$, factor exposures
$y_t=F^\top x_t$, and (on the SOCP) an auxiliary risk $t_t$.
Turnover is $c\|x_t-x_{t-1}\|_1$ with $x_{-1}=x_0=0$ and
$c=0.002$.
Risk aversion is $\gamma=1$.
Two canonical QCPs appear throughout the examples; the CVXPY
models below are the instances timed in \Cref{sec:exp-portfolio}.

\paragraph{QP.}
\[
\begin{aligned}
\min\quad&
\sum_{t=0}^{T-1}\Bigl(
-r_t\transp x_t + \gamma\bigl(y_t\transp\Omega y_t + x_t\transp D x_t\bigr)
+ c\|x_t-x_{t-1}\|_1\Bigr)\\
\text{s.t.}\quad&
\mathbf{1}^\top x_0=1,\quad
\mathbf{1}^\top x_t=\mathbf{1}^\top x_{t-1},\quad
0\le x_t\le 0.1,\quad
0\le y_t\le 0.1,\quad
y_t=F^\top x_t .
\end{aligned}
\]
A quadratic form on $y$ plus a diagonal quadratic on $x$
keeps $\dataP$ sparse in the factor metric and DPP-compliant if only
$r$ is a parameter \cite{agrawal2018rewriting}.
CVXPY is asked for a native quadratic objective, not an SOC expansion.

The turnover term $c\|x_t-x_{t-1}\|_1$ is the trading friction: it
prices each one-way trade at the linear transaction tax $c=0.002$
(two per mille of turnover), charged symmetrically on buys and sells,
so the optimizer trades off a rebalance against the expected return
it produces.
Because multi-period rebalancing is a frequent application of the
layer, the model is written as a $T$-period problem in both the
formulation above and the listings below: the state $x_t$ carries
across periods, the turnover of period $t$ is measured against the
previous holdings $x_{t-1}$, and the budget is kept feasible by
$\mathbf{1}^\top x_t=\mathbf{1}^\top x_{t-1}$ rather than by
re-normalizing each period in isolation.
The $T{=}1$ rows of the timing tables are the same model specialized
to a single rebalance.

\begin{lstlisting}[caption={Multi-period factor-model QP (CVXPY).},label={lst:portfolio-qp}]
# Multi-period factor-model QP (native quadratic objective).
# Data: D_ii ~ U[0,1] * sqrt(k), F ~ N(0,1) with 50% zeros,
# Omega = G G^T / k, mu_{i t} = (3 + 9 U[0,1]) / 100.
# Variables are indexed by period t: x_t holdings, y_t factor
# exposures; the state carries across periods through the turnover
# term, which prices every one-way trade at the 0.2% transaction tax.
mu = cp.Parameter((n, T))
x = cp.Variable((n, T), nonneg=True)
y = cp.Variable((k, T), nonneg=True)
cons = [x <= 0.1, y <= 0.1,
        cp.sum(x[:, 0]) == 1]
for t in range(1, T):   # budget rolls over: sum(x_t) == sum(x_{t-1})
    cons += [cp.sum(x[:, t]) == cp.sum(x[:, t - 1])]
for t in range(T):
    cons += [y[:, t] == F.T @ x[:, t]]
terms = []
for t in range(T):
    dx = x[:, t] - (x0 if t == 0 else x[:, t - 1])  # trade vs prior holdings
    terms += [-mu[:, t] @ x[:, t]
              + 0.002 * cp.norm(dx, 1)  # 0.2% tax on buys and sells
              + cp.quad_form(y[:, t], Omega)
              + cp.sum(cp.multiply(D_diag, cp.square(x[:, t])))]
prob = cp.Problem(cp.Minimize(cp.sum(terms)), cons)
\end{lstlisting}

\paragraph{SOCP.}
Replace the quadratic risk by two second-order cones per period:
$\|[L_\Omega y_t;\, D^{1/2}x_t]\|_2\le 0.15$ (hard cap) and
$\|[L_\Omega y_t;\, D^{1/2}x_t]\|_2\le t_t$ (objective).
Factor exposures $y_t$ are free.
The cone product is a zero cone (budget and factor links), a
nonnegative cone (bounds and turnover), and $2T$ SOCs of dimension
$k+n+1$.

\begin{lstlisting}[caption={Multi-period factor-model SOCP (CVXPY).},label={lst:portfolio-socp}]
# Multi-period factor-model SOCP. Risk is an SOC, not a quadratic.
# Hard cap ||risk||_2 <= 0.15 and an auxiliary t in the objective.
# As in the QP, the holdings state carries across periods and each
# one-way trade is charged the 0.2% transaction tax via turnover.
L = np.linalg.cholesky(Omega)
mu = cp.Parameter((n, T))
x = cp.Variable((n, T), nonneg=True)
y = cp.Variable((k, T))          # free factor exposures
t = cp.Variable(T)
cons = [x <= 0.1, cp.sum(x[:, 0]) == 1]
for p in range(1, T):   # budget rolls over: sum(x_p) == sum(x_{p-1})
    cons += [cp.sum(x[:, p]) == cp.sum(x[:, p - 1])]
for p in range(T):
    cons += [y[:, p] == F.T @ x[:, p]]
    risk = cp.hstack([L @ y[:, p], cp.multiply(np.sqrt(D_diag), x[:, p])])
    cons += [cp.norm(risk, 2) <= 0.15, cp.norm(risk, 2) <= t[p]]
terms = []
for p in range(T):
    dx = x[:, p] - (x0 if p == 0 else x[:, p - 1])  # trade vs prior holdings
    terms += [-mu[:, p] @ x[:, p] + 0.002 * cp.norm(dx, 1) + t[p]]  # 0.2% tax
prob = cp.Problem(cp.Minimize(cp.sum(terms)), cons)
\end{lstlisting}

In both cases the standard-form linear cost contains $q=-r$ on the
weight block, so
$\partial L/\partial r = -\partial L/\partial q$ on those
coordinates---the chain rule from \eqref{eq:grads} back to expected
returns, which is \eqref{eq:vjp} with $\theta$ specialized to $r$.
This is the ingredient an end-to-end predict-and-optimize loop
needs: a network that predicts $r$ (or $(r,\Sigma)$) can be
composed with a QP/SOCP layer that returns $x^\star$ and trained
on, e.g., a Sharpe-type loss, without unrolling the IPM.
A solver-based layer, rather than a neural allocation head, is
what keeps the hard constraints of the portfolio---budget, bounds,
factor limits, SOC risk---feasible to IPM tolerance while still
returning $\partial L/\partial r$: a feed-forward head can only
penalize violations, not certify feasibility.
We do not train such a network in this paper; \Cref{sec:evaluation}
measures the layer itself (forward/backward times, gradient checks
against finite differences and against CVXPYLayers), not an
end-to-end learning curve.
Structural tensors $(F,\Omega,D,x_{\max})$ stay constant across a
rebalance or a training step, so values-only updates and
\batchdiff{} hit the cached-pattern path; that is why this
protocol treats portfolio QP/SOCP as the primary
workload.
Every portfolio number in this paper is produced from this single
modelling family: the timing tables of \Cref{sec:evaluation}, the
parametric re-solve, warm-start, differentiable-layer, batch, and
ablation suites, and the wrapper-isolation table all instantiate
exactly the models of
\Cref{lst:portfolio-qp,lst:portfolio-socp}, with the single-instance
suites using their $T{=}1$ specialization.
No separate or simplified portfolio model is timed anywhere in the
paper, so the reader can reproduce any table from these two listings.

\subsection{Reporting rules}
\label{sec:exp-rules}

\paragraph{Hardware and software.}
AMD Ryzen~9 5900X (12 cores / 24 threads; 3.7\,GHz base,
4.8\,GHz boost),
NVIDIA RTX A6000 (SM~8.6, 46\,GiB, driver 580.173.02, CUDA~13.0),
cuDSS~0.8.0, Clarabel.rs~0.11.1, Mosek~11.2.2, Gurobi~13.0.2,
CVXPY~1.9.1, CVXPYLayers~1.2.0, PyTorch~2.12, Python~3.14.4.
Thread budget, set by the harness: Clarabel.rs and \cpu{} run a
single-threaded IPM (QDLDL), with
\texttt{OMP\_NUM\_THREADS=1} so that the OpenMP residual kernels of
\Cref{sec:cpu-opts} do not use additional cores;
the portfolio tables report Mosek and Gurobi at both 1 thread
and 8 threads;
\gpu{} uses the device for the KKT solve, letting cuDSS run root
supernodes on four host threads; CuClarabel is the device IPM on one
Julia host thread.
The earlier Maros--M\'esz\'aros commercial columns used eight
threads only.
All \gpu{} numbers in
\Cref{tab:portfolio-first,tab:portfolio-fwd,tab:portfolio-dpp,tab:warmstart,tab:warmstart-seq,tab:diff,tab:cvxlayer-acc,tab:cvxlayer-perf,tab:ablate,tab:maros-sub}
use the default hybrid AUTO mode (\Cref{sec:gpu-cudss}); hybrid
on/off is numerically equivalent, not bitwise identical, and
\Cref{tab:batch} reports both hybrid AUTO and device-only.

\paragraph{Timing.}
Two quantities are reported, and they are never mixed within one
table.
\emph{First execution} (\Cref{tab:portfolio-first} only) is the
first solver solve on a fresh object, with construction, CVXPY
canonicalization, and data updates excluded from the tabled
milliseconds (cuDSS analysis sits in construction when
\texttt{prefactor\_at\_setup} is on; Gurobi model build and
Python/CUDA lazy init sit outside the timed call).
The plotted value in \Cref{fig:solver-first-subseq} is the
tabled first solver solve of \Cref{tab:portfolio-first}.
\emph{Subsequent mean} is the arithmetic mean of ten sequential
re-solves of the same numerical data on a reused object
(\Cref{tab:portfolio-fwd} and all other steady-state rows).
CVXPY warm start is on for Mosek/Gurobi on the first solve and the
ten re-solves; Clarabel.rs, \cpu{}, and \gpu{} run each re-solve from the
default start with smoothing off.
The batch and ablation suites instead report the median of ten to
twenty steady steps.
Mosek and Gurobi subsequent times are CVXPY solve times.
The timing paths are not identical: \gpu{} rows in the main tables
are timed through the \diffsolver{} object (the training-relevant
path, including a host-side solution copy), \cpu{} rows
through the native forward solver, and Mosek/Gurobi through CVXPY.
At small $n$ these path differences are a visible fraction of the
reported milliseconds, so the small-$n$ ranking should be read as
approximate.
\Cref{tab:gpu-path} isolates the wrapper cost on two QP sizes by
timing the native GPU solver against
\diffsolver{} on the same data.
We do not resample $F,\Omega,\mu$ in the main tables (seed 42).
A three-seed repeat of the isolated-CPU protocol of
\Cref{tab:portfolio-fwd} on the $n{=}5000$ QP
(seeds $42,7,99$) gives subsequent means
$134$/$133$/$143$\,ms (iterations $11$/$11$/$12$):
a coefficient of variation of about $4\%$.
The Table~7 seed-42 cell is $150$\,ms (11 iterations) from an
earlier draw of that same protocol; we treat $150$ as one
isolated draw, not a tighter point estimate than the
three-seed band.
A GPU three-seed on the same instance
($24.7$/$25.8$/$23.8$\,ms, iterations $13$/$14$/$13$)
has a coefficient of variation of about $4\%$.
Those checks cover only this isolated QP; the main portfolio
tables, the SOCP warm-start suite, and the CVXPYLayers comparison
remain single-seed (42).
Percentage-level claims at or below that scatter---including the
6\% CPU $q$-update change of \Cref{tab:portfolio-dpp}---are
indicative, not decisive.
Within one run the mean--median gap is small
(1\% on the $n{=}5000$ QP); a second run of the same instance in
\Cref{tab:ablate} differs by about 7\%.

\begin{table}[t]
\centering
\caption{Table-to-protocol map.
``Ten'' means the mean of ten sequential re-solves on a reused
object; ``twenty'' means the median of twenty steady steps.
The factor family is \Cref{sec:portfolio-layer}
($k{=}50$ for $n\le 5000$; $k{=}100$ at $n{=}10000$).
The four AlmostSolved SOCP cells re-timed at two refinement steps
(\Cref{tab:portfolio-first}) are inherited by
\Cref{tab:portfolio-fwd}.}
\label{tab:protocol-map}
\small
\begin{tabular}{@{}llll@{}}
\toprule
Table & Quantity & Repeats & Notes \\
\midrule
\Cref{tab:portfolio-first} & first solve & one & isolated CPU; other columns original campaign \\
\Cref{tab:portfolio-fwd} & subsequent mean & ten & same splice as \Cref{tab:portfolio-first} \\
\Cref{tab:portfolio-gpu-obj} & ten GPU objectives & ten & default start \\
\Cref{tab:gpu-path} & subsequent mean & twenty & $k{=}50$ at both sizes \\
\Cref{tab:portfolio-dpp} & first / steady $q$-update & one / twenty & $n{=}5000$, $k{=}50$ \\
\Cref{tab:warmstart,tab:warmstart-seq} & iterations & five / twenty steps & median / mean \\
\Cref{tab:warmstart-equil} & iterations & five & median \\
\Cref{tab:diff,tab:eps-sweep} & VJP vs FD & eight coords & interior filter \\
\Cref{tab:cvxlayer-acc,tab:cvxlayer-perf} & acc.\ / time & ten / five & CVX includes canon. \\
\Cref{tab:batch} & median step & ten & $n{=}1000$, $k{=}50$ \\
\Cref{tab:ablate} & median solve & twenty & $n{=}5000$, $k{=}50$ \\
\Cref{tab:maros-sub} & subsequent mean & ten & YimingYAN MAT \\
\Cref{tab:n20000} & subsequent mean & ten & whole table re-timed; CPU isolated; $k{=}50$ \\
\Cref{tab:residuals} & KKT residuals & one & Clarabel form + Mosek native \\
\Cref{tab:vjp-ratios} & VJP/FD ratios & 24 coords & plus $P,A$ \\
\Cref{tab:soc-cap} & SOC cap $2{\times}2$ & five / twenty & CPU/GPU on/off \\
\bottomrule
\end{tabular}
\end{table}

\paragraph{Scope of comparison.}
\texttt{diffqcp} is not timed: it is a solver-agnostic JAX operator
and not a drop-in backend for the factor-model suite, so a
meaningful comparison is a separate study rather than a missing row.
CuClarabel is timed in
\Cref{tab:portfolio-first,tab:portfolio-fwd,tab:maros-sub} on the same machine.
Those times start at the second setup+solve so Julia JIT is not
in the sample (\Cref{sec:exp-maros}).
CuClarabel ran against CUDSS.jl~0.6 with the CUDA math library
\texttt{libcudss}~0.7.1 artifact (CUDA~13); \gpu{} keeps the
conda \texttt{libcudss}~0.8.0 used in the rest of the paper.
That library split is a caveat on the Maros subset, not a
version-matched comparison: a roughly 2\% factorize/solve
advantage of cuDSS~0.8.0 over 0.7.1 is the allowance we make
for the version gap (we do not re-measure it here), and it is
an order of magnitude smaller than the subsequent-time gaps on
that subset (e.g.\ CONT-050 20\,ms versus 58\,ms).
The ranking on those 17 timed QPs is therefore not overturned
by the library version, but it is not version-matched either.
Statements about CuClarabel in
\Cref{sec:related,sec:gpu,sec:diff} remain feature-level, and we
claim no speedup over it.
The timed baselines are Clarabel.rs (same algorithm, reference
Rust implementation), Mosek, and Gurobi.

\paragraph{Tolerances.}
Gap and feasibility tolerances are the official Clarabel
recommended binary64 defaults
\texttt{tol\_gap\_abs}=\texttt{tol\_gap\_rel}=\texttt{tol\_feas}=$10^{-8}$,
with \texttt{max\_iter}=$200$
\cite{clarabelsettings}.
Chordal decomposition is off (official default on).
Iterative refinement uses a budget of two steps against the
official default of ten
\cite{clarabelsettings}:
disabled on QPs for
Clarabel.rs / \cpu{} / \gpu{}; on SOCPs, Clarabel.rs enables
refinement, \cpu{} uses AUTO, and \gpu{} uses AUTO2
(\Cref{sec:gpu-cudss}).
The remaining official knobs that the protocol touches, and the
CPU/GPU-only overrides, are in \Cref{tab:settings}.
Mosek QPs are sent as quadratic objectives; CVXPY
does not expand a quadratic form into a conic program.
Mosek and Gurobi use $10^{-8}$ interior-point / barrier tolerances.
Because the Clarabel.rs / \cpu{} / \gpu{} IR policies differ on SOCPs, SOCP
iteration counts are not directly comparable across solvers
(e.g.\ 29 versus 25 versus 25 at $n{=}20000$ in
\Cref{tab:n20000}); the times remain comparable on that
companion campaign, since refinement is part of the measured
solve.
The $n{=}20000$ size is not in
\Cref{tab:portfolio-first,tab:portfolio-fwd}: those tables stop
at $n{=}10000$, $k{=}100$.

\paragraph{Metrics.}
Wall-clock milliseconds, IPM/barrier iterations, primal objective,
relative objective gap versus Mosek, and
$\|x-x_{\mathrm{Mosek}}\|_2/\|x_{\mathrm{Mosek}}\|_2$ on the
$n$-vector of asset weights (Clarabel canonical $x_{1:n}$ coincides
with this vector on these instances).

%% file: sections/10_evaluation.tex
\section{Evaluation}
\label{sec:evaluation}

All numbers in this section follow the protocol of
\Cref{sec:experiments} on the hardware listed there.

\subsection{Forward: factor-model portfolios}
\label{sec:exp-portfolio}

The timed instances are the multi-period CVXPY models of
\Cref{lst:portfolio-qp,lst:portfolio-socp}
(seed~42; $D_{ii}\sim U[0,1]\sqrt{k}$; $F$ standard normal with a
$50\%$ zero mask; $\Omega=GG^\top/k$;
$\mu_{it}=(3+9U[0,1])/100$).
Sizes are $n{=}200$ with $k{=}50$, $T{=}1$;
$n\in\{500,1000,5000\}$ with $k{=}50$ and $T\in\{1,2,3,4\}$;
$n{=}2000$, $k{=}50$, $T{=}1$; $n{=}10000$, $k{=}100$, $T{=}1$;
and a companion $n{=}20000$, $k{=}50$, $T{=}1$ accuracy row
(\Cref{tab:portfolio-acc}) with isolated timings in
\Cref{tab:n20000}.
The $n{=}100$, $k{=}50$ pair is omitted from the reported sweep
(the corresponding SOCP is primal infeasible under the written
$0.15$ cap).
QP uses a native quadratic objective (not an SOC expansion).
SOCP uses a hard risk cap $0.15$ plus an auxiliary $t_t$.
IR is off on QPs and on for SOCPs
(Clarabel.rs / CuClarabel True, \cpu{} AUTO, \gpu{} AUTO2;
the \Cref{sec:exp-rules} budget is two refinement steps).
Four SOCP solves that returned AlmostSolved at a one-step budget
were re-timed at $\mathrm{IR\_MAX\_ITER}{=}2$ and then solved:
\gpu{} and CuClarabel at $n{=}1000$, $T{=}4$, and
\cpu{} at $n{=}5000$, $T{=}3$ and $T{=}4$.
\Cref{tab:portfolio-first} is the first solver solve on a fresh
object (construction, CVXPY canonicalization, and data updates
excluded);
\Cref{tab:portfolio-fwd} is the mean of ten sequential re-solves
on the reused object
(CVXPY warm start on the first solve and the ten
re-solves; \cpu{}/\gpu{} run each re-solve from the default start,
with no warm start).
\gpu{} is the native forward solver.
A dagger marks a status other than solved/optimal.
CuClarabel's GPU path copies $P$ as-is into a full KKT block
(official examples and \texttt{update\_P!} require a symmetric
CSR, not the upper triangle alone) \cite{cuclarabeljl}; we therefore
pass the symmetrized $P$.

\input{tables/tab_portfolio_first}
\input{tables/tab_portfolio_fwd}

cuDSS reductions are not associative, so ten cold \gpu{}
re-solves need not return bit-identical objectives.
\Cref{tab:portfolio-gpu-obj} records that ten-point sample
(from the default start, no warm start) against Mosek and Gurobi
on the same seed-42 data.
The sample standard deviation is at most
$5.3\times10^{-14}$ on QPs and $9.2\times10^{-10}$ on feasible
SOCPs (relative standard deviation at most $1.2\times10^{-8}$).
On every compared instance the sample width is smaller than
$|f_{\mathrm{Mosek}}-f_{\mathrm{Gurobi}}|$.
Thirteen of twenty-seven samples lie in the closed commercial
interval; the other fourteen lie just beyond Gurobi (a slightly
better objective) by at most $0.22$ of that same gap
(the $n{=}500$, $T{=}2$ SOCP).
The GPU run-to-run scatter is therefore smaller than the
disagreement already present between the two commercial
solvers.

\input{tables/tab_portfolio_gpu_obj}

Mosek and Gurobi appear as 1-thread and 8-thread subcolumns.
On subsequent $T{=}1$ QPs, Clarabel.rs is fastest at $n{=}200$
(3.2\,ms), ahead of 1-thread Mosek (4.4\,ms) and Gurobi
(5.2\,ms); 8 threads do not reduce that small QP (Mosek 4.6\,ms).
On $T{=}1$ QPs the isolated \cpu{} column is ahead of the
original-campaign \gpu{} and Clarabel.rs cells at
$n{=}1000$ (12.3\,ms versus 14.2\,ms and 16.2\,ms); that
15\% gap is of the same order as campaign scatter, so we do
not treat $n{=}1000$ as a hard crossover.
\gpu{} leads from $n{=}2000$ (16.4\,ms versus 24.9\,ms
isolated \cpu{}) through
$n{=}10000$, $k{=}100$ (59\,ms versus 439/406\,ms Mosek
at 1/8 threads, 698/457\,ms Gurobi, 1.09\,s \cpu{},
1.59\,s Clarabel.rs).
Raising $T$ at fixed $(n,k)$ grows the KKT roughly linearly:
the $n{=}5000$ QP subsequent means on \gpu{} are
24.2\,ms ($T{=}1$), 59\,ms ($T{=}2$), 89\,ms ($T{=}3$),
118\,ms ($T{=}4$), while Clarabel.rs goes
98\,ms $\to$ 1.46\,s and Gurobi
122\,ms $\to$ 2.08\,s at 1 thread
(130\,ms $\to$ 1.61\,s at 8 threads).
On feasible SOCPs, 1-thread Mosek is slightly ahead at
$n{=}1000$ (30\,ms versus 32\,ms \gpu{});
\gpu{} leads from $n{=}5000$ (47\,ms) to
$n{=}10000$ (117\,ms versus 610/562\,ms Mosek at 1/8 threads
and 4.73/4.55\,s Gurobi).
CuClarabel is solved on the QPs (51\,ms at $n{=}1000$,
76\,ms at $n{=}5000$, 136\,ms at $n{=}10000$) and on those
SOCPs (190\,ms at $n{=}5000$, 346\,ms at $n{=}10000$);
it is slower than native \gpu{} on this machine; we claim no
general speedup over it.
Gurobi's QCP barrier is not competitive on the SOC form
(0.96/0.90\,s at $n{=}5000$ for 1/8 threads;
4.73/4.55\,s at $n{=}10000$; 11.6/9.61\,s at
$n{=}5000$, $T{=}4$).
Tabled first times are the first solver solve only
(10.2\,ms at the $n{=}200$ QP, 59\,ms at $n{=}10000$ on \gpu{}),
not cuDSS analysis.
Isolated \cpu{} first solves can still show a one-time allocator
or thread-pool spike (for example a first call several times the
subsequent mean); the subsequent column is the ranking quantity.
The opposite can also occur: at the $n{=}10000$ QP the isolated
first call (885\,ms) is faster than the subsequent mean
(1.09\,s), so a single first sample is not a lower bound.
That crossover is for this sparse factor-model family; we do not
claim it on unstructured QP corpora.
On solved $T{=}1$ instances, Clarabel.rs / \cpu{} / \gpu{}
relative objectives agree to $10^{-8}$; the gap versus Mosek is
about $5{\times}10^{-8}$ at $n{=}200$ and $4{\times}10^{-5}$ at
$n{=}10000$.

\begin{table}[t]
\centering
\caption{Native GPU solver versus the
\diffsolver{} wrapper used in the main tables, on the same
multi-period factor-model QP at $T{=}1$
($k{=}50$ at both $n{=}500$ and $n{=}5000$),
subsequent mean of 20 cold solves (one standard deviation in
parentheses).
The wrapper cost is negligible on this path.
The $n{=}5000$ \diffsolver{} column ($23.9$\,ms) lies within one
standard deviation of the $24.2$\,ms cold mean of
\Cref{tab:portfolio-fwd}; the native $25.6$\,ms is higher than the
wrapper by slightly more than that spread, a residual path
difference we have not isolated.}
\label{tab:gpu-path}
\small
\begin{tabular}{@{}lcc@{}}
\toprule
Instance & Native Solver (ms) & \diffsolver{} (ms) \\
\midrule
QP $n{=}500$ & $9.44$ $(0.18)$ & $9.46$ $(0.18)$ \\
QP $n{=}5000$ & $25.6$ $(0.67)$ & $23.9$ $(0.22)$ \\
\bottomrule
\end{tabular}
\end{table}

The small-$n$ GPU times in \Cref{tab:portfolio-fwd} are therefore
the IPM/cuDSS cost, not an artifact of the training wrapper
(\Cref{tab:gpu-path}).

\begin{figure}[t]
\centering
\includegraphics[width=0.92\linewidth]{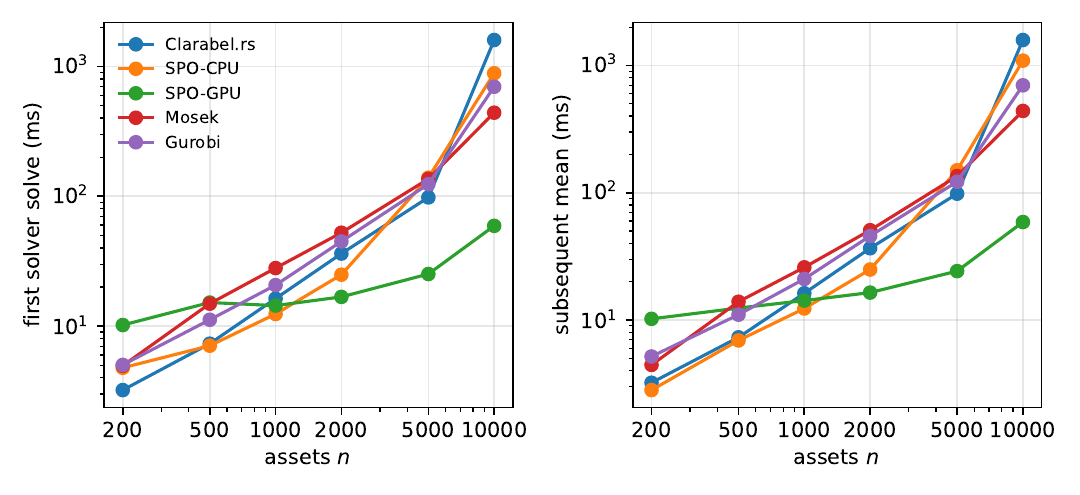}
\caption{The $T{=}1$ QP of the multi-period factor-model family:
tabled first solver solve from \Cref{tab:portfolio-first} (left)
and subsequent mean from \Cref{tab:portfolio-fwd} (right).
The \cpu{} curve is the isolated column; the other solvers are
the original multi-solver campaign.
As in the main tables, the $n{=}100$ pair is omitted from both
panels.}
\label{fig:solver-first-subseq}
\end{figure}

\begin{table}[t]
\centering
\caption{Relative objective difference versus Mosek, and
$\|x-x_{\mathrm{Mosek}}\|_2/\|x_{\mathrm{Mosek}}\|_2$ on asset
weights, on the same instances.
Mosek is a \emph{comparison column}, not a residual oracle:
from $n{=}1000$ Clarabel.rs, \cpu{}, \gpu{}, and Gurobi agree with
each other more tightly than any of them agrees with Mosek
(\Cref{tab:residuals}).
Pairwise relative objectives among Clarabel.rs, \cpu{}, and \gpu{}
(computed from the same campaign's objectives, not shown
as extra columns) are at most $5.3\times10^{-8}$ on QPs and
$2.2\times10^{-7}$ on the full sweep.
All five solvers satisfy $\mathbf{1}^\top x=1$ to at least
$10^{-8}$.
The $n{=}100$ rows are omitted.}
\label{tab:portfolio-acc}
\footnotesize
\setlength{\tabcolsep}{2.4pt}
\begin{tabular}{@{}lcccc@{}}
\toprule
Instance & Clarabel.rs & \cpu{} & \gpu{} & Gurobi \\
\midrule
QP $n{=}200$
  & $8.3{\times}10^{-9}/4.3{\times}10^{-6}$
  & $7.4{\times}10^{-9}/4.4{\times}10^{-6}$
  & $8.4{\times}10^{-9}/4.3{\times}10^{-6}$
  & $4.7{\times}10^{-8}/2.5{\times}10^{-4}$ \\
QP $n{=}500$
  & $3.1{\times}10^{-7}/7.3{\times}10^{-5}$
  & $3.0{\times}10^{-7}/8.1{\times}10^{-5}$
  & $3.1{\times}10^{-7}/7.3{\times}10^{-5}$
  & $2.8{\times}10^{-7}/2.9{\times}10^{-4}$ \\
QP $n{=}1000$
  & $1.6{\times}10^{-6}/3.3{\times}10^{-4}$
  & $1.6{\times}10^{-6}/3.1{\times}10^{-4}$
  & $1.6{\times}10^{-6}/3.3{\times}10^{-4}$
  & $1.6{\times}10^{-6}/2.7{\times}10^{-4}$ \\
QP $n{=}2000$
  & $2.5{\times}10^{-5}/1.2{\times}10^{-3}$
  & $2.5{\times}10^{-5}/1.2{\times}10^{-3}$
  & $2.5{\times}10^{-5}/1.2{\times}10^{-3}$
  & $2.5{\times}10^{-5}/1.2{\times}10^{-3}$ \\
QP $n{=}5000$
  & $5.9{\times}10^{-5}/2.5{\times}10^{-3}$
  & $5.9{\times}10^{-5}/2.6{\times}10^{-3}$
  & $5.9{\times}10^{-5}/2.5{\times}10^{-3}$
  & $5.9{\times}10^{-5}/2.3{\times}10^{-3}$ \\
QP $n{=}10000$
  & $7.9{\times}10^{-5}/3.2{\times}10^{-3}$
  & $7.9{\times}10^{-5}/3.1{\times}10^{-3}$
  & $7.9{\times}10^{-5}/3.2{\times}10^{-3}$
  & $7.9{\times}10^{-5}/3.1{\times}10^{-3}$ \\
QP $n{=}20000$
  & $7.8{\times}10^{-5}/1.8{\times}10^{-3}$
  & $7.8{\times}10^{-5}/1.8{\times}10^{-3}$
  & $7.8{\times}10^{-5}/1.8{\times}10^{-3}$
  & $7.8{\times}10^{-5}/1.7{\times}10^{-3}$ \\
\midrule
SOCP $n{=}200$
  & $5.9{\times}10^{-8}/9.1{\times}10^{-6}$
  & $1.1{\times}10^{-7}/1.6{\times}10^{-6}$
  & $5.9{\times}10^{-8}/9.1{\times}10^{-6}$
  & $6.6{\times}10^{-10}/1.1{\times}10^{-6}$ \\
SOCP $n{=}500$
  & $6.6{\times}10^{-7}/2.1{\times}10^{-4}$
  & $8.6{\times}10^{-7}/2.2{\times}10^{-4}$
  & $8.7{\times}10^{-7}/2.2{\times}10^{-4}$
  & $8.5{\times}10^{-7}/2.4{\times}10^{-4}$ \\
SOCP $n{=}1000$
  & $1.9{\times}10^{-6}/8.2{\times}10^{-4}$
  & $1.9{\times}10^{-6}/8.2{\times}10^{-4}$
  & $1.9{\times}10^{-6}/8.2{\times}10^{-4}$
  & $1.9{\times}10^{-6}/8.0{\times}10^{-4}$ \\
SOCP $n{=}2000$
  & $5.7{\times}10^{-7}/3.1{\times}10^{-4}$
  & $6.1{\times}10^{-7}/3.1{\times}10^{-4}$
  & $5.7{\times}10^{-7}/3.1{\times}10^{-4}$
  & $6.2{\times}10^{-7}/2.3{\times}10^{-4}$ \\
SOCP $n{=}5000$
  & $1.2{\times}10^{-5}/2.2{\times}10^{-3}$
  & $1.2{\times}10^{-5}/2.2{\times}10^{-3}$
  & $1.2{\times}10^{-5}/2.2{\times}10^{-3}$
  & $1.2{\times}10^{-5}/2.2{\times}10^{-3}$ \\
SOCP $n{=}10000$
  & $2.9{\times}10^{-5}/2.7{\times}10^{-3}$
  & $2.8{\times}10^{-5}/2.7{\times}10^{-3}$
  & $2.9{\times}10^{-5}/2.7{\times}10^{-3}$
  & $2.9{\times}10^{-5}/2.7{\times}10^{-3}$ \\
SOCP $n{=}20000$
  & $1.0{\times}10^{-4}/2.1{\times}10^{-3}$
  & $1.0{\times}10^{-4}/2.1{\times}10^{-3}$
  & $1.0{\times}10^{-4}/2.1{\times}10^{-3}$
  & $1.0{\times}10^{-4}/2.0{\times}10^{-3}$ \\
\bottomrule
\end{tabular}
\end{table}

\input{tables/tab_n20000}
\input{tables/tab_residuals}

On QPs, Clarabel.rs, \cpu{}, and \gpu{} agree pairwise to at most
$5.3\times10^{-8}$ in relative objective (the \cpu{}/Clarabel.rs
maximum, at $n{=}10000$).
Across the whole QP+SOCP sweep the three stay within
$2.2\times10^{-7}$ (the maximum is \gpu{} versus Clarabel.rs on
the $n{=}500$ SOCP).
Versus native Mosek QP the gap grows
from $10^{-8}$--$10^{-7}$ at $n\le 500$ to $7.8{\times}10^{-5}$ at
$n{=}20000$.
From $n{=}1000$ onward that offset is shared by Clarabel.rs,
\cpu{}, \gpu{}, \emph{and Gurobi}: the four objectives agree with
each other to $\sim 10^{-8}$ relative, while Mosek's quadratic
iterate is strictly worse (higher, since we minimize).
At $n{=}20000$ the four report $-0.11464358$ and Mosek reports
$-0.11463462$.
Mosek is not the short-iteration outlier: it uses 16 barrier
iterations against \cpu{}'s 13, \gpu{}/Clarabel.rs's 18, and
Gurobi's 20
(\Cref{tab:n20000}), after we set every Mosek interior-point
gap/feasibility parameter to $10^{-8}$ (\Cref{sec:exp-rules}).
The four-way agreement including Gurobi shows the split is
Mosek-specific, not a Clarabel KKT bug.
The remaining explanation is that the two codes terminate on
different merit functions: Clarabel on homogeneous-embedding
residuals after Ruiz equilibration, Mosek on its native QP
barrier with its own scaling.
Those tests can both read as ``$10^{-8}$'' and still stop at
distinct $(x,s,z)$.
\Cref{tab:residuals} reports the Clarabel-form residuals
\eqref{eq:kkt} on the same standard form: all three
Clarabel.rs, \cpu{}, and \gpu{} sit at $10^{-9}$--$10^{-12}$, two to five
orders inside the $10^{-8}$ protocol, including at $n{=}20000$.
Those residuals, not the objective gap versus Mosek, are the
accuracy certificate.
A Mosek Task on the same $(P,q,A,b)$ reports
\texttt{solsta.optimal} with native relative gap
$7.5{\times}10^{-10}$ at $n{=}5000$ and $2.0{\times}10^{-8}$ at
$n{=}20000$ (\Cref{tab:residuals}); the latter sits about
$2\times$ above the written $10^{-8}$ protocol.
We do not convert that iterate into Clarabel-form $(s,z)$.
On SOCPs the objective gaps are $10^{-8}$--$10^{-4}$ and the three
Clarabel.rs / \cpu{} / \gpu{} solvers stay within a few $10^{-9}$ of each other;
\gpu{} terminates as solved (AUTO2).
Portfolio weights stay within $10^{-2}$ relative $\ell_2$ of Mosek
at $n{=}20000$.
In this section $\mu$ denotes the expected-return vector (the
objective block $q_{1:n}=-\mu$ of \Cref{sec:portfolio-layer}); the
IPM centrality is written $\mubar$.

\begin{figure}[t]
\centering
\includegraphics[width=0.96\linewidth]{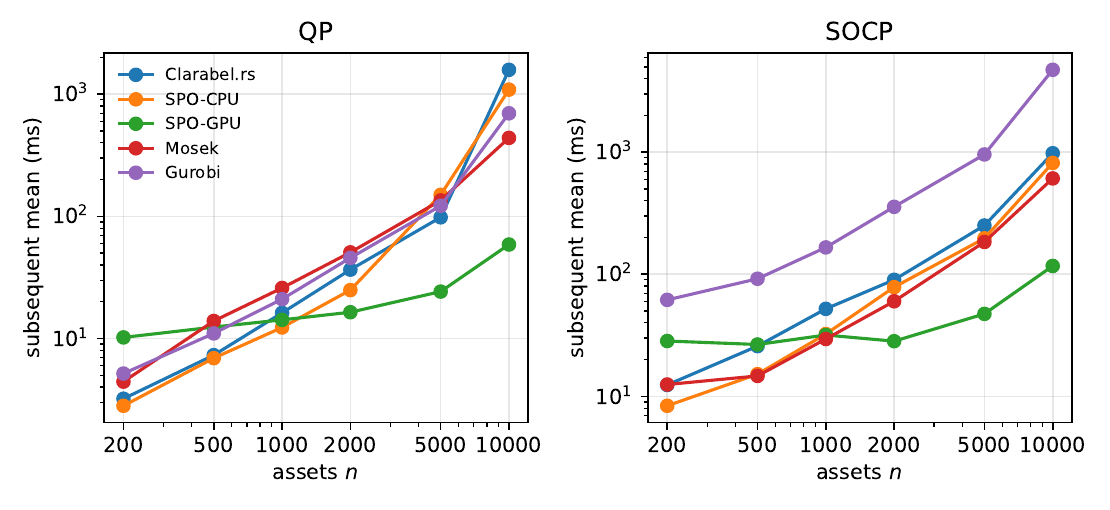}
\caption{Subsequent-mean solve time versus number of assets for the
QP (left) and SOCP (right) families, drawn from
\Cref{tab:portfolio-fwd}
($k{=}50$ up to $n{=}5000$, and $k{=}100$ at $n{=}10000$).
The \cpu{} curve is the isolated column; the other solvers are
the original multi-solver campaign.
As in the main tables, the $n{=}100$ pair is omitted.}
\label{fig:portfolio-scaling}
\end{figure}

\begin{table}[t]
\centering
\caption{Parametric re-solves on the $n{=}5000$, $k{=}50$ QP:
cold first solve, then median of 20 $q$-updates that
perturb only the asset block $q_{1:n}=-\mu$ (scale $0.01$).
Iteration counts in parentheses where recorded.
Clarabel.rs~0.11.1 has no smoothing-operator setting and this path
did not record its iteration count.}
\label{tab:portfolio-dpp}
\small
\begin{tabular}{@{}lccc@{}}
\toprule
Solver & First solve (ms) & Steady $q$-update & Steady + warm start \\
\midrule
\cpu{} & 131 (11) & 133 (11) & 125 (11) \\
\gpu{} & 41.7 (13) & 26.1 (14) & 22.7 (13) \\
Clarabel.rs & 112 & 106 & --- \\
\bottomrule
\end{tabular}
\end{table}

Smoothing on the $n{=}5000$ QP cuts median $q$-update time
from 133\,ms to 125\,ms on \cpu{} (iterations stay at 11) and
from 26.1\,ms to 22.7\,ms on \gpu{} (13\%), matching a
one-iteration drop (14$\to$13) on \gpu{}.
The CPU 6\% wall-clock change is of the same order as the
three-seed scatter of \Cref{sec:exp-rules} (about 5\% on this
size) and is not treated as a measured speedup.
The GPU first-solve 41.7\,ms includes construction of the
\diffsolver{} wrapper; the 26.1\,ms steady column is the subsequent
$q$-update path, not a copy of
\Cref{tab:portfolio-fwd}'s cold mean.
This parametric protocol is reported at one size; the $n{=}1000$
sweep of \Cref{sec:exp-warm} is the $n$-resolved iteration
comparison for updates of $q$.
Unstructured noise on the full canonical $q$ (including equality
multipliers) can mark the QP dual infeasible; all parametric rows
perturb only $-\mu$.

\subsection{Warm start}
\label{sec:exp-warm}

The smoother of \Cref{sec:warmstart} is a claim about
\emph{consecutive} nearby QCPs, not about a single solve.
We therefore test two protocols on the same $n{=}1000$, $k{=}50$
factor-model family, in both the QP and the SOCP formulation, on
\cpu{} and \gpu{}.
The SOCP carries two second-order cones of dimension $1051$ each.
Cold start is a default start on a reused solver object
(smoothing off); smoothing is Algorithm~1 of
\cite{chen2025warmstart} with the documented modifications of
\Cref{sec:smoothing-impl}
($\tau=1$ slice, feasibility-only $\mubar_0$, Theorem~4.1 for SOC).
The solver has no unsmoothed-restore mode, so we do not report a
``raw reuse'' row: disabling smoothing is cold start.
Following \cite{chen2025warmstart}, the reduction ratio on a batch of
problems is the geometric mean
\begin{equation}
\label{eq:Riter}
R
=
\exp\Bigl(\tfrac1T\sum_{t=1}^{T}
\log\frac{\mathrm{iters}^{\mathrm{smooth}}_t}{\mathrm{iters}^{\mathrm{cold}}_t}\Bigr),
\end{equation}
so $R<1$ is a reduction.
Only $-\mu$ is updated; $x_0=0$ is held
fixed, as in \Cref{tab:portfolio-dpp}.
Wall time tracks the iteration count on both backends; we quote
milliseconds in the text and keep the tables in iterations.

\begin{figure}[t]
\centering
\includegraphics[width=0.96\linewidth]{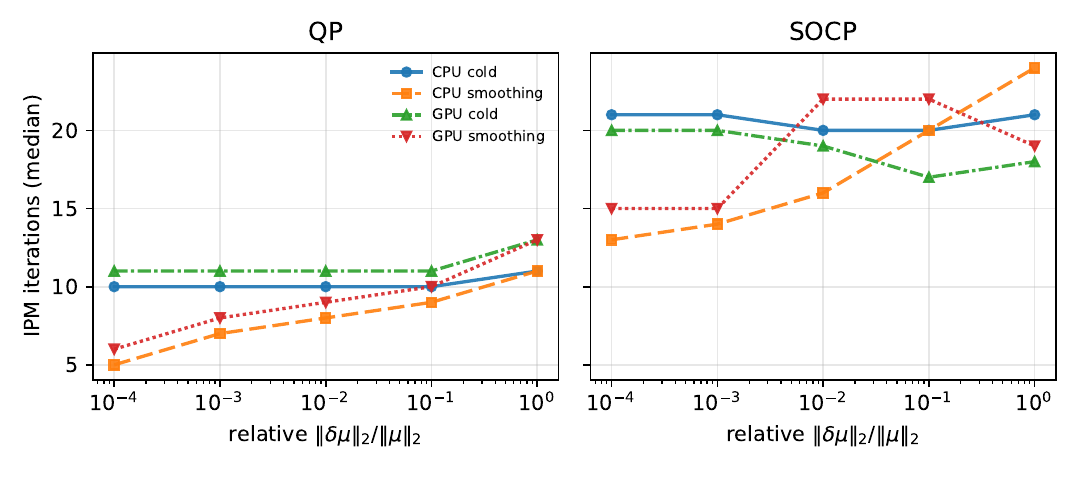}
\caption{Median IPM iterations versus relative $\mu$ perturbation
(\Cref{tab:warmstart}).
Left: QP, a monotone saving that vanishes at a jump of size $1$;
the CPU and GPU curves track each other within one iteration
across the sweep.
Right: SOCP, a small-jump saving on the CPU and a GPU iteration
increase at $10^{-2}$.}
\label{fig:exp-warmstart}
\end{figure}

\paragraph{One-shot perturbation.}
A cached solve of the original instance is followed by one relative
jump $\mu\leftarrow\mu+\delta$ with
$\|\delta\|_2/\|\mu\|_2\in\{10^{-4},10^{-3},10^{-2},10^{-1},1\}$.
\Cref{tab:warmstart} reports the median of five independent draws;
\Cref{fig:exp-warmstart} plots the same numbers.
On the QP (zero and nonnegative cones only) the two backends track
each other to within one iteration (\Cref{tab:warmstart}).
Smoothing halves the CPU count at a jump of $10^{-4}$
($10\to 5$, $R=0.50$) and cuts the GPU count from $11$ to $6$
($R\approx 0.55$).
The saving then shrinks monotonically with the jump
(CPU $10\to 7/8/9$; GPU $11\to 8/9/10$) and vanishes at a jump of
size $1$ (CPU $11\to 11$, GPU $13\to 13$), which is the prediction of
\cite{chen2025warmstart}.

On the SOCP the behaviour is cone- and backend-dependent, not a
uniform loss.
At $10^{-4}$, Theorem~4.1 with a tiny feasibility $\mubar_0$
reduces the iteration count on both solvers (\cpu{} $21\to 13$,
$R\approx 0.62$; \gpu{} $20\to 15$, $R=0.75$).
At $10^{-3}$ the CPU still saves seven iterations ($21\to 14$)
while the GPU saves five ($20\to 15$).
At $10^{-2}$---the DPP/rebalance scale---the CPU still saves four
iterations ($20\to 16$) but the GPU requires three more
($19\to 22$): this is the jump size at which the dual residual
that sets $\mubar_0$ starts to lift $s_0$ enough to disturb the
four-tier SOC scaling, even with the relative cap of
\Cref{sec:smoothing-impl}.
At $10^{-1}$ the CPU is neutral ($20\to 20$) and the GPU is worse
($17\to 22$); at a jump of size $1$ the CPU is worse ($21\to 24$)
and the GPU is essentially neutral ($18\to 19$).
The earlier one-shot protocol that restored a complementary SOC at
$\lambda=1$ (Appendix~B) or that folded the equilibrated duality
gap into $\mubar_0$ increased the count at every jump; those
failure modes are not the SOC itself.

\begin{table}[t]
\centering
\caption{One-shot warm start: median IPM iterations over five
independent relative perturbations of $\mu$ on the $n{=}1000$,
$k{=}50$ QP and SOCP.
Cold is a default start on a reused object; smoothing is
Algorithm~1 with the coordinate choices of
\Cref{sec:smoothing-impl}.
On the QP the two backends track each other closely, and the
saving shrinks as the jump grows.
On the SOCP a small jump helps; at $10^{-2}$ the CPU still saves
iterations while the GPU requires a few more.}
\label{tab:warmstart}
\small
\begin{tabular}{@{}llccccc@{}}
\toprule
Family & Method & $10^{-4}$ & $10^{-3}$ & $10^{-2}$ & $10^{-1}$ & $1$ \\
\midrule
QP & \cpu{} cold & 10 & 10 & 10 & 10 & 11 \\
QP & \cpu{} smoothing & 5 & 7 & 8 & 9 & 11 \\
QP & \gpu{} cold & 11 & 11 & 11 & 11 & 13 \\
QP & \gpu{} smoothing & 6 & 8 & 9 & 10 & 13 \\
\midrule
SOCP & \cpu{} cold & 21 & 21 & 20 & 20 & 21 \\
SOCP & \cpu{} smoothing & 13 & 14 & 16 & 20 & 24 \\
SOCP & \gpu{} cold & 20 & 20 & 19 & 17 & 18 \\
SOCP & \gpu{} smoothing & 15 & 15 & 22 & 22 & 19 \\
\bottomrule
\end{tabular}
\end{table}

\paragraph{Sequential rebalancing.}
The intended application is a chain: $\mu_t=\mu_{t-1}+\delta_t$ with
$\|\delta_t\|_2/\|\mu_{t-1}\|_2$ fixed, one solver object, $T=20$
steps after the original solve.
\Cref{tab:warmstart-seq} summarizes mean iterations and $R$
at relative steps $10^{-2}$ and $1$.
At relative step $10^{-2}$ the QP saving is immediate and stable,
and larger than the one-shot at the same scale because every step
stays near the previous active set
(\cpu{} $10.0\to 8.0$, $R=0.80$, $12.3\to 9.2$\,ms;
\gpu{} $11.0\to 9.0$, $R=0.82$, $15.0\to 11.6$\,ms).
The CPU SOCP chain exhibits the same pattern at a smaller margin
($20.0\to 15.3$, $R=0.77$, $30.7\to 22.5$\,ms), with the smoothed
count remaining in the range 15--17 from the first step.
The GPU SOCP chain does \emph{not} improve: mean iterations rise
from $18.0$ to $21.7$ ($R=1.20$, $30.2\to 35.8$\,ms), and the
smoothed per-step trace lies above the cold trace throughout the
chain.
That is the sequential counterpart of the one-shot GPU iteration
increase at $10^{-2}$: the relative cap on $\mubar_{\mathrm{soc}}$
stops a
$10\%$ lift of $s_0$ from turning into a large iteration spike, but
it does not make sparse GPU SOC scaling cheaper than a cold unit
start on this cone.
We do not claim a match to the sequential SOCP portfolios of
\cite{chen2025warmstart} ($R_{\mathrm{iter}}\approx 0.63$).

A relative step of $1$ is a stress test.
The walk makes the QPs harder (cold climbs from 10 to $\sim$19
iterations); smoothing no longer helps on average
(\cpu{} $18.2\to 19.0$, $R=0.99$; \gpu{} $19.1\to 18.7$, $R=0.95$),
and later steps are noisy and the smoothed iterates can
overshoot cold.
The same large walk on the SOCP is at best neutral on the CPU
($22.8\to 23.1$, $R=1.02$) and worse on the GPU
($25.9\to 27.4$, $R=1.06$), with late GPU steps drifting to 24--27
iterations.

\paragraph{A quantitative account of the SOCP failure.}
Theorem~4.1 replaces the Lorentz scalar by
$s_0\leftarrow\sqrt{\mubar_{\mathrm{soc}}+\|c_1\|_2^2/4}$ with
$\mubar_{\mathrm{soc}}=\lambda\mubar_0$ and
$\lambda=s_0^\star/z_0^\star$
(the SOC smoothing kernel, on both backends).
At complementarity $\|c_1\|_2^2/4\approx s_0^2$, so a
feasibility-sized $\mubar_0$ after a relative $\mu$ jump of
$10^{-2}$ produces a relative lift
$\sqrt{1+\mubar_{\mathrm{soc}}/s_0^2}-1$ of a few percent on a
cone of dimension $1051$.
The GPU path additionally \emph{caps}
$\mubar_{\mathrm{soc}}\le 10^{-4}\max(\|c_1\|_2^2/4,s_0^2)$
because an uncapped $\lambda\mubar_0$ was observed to lift $s_0$
by $\sim$10\% and to add Mehrotra iterations in the four-tier
sparse SOC scaling; the CPU uses the matching floor only
(\Cref{sec:smoothing-impl}).
That cap is an algorithmic variant of Theorem~4.1, not a
device-side rewrite of the same map: CPU and GPU therefore
execute different SOC warm starts.
\Cref{tab:soc-cap} isolates the cap on the CPU: sequential
$R$ moves only from $0.77$ (floor) to $0.85$ (cap), still a
reduction, while the GPU capped chain stays at $R=1.20$.
The $2{\times}2$ of \Cref{tab:soc-cap} separates the operator
from the backend: CPU stays at $R<1$ with or without the cap
($0.77$ floor, $0.85$ cap); GPU stays at $R>1$ with or without
it, and turning the cap off makes the sequential GPU cell
worse ($R=1.34$ versus $1.20$).
The cap is therefore neither a sufficient nor a necessary cause
of the GPU failure; the remaining effect is the GPU path's
interaction of Theorem~4.1's lift of $s_0$ with four-tier
sparse SOC scaling.
The cap prevents a large spike but does not make the scaled
GPU SOC cheaper than a cold unit start: one-shot GPU iterations
rise from 19 to 22 at jump $10^{-2}$ (\Cref{tab:warmstart}), and
the sequential chain remains above cold
(\Cref{tab:warmstart-seq});
the isolation table is \Cref{tab:soc-cap}.
A decision rule on this family is therefore:
\begin{itemize}[leftmargin=1.2em,itemsep=0.1em]
\item QP, any jump we tested: smoothing on;
\item SOCP, $\|\Delta\mu\|_2/\|\mu\|_2\lesssim 10^{-3}$: on;
\item SOCP, CPU sequential rebalance at $10^{-2}$: on;
\item SOCP, GPU sequential rebalance at $10^{-2}$, or any
      $O(1)$ jump on a large SOC: off.
\end{itemize}

\paragraph{A parameter-agnostic problem distance.}
The rule above is stated in terms of the single perturbation we
swept, $\|\Delta\mu\|_2/\|\mu\|_2$.
In a multi-period application the return vector is not the only
input that changes: the initial holdings $w_0$ roll over every
period, and the factor loadings are re-estimated.
A per-parameter sweep of the kind in \Cref{fig:exp-warmstart} does
not determine whether a warm start is beneficial once several
inputs change simultaneously.
We therefore replace the parameter axis by the problem distance
$d_{\mathrm{KKT}}$ of \Cref{def:dkkt}, which aggregates every
perturbation into one scalar and is independent of whether the
change came from $\mu$, from $w_0$, or from the constraint
matrices.
As a geometric complement we record the fraction $d_{\mathrm{act}}$
of holdings whose active bound flips, which detects a changed
active set even when the residual stays small.

\begin{figure}[t]
\centering
\includegraphics[width=0.96\linewidth]{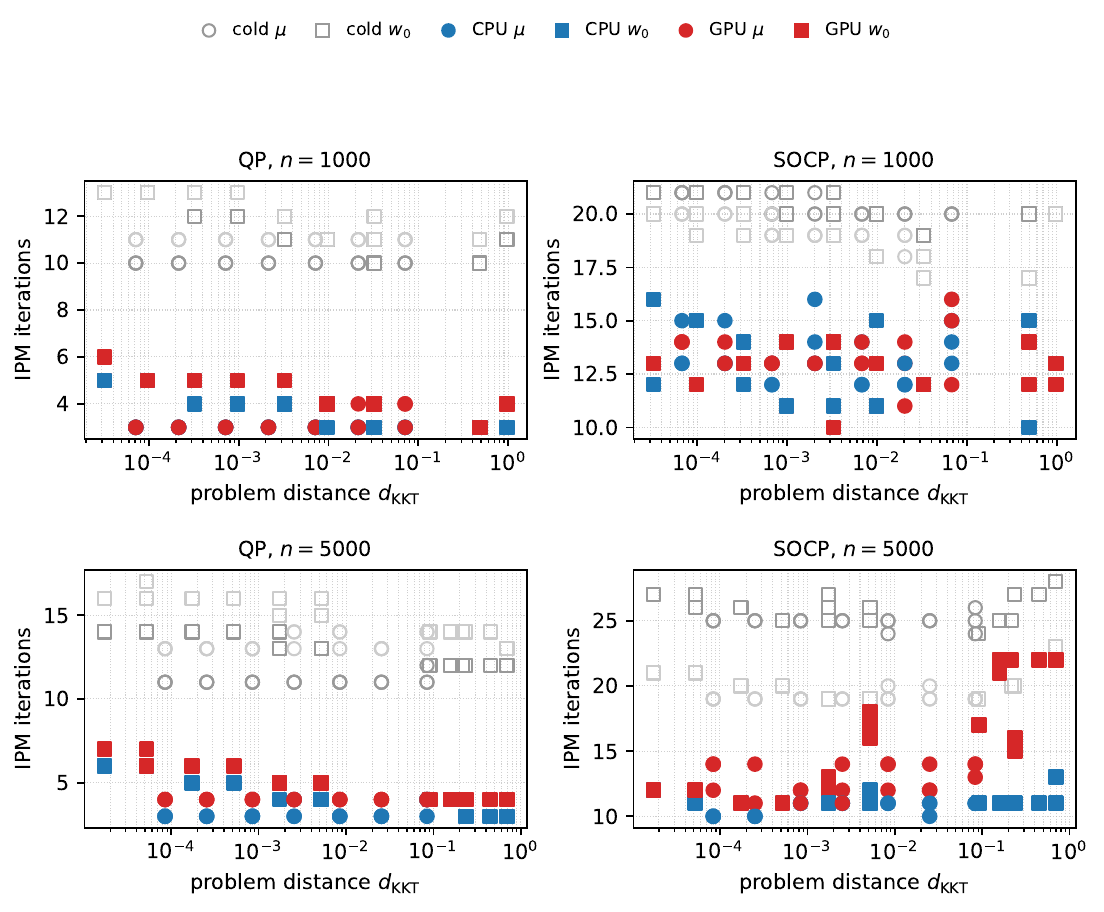}
\caption{IPM iterations versus the problem distance
$d_{\mathrm{KKT}}$ of \eqref{eq:dkkt} on the $k{=}50$ family at
$n{=}1000$ (top) and $n{=}5000$ (bottom), perturbing either the
return vector $\mu$ (circles) or the initial holdings $w_0$
(squares), on the CPU (blue) and GPU (red) backends.
Marker \emph{shape} encodes the perturbation source and marker
\emph{colour} the backend; open grey markers are cold restarts
(open circles for $\mu$, open squares for $w_0$), filled coloured
markers are smoothed warm starts.
At both sizes the two sources interleave on the shared axis with no
systematic offset, so the distance---not the identity of the
perturbed parameter---organises the warm-start cost, and the CPU
smoothed trace stays well below its cold trace throughout.
The right-most points are $w_0$ perturbations whose residual is an
order of magnitude larger than their nominal step would suggest;
$d_{\mathrm{KKT}}$ identifies them.
At $n{=}1000$ the only place the GPU smoothing approaches its cold
trace is the cluster of large-$\mu$ steps near
$d_{\mathrm{KKT}}\approx7\cdot10^{-2}$, so there the GPU SOCP
overhead is driven by the objective component alone.
At $n{=}5000$ the GPU trace also rises toward cold on the
large-$d_{\mathrm{KKT}}$ $w_0$ outliers (bottom right), while the
CPU still saves $13$--$17$ iterations there: the larger sparse SOC
system makes the GPU warm start sensitive to the geometric
component as well.}
\label{fig:exp-warmstart-dist}
\end{figure}

\Cref{fig:exp-warmstart-dist} sweeps both $\mu$ and $w_0$ on the
$T{=}1$ family and plots iterations against $d_{\mathrm{KKT}}$.
Three observations support using the distance as the decision
variable.
First, the two sources \emph{interleave} on the shared axis: at
$d_{\mathrm{KKT}}\approx10^{-3}$ a $\mu$ step and a $w_0$ step give
indistinguishable smoothed iteration counts, so a single threshold
in $d_{\mathrm{KKT}}$ governs both.
Second, the two metrics separate cause: $\mu$ steps leave
$d_{\mathrm{act}}\approx0$ until the relative step reaches
$10^{-1}$, whereas $w_0$ steps flip $1.5$--$2\%$ of the active
bounds at every magnitude---$w_0$ moves the geometry, $\mu$ mostly
does not---yet the smoothed count tracks $d_{\mathrm{KKT}}$, not
$d_{\mathrm{act}}$.
Third, the residual is an indicator rather than a failure signal:
the $w_0$ perturbations whose $d_{\mathrm{KKT}}$ is an order of
magnitude above their nominal step (the right-most squares) are
precisely the cases in which the cached point lands near a
constraint corner, and smoothing still saves five to ten iterations
there on \emph{both} backends.
Fourth, the distance localises the GPU SOCP overhead of
\Cref{tab:warmstart,tab:warmstart-seq}: the only cluster where the
GPU smoothed trace rises to meet its cold trace is the large-$\mu$
steps near $d_{\mathrm{KKT}}\approx7\cdot10^{-2}$, where the saving
falls from five--eight iterations to zero--three.
The comparably large $d_{\mathrm{KKT}}$ produced by $w_0$ does not
trigger it, which isolates the cause to the objective-side lift of
$s_0$ that the four-tier sparse SOC scaling must absorb
(\Cref{sec:smoothing-impl}), not to a changed active set.
The single-axis sweep of \Cref{fig:exp-warmstart} cannot separate
these two, because it never varies $w_0$.
The per-parameter thresholds of the rule above therefore reduce to
a heuristic on this family only: apply a warm start
when $d_{\mathrm{KKT}}\lesssim 10^{-2}$, and treat
$d_{\mathrm{KKT}}\gtrsim 5{\times}10^{-2}$---from any source---as
evidence that the cached point no longer lies near the new central
path (the GPU SOCP trace meets cold near
$d_{\mathrm{KKT}}\approx 7{\times}10^{-2}$ in
\Cref{fig:exp-warmstart-dist}).
We do not claim this threshold on other models.
On the GPU, the benefit should be discounted when the
objective ($\mu$) component of the perturbation is the large one.

\paragraph{Dependence on equilibration.}
All warm-start numbers above use the solver default
\texttt{equilibrate\_enable}\,=\,true.
Because this factor-model family is already well scaled (the ratio
of the largest to smallest row norm of $A$ is about $25$, and the
column norms span a factor of about $4.5$), equilibration has little
effect on conditioning, and one might expect disabling it to be
neutral.
\Cref{tab:warmstart-equil} shows that it is not, and the effect
splits along cone families.
On the QP the saving survives with equilibration off, merely
shrinking (CPU $R$: $0.50\!\to\!0.64$ at a $10^{-4}$ step): the
extra equilibration sweeps add a few cold iterations, so part of
the reported QP gap is scaling overhead, but the smoothed start
still reduces the iteration count.
On the SOCP the saving \emph{reverses}: with equilibration off the
cold count drops from $21$ to $14$ while the smoothed count climbs
above it ($R$ flips from $0.63$ to $1.43$ at $10^{-4}$, and stays
above $1$ across the sweep).
The cause is that the SOC smoothing operator of Theorem~4.1---the
replacement of the Lorentz scalar and the $\bar\mu_{\mathrm{soc}}$
formula---is computed in the equilibrated coordinate system.
This dependence is visible in the theory itself: the residual bound
\eqref{eq:socbound} of \Cref{thm:cgj}
is stated in the coordinate system in which the operator is
evaluated, and it depends explicitly on the scale-dependent
quantities $\|A\|_{\infty}$, the in-cone ratio
$\lambda=s_{0}^{\star}/z_{0}^{\star}$, and the Lorentz scalar
$s_{0}^{\star}$ (\Cref{rem:cgj-scope}).
Equilibration is a diagonal rescaling
$(\hat P,\hat q,\hat A,\hat b)=(D_cPD_c,\,D_cq,\,D_rAD_c,\,D_rb)$,
so it changes exactly these quantities---$\|A\|_{\infty}$ becomes
$\|D_rAD_c\|_{\infty}$ and the smoothed $(s^{0},z^{0})$ inherit the
scaling.
With equilibration disabled, the cached SOC point no longer matches
the scale at which the smoothing operator assumes it was computed,
and the warm start lands farther from the central path than a
default start.
The original analysis of \cite{chen2025warmstart} is carried out in
a fixed coordinate system and does not address this interaction with
a preceding equilibration step; the reversal in
\Cref{tab:warmstart-equil} is the practical consequence.
The $d_{\mathrm{KKT}}$ distance itself is invariant to equilibration
(\Cref{prop:dkkt-invar}), so
\eqref{eq:dkkt} still predicts the \emph{cold} cost correctly; what
changes is whether the smoothed construction actually reaches the
point that the distance assumes.
We therefore read the SOCP warm-start benefit as conditional on
running with the solver's default equilibration, and we report both
settings for transparency rather than as a recommendation to disable
it: on ill-scaled problems equilibration is typically essential, and
removing it would change the sign of these cold baselines.

The qualitative pattern is independent of scale.
\Cref{tab:warmstart-equil} repeats the test at $n{=}5000$ alongside
$n{=}1000$.
The QP saving survives with equilibration off at both sizes, and
the SOCP reversal appears at both sizes (CPU $R$ at a $10^{-2}$
step is $1.71$ at $n{=}1000$ and $1.26$ at $n{=}5000$).
Two scale-dependent shifts remain.
First, disabling equilibration removes more cold iterations as $n$
grows (the SOCP cold count falls $21\to14$ at $n{=}1000$ but
$25\to18$ at $n{=}5000$), consistent with the data remaining no
worse conditioned while the equilibration sweeps themselves become
more expensive.
Second, the GPU SOCP overhead appears earlier and on a larger
fraction of the sweep: at $n{=}5000$ the GPU trace already exceeds
cold at a $10^{-3}$ step, and---unlike at $n{=}1000$---it also
degrades on the large-$d_{\mathrm{KKT}}$ $w_0$ outliers of
\Cref{fig:exp-warmstart-dist}, where the CPU still saves $13$--$17$
iterations.
Larger sparse SOC systems therefore make the GPU warm start more
sensitive to \emph{both} components of the perturbation, not only
the objective one, while the CPU behaviour is essentially
scale-independent.

\begin{table}[t]
\centering
\caption{Warm-start sensitivity to equilibration on the $k{=}50$
family at $n{=}1000$ and $n{=}5000$: median IPM iterations,
cold\,$\to$\,smooth, with \texttt{equilibrate\_enable} on (the
default used elsewhere) and off.
The QP saving survives with equilibration off at both sizes, only
shrinking.
The SOCP saving reverses at both sizes, because the SOC smoothing
operator is computed in equilibrated coordinates; the well-scaled
data makes the off-diagonal cold start artificially cheap, which is a
property of this family, not a general recommendation.
The GPU SOCP overhead starts earlier and extends to the large-$w_0$
outliers at the larger size, while the CPU behaviour is
scale-independent.}
\label{tab:warmstart-equil}
\small
\input{tables/tab_warmstart_equil}
\end{table}

\begin{table}[t]
\centering
\caption{Sequential rebalancing: $T=20$ random-walk updates of $\mu$
on the $n{=}1000$, $k{=}50$ family.
Mean IPM iterations and geometric-mean reduction ratio $R$ of
\eqref{eq:Riter}.
Relative step $10^{-2}$ is the DPP/rebalance regime; step $1$ is a
stress test.
GPU SOCP at $10^{-2}$ is the one cell with $R>1$.}
\label{tab:warmstart-seq}
\small
\begin{tabular}{@{}llcccc@{}}
\toprule
Family & Step & Backend & Mean cold & Mean smooth & $R$ \\
\midrule
QP   & $10^{-2}$ & \cpu{} & 10.0 & 8.0  & 0.80 \\
QP   & $10^{-2}$ & \gpu{} & 11.0 & 9.0  & 0.82 \\
SOCP & $10^{-2}$ & \cpu{} & 20.0 & 15.3 & 0.77 \\
SOCP & $10^{-2}$ & \gpu{} & 18.0 & 21.7 & 1.20 \\
\midrule
QP   & $1$       & \cpu{} & 18.2 & 19.0 & 0.99 \\
QP   & $1$       & \gpu{} & 19.1 & 18.7 & 0.95 \\
SOCP & $1$       & \cpu{} & 22.8 & 23.1 & 1.02 \\
SOCP & $1$       & \gpu{} & 25.9 & 27.4 & 1.06 \\
\bottomrule
\end{tabular}
\end{table}

\subsection{Differentiable layer}
\label{sec:exp-diff}

\paragraph{Accuracy.}
The implicit-function Jacobian is differentiable only on regions of
constant projection type.
We therefore check $\partial(\mathbf{1}^\top x)/\partial q$ from the
in-solver VJP against finite differences of the same scalar, with
$\varepsilon=10^{-6}$ (a step of
$10^{-7}$ lies inside the $10^{-8}$ IPM tolerance).
Cosine similarity is scale-invariant, so \Cref{tab:diff} also reports
the relative error $\|\hat g-g_{\mathrm{FD}}\|_2/\|g_{\mathrm{FD}}\|_2$
on the sampled coordinates.
Interior coordinates are those with $x_i\in(10^{-6},0.1-10^{-6})$;
a coordinate is dropped if the two one-sided differences disagree
(a kink of the solution map).
On QPs that filter removes nothing and both backends match FD in
direction (cosine $0.9998$--$1.000$) with relative error
$8{\times}10^{-4}$ at $n{=}500$ and $3$--$5{\times}10^{-2}$ at
$n{=}5000$---an amplitude discrepancy that cosine similarity
conceals, and large enough to affect a poorly scaled training
loss.
On the $n{=}500$ SOCP the CPU and GPU matches are comparable
(cosine $1.000$ and $0.999$, relative error $5.2{\times}10^{-2}$ and
$5.4{\times}10^{-2}$).
We do not treat the SOCP figures as a complete adjoint
validation.

Bound assets ($x_i\le 10^{-6}$ or $x_i\ge 0.1$) are the
predict-then-optimize regime for which the layer is intended.
There the VJP is numerically zero ($10^{-12}$--$10^{-9}$), which is
the IFT Jacobian of a strictly active bound: $\mathbf{1}^\top x$
does not move with $q_i$ to first order.
A one-sided FD of size $10^{-6}$ is \emph{not} a subgradient oracle
here---it is dominated by IPM noise and active-set variation, so
cosine against that FD is not a meaningful check, and we do not
report it as one.
The layer therefore implements the smooth-branch Jacobian and
returns a particular element of the Clarke subdifferential at a
kink (typically $0$ on a strictly active bound).
It does not expose the full subdifferential.
Finite-difference checks of $\partial L/\partial P$ and
$\partial L/\partial A$ on a small QP are in
\Cref{tab:vjp-ratios}.

\paragraph{Sources of the relative error.}
The check is a VJP of $\mathbf{1}^\top x$ against a central
finite difference of the same scalar with
$\varepsilon=10^{-6}$, while the IPM gap/feasibility tolerance is
$10^{-8}$.
A textbook central difference of a map evaluated to absolute
accuracy $\delta$ trades $O(\varepsilon^2)$ truncation against
$O(\delta/\varepsilon)$ solver noise.
\Cref{tab:eps-sweep} does \emph{not} show that tradeoff: the
relative error is constant to well under one part in $10^{2}$
across four decades of $\varepsilon$ on every configuration.
A noise term $O(\delta/\varepsilon)$ would rise by $10^{3}$ from
$10^{-4}$ to $10^{-7}$; a truncation term $O(\varepsilon^2)$ would
fall by $10^{6}$.
Neither appears.
The residual is therefore a systematic, $\varepsilon$-independent
bias between the reduced two-block VJP and this FD oracle---an
$O(\delta)$ solver-accuracy floor and/or a stack difference, not
an FD-step artefact---and the $\varepsilon{=}10^{-6}$ of
\Cref{tab:diff} lies inside that flat regime.
Per-coordinate ratios $\hat g_i/g_{\mathrm{FD},i}$
(\Cref{tab:vjp-ratios}) have median near $1$ and an IQR of a few
percent, but tails reach $O(1)$ (one $n{=}5000$ sample at $0.017$).
That is coordinate scatter, not a single global scale.
Finite-difference checks of $\partial L/\partial P$ and
$\partial L/\partial A$ on a small QP are reported in the same
table.
We therefore interpret the check as an order-of-magnitude
consistency argument that the descent direction is correct
(cosine $1.000$ on every QP row), not a proof that the VJP is
exact to better than $10^{-2}$.

\input{tables/tab_eps_sweep}
\input{tables/tab_vjp_ratios}
Second, the CPU and GPU relative errors at $n{=}500$ QP agree
($8.8{\times}10^{-4}$ versus $7.6{\times}10^{-4}$), so the gap is
not a KLU-versus-cuDSS artefact.
Third, the GPU SOCP row ($5.4{\times}10^{-2}$ after dropping
kinked coordinates; \Cref{tab:diff}) is the one case in which $\Dpi$ lies near an
SOC face: the projection Jacobian is discontinuous across that
face, and we do not treat that row as a complete adjoint
validation (\Cref{sec:dpi}).
For a predict-then-optimize loop the practical consequence is a
step-size factor, not a wrong direction: a $10^{-2}$ amplitude
error is absorbed by a well-scaled learning rate and is material
only if the training loss is poorly scaled, as already noted
above.
Strictly active bounds contribute a zero VJP---the correct
Clarke element---so SGD does not try to move those coordinates to
first order.

\paragraph{Speed.}
Forward, backward, and total wall time for one already-constructed
layer.
CVXPYLayers runs on CPU through \texttt{diffcp}/SCS and is shown only
at $n{=}500$; at $n{=}5000$ that path was skipped as too slow
(already 64\,ms total at $n{=}500$).
\cpu{} \diffsolver{} at $n{=}5000$ uses KLU on the unsymmetric
$(n{+}m)$ adjoint (the dense LU shortcut applies only for
$n{+}m\le 800$) and is included.

\begin{table}[t]
\centering
\caption{Differentiable layer on the same factor-model family.
Times in milliseconds for an already-constructed layer, except
CVXPYLayers, whose 100\,ms total at $n{=}500$ \emph{includes}
CVXPY canonicalization and a \texttt{diffcp}/SCS forward solve
(\Cref{sec:exp-rules}); that row is a stack baseline, not a
like-for-like IPM comparison.
Cosine and relative error are versus FD of $\mathbf{1}^\top x$ on
eight interior $q$-coordinates whose one-sided differences agree
($\varepsilon=10^{-6}$).
The $n{=}5000$ CVXPYLayers row is skipped (too slow).}
\label{tab:diff}
\footnotesize
\setlength{\tabcolsep}{3pt}
\small
\begin{tabular}{@{}llccccc@{}}
\toprule
Problem & Layer & Fwd & Bwd & Total & cosine & rel.\ err. \\
\midrule
QP $n{=}500$ & CVXPYLayers & 34.5 & 65.8 & 100.4 & --- & --- \\
QP $n{=}500$ & \cpu{} \diffsolver{} & 2.35 & 0.30 & 2.65 & $1.000$ & $8.8{\times}10^{-4}$ \\
QP $n{=}500$ & \gpu{} \diffsolver{} & 9.48 & 0.45 & 9.97 & $1.000$ & $7.6{\times}10^{-4}$ \\
QP $n{=}5000$ & \cpu{} \diffsolver{} & 135 & 15.9 & 151 & $1.000$ & $3.0{\times}10^{-2}$ \\
QP $n{=}5000$ & \gpu{} \diffsolver{} & 24.1 & 1.02 & 25.0 & $1.000$ & $5.1{\times}10^{-2}$ \\
SOCP $n{=}500$ & \cpu{} \diffsolver{} & 5.47 & 0.43 & 5.91 & $1.000$ & $5.2{\times}10^{-2}$ \\
SOCP $n{=}500$ & \gpu{} \diffsolver{} & 19.8 & 1.10 & 20.9 & $0.999$ & $5.4{\times}10^{-2}$ \\
\bottomrule
\end{tabular}
\end{table}

The GPU backward pass is a fraction of a millisecond after the
forward factorization is warm (0.45\,ms at $n{=}500$, 1.02\,ms at
$n{=}5000$ on the QP).
At $n{=}500$, \cpu{} \diffsolver{} is the fastest already-constructed
layer on this host.
The $38{\times}$ gap versus CVXPYLayers is not a like-for-like IPM
comparison: that row includes canonicalization and a different
forward solver.
That $n{=}500$ row also used \texttt{diffcp}/SCS and skipped
$n{=}5000$ as too slow; for the training-stack question,
\Cref{sec:exp-cvxlayer} replaces it with CVXPYLayers driven by
Clarabel on 1000- and 5000-asset instances.

\subsection{CVXPYLayers as the layer baseline}
\label{sec:exp-cvxlayer}

The comparison of practical interest is not VJP-versus-FD on a
scalar, but whether an in-solver layer reproduces the
\emph{solution and the $\mu$-gradient} of the stack a training
loop would call---CVXPYLayers driven by Clarabel---and at what
wall-clock.
The instances are the $T{=}1$ specialization of the multi-period
factor QP and SOCP of \Cref{sec:exp-portfolio}---the same model
timed everywhere else in this evaluation---with $k{=}50$
and $n\in\{1000,5000\}$ (seed~42):
\[
\begin{aligned}
\min\quad&
-\mu^\top x+\gamma\bigl(y^\top\Omega y+x^\top Dx\bigr)
+c_{\mathrm{tc}}\|x-x_0\|_1\\
\text{s.t.}\quad&
\mathbf{1}^\top x=1,\quad
y=F^\top x,\quad
0\le x,y\le 0.1,
\end{aligned}
\]
with the quadratic risk replaced by an SOC on the SOCP and
transaction-cost rate $c_{\mathrm{tc}}=0.002$
(the $\ell_1$ coefficient; not the homogeneous embedding scalar
$\tauvar$).
Because the turnover term enters as a norm rather than through
explicit slack variables, the canonical $n{=}5000$ QP has $15100$
variables and $25201$ constraints---exactly the $T{=}1$ instance of
\Cref{tab:portfolio-fwd}, so the $n{=}5000$ \cpu{} forward times
here match that table up to protocol.
Tolerances match the rest of the section ($10^{-8}$).
CVXPYLayers and the $n{=}5000$ \cpu{} timings were taken in
processes that never initialize RMM \cite{rmm2024} or the GPU solver
(\Cref{sec:exp-rules} thread budget: one OpenMP/MKL thread).
The GPU rows keep the hybrid AUTO path of the first campaign.
Upstream $\partial L/\partial x$ is a fixed Gaussian vector on
the asset block (seed~123); $d\mu=-\partial L/\partial q$ on
those coordinates is the quantity a return predictor would
use.

\paragraph{Setup.}

\paragraph{Accuracy versus CVXPYLayers.}
\Cref{tab:cvxlayer-acc} takes CVXPYLayers as the oracle.
Every pair is directionally identical (cosine $1.000$ on $x$ and
on $d\mu$).
Asset-weight relative error is $10^{-5}$ on the QPs and
$10^{-6}$--$10^{-10}$ on the SOCPs; the GPU SOCP $x$ matches
CVXPYLayers to $8{\times}10^{-10}$, consistent with both stacks
calling the Clarabel IPM on the same cone product.
The $d\mu$ residual is $8{\times}10^{-5}$--$9.5{\times}10^{-5}$ at
$n{=}1000$ and rises to $8.3{\times}10^{-3}$ on the $n{=}5000$ QP
and $5.6{\times}10^{-4}$ ($5.0{\times}10^{-4}$ on the GPU) on the
$n{=}5000$ SOCP.
That amplitude gap is the same order as the VJP-versus-FD residual
of \Cref{tab:diff} at $n{=}5000$, and CPU/GPU agree with each
other more tightly than either agrees with CVXPYLayers, so the
gap is a stack difference (homogeneous \texttt{diffcp} adjoint
versus the reduced two-block $F^\top$), not a CPU/GPU split.
Objectives agree to about $10^{-8}$.
For a predict-then-optimize loop the implication is the
same as in \Cref{sec:exp-diff}: the descent direction on $\mu$
matches; a $10^{-3}$--$10^{-2}$ amplitude factor is absorbed by
a well-scaled learning rate.

\input{tables/tab_cvxlayer_acc}

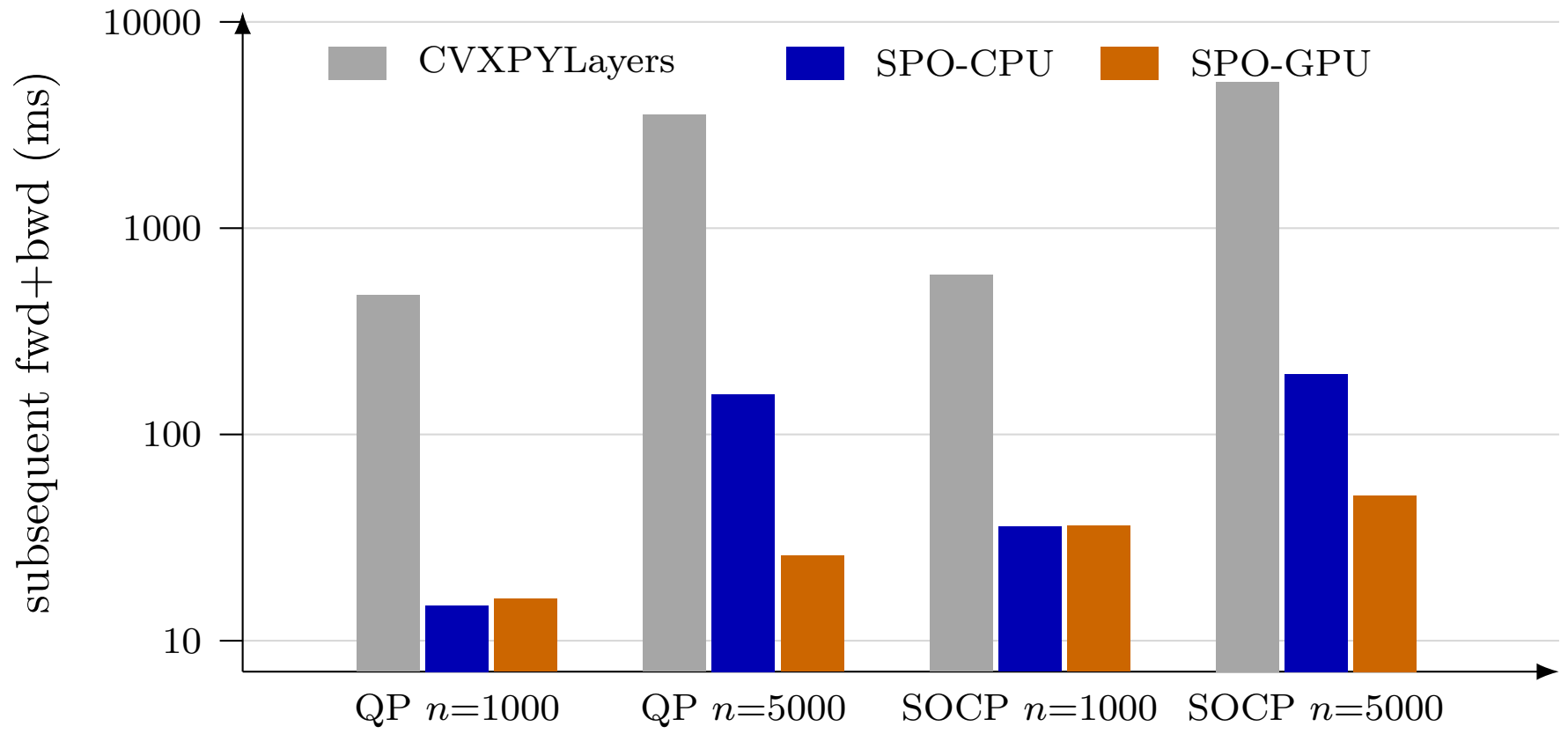
\begin{figure}[t]
\centering
\input{figs/fig_cvxlayer_speedup}
\caption{Subsequent-mean forward+backward time from
\Cref{tab:cvxlayer-perf}, logarithmic vertical axis.
CVXPYLayers includes per-call canonicalization; \diffsolver{}
does not.}
\label{fig:exp-cvxlayer}
\end{figure}

\paragraph{Speed versus CVXPYLayers.}
\Cref{tab:cvxlayer-perf} and \Cref{fig:exp-cvxlayer} report
subsequent-mean forward, backward, and total time.
This is the comparison relevant to a training step: CVXPYLayers
re-canonicalizes and runs a \texttt{diffcp} backward on every
call, while \diffsolver{} is already constructed and factors
only the reduced adjoint.
At $n{=}1000$ the GPU total is $16.0$\,ms against $474$\,ms
($30{\times}$) on the QP and $36.0$\,ms against $594$\,ms
($17{\times}$) on the SOCP.
At $n{=}5000$---a size at which a CPU cone layer has been the
limiting factor---the GPU total is $25.8$\,ms against $3.55$\,s
($137{\times}$) on the QP, and $50.6$\,ms against $5.12$\,s
($101{\times}$) on the SOCP.
The in-solver CPU layer already removes most of the
canonicalization and \texttt{diffcp} overhead
($23{\times}$ / $26{\times}$ at $n{=}5000$); the remaining
order of magnitude is the GPU IPM plus a sub-millisecond to
few-millisecond adjoint (1.03\,ms on the QP, 4.39\,ms on the
SOCP) against a $2.2$--$4.1$\,s CVXPYLayers backward.
Those one-to-two-order totals are the wall-clock evidence for
the claim of \Cref{sec:intro}: a 5000-name subsequent step that
required seconds becomes tens of milliseconds, which is the
difference between a research prototype and a fine-tuning loop
that can run over a trading day.
We do not train that loop here.

\input{tables/tab_cvxlayer_perf}

\subsection{Batch training}
\label{sec:exp-batch}

Same-pattern $n{=}1000$, $k{=}50$ QPs with $B\in\{1,4,8,16\}$,
uniform-batch width and worker count $\min(B,8)$.
Each step is a $q$-update (asset block only), forward, and
backward on already-constructed solvers.
\Cref{tab:batch} is the median of ten such steps.
This is \emph{not} construction time.

\begin{table}[t]
\centering
\caption{\batchdiff{} versus a Python loop of already-constructed
\diffsolver{} objects.
Median fwd+bwd step time in milliseconds.
Hybrid AUTO lets cuDSS place root supernodes on the host;
device-only keeps the factor on the GPU.
An RMM pool of about 19\,GB is reserved at import; per-$B$
increments were not isolated and are not reported.}
\label{tab:batch}
\small
\begin{tabular}{@{}lccc@{}}
\toprule
$B$ & Python loop (ms) & \batchdiff{} hybrid AUTO & device-only \\
\midrule
1 & 15.8 & 17.9 & 20.2 \\
4 & 63.2 & 69.3 & 69.1 \\
8 & 127 & 128 & 127 \\
16 & 260 & 264 & 277 \\
\bottomrule
\end{tabular}
\end{table}

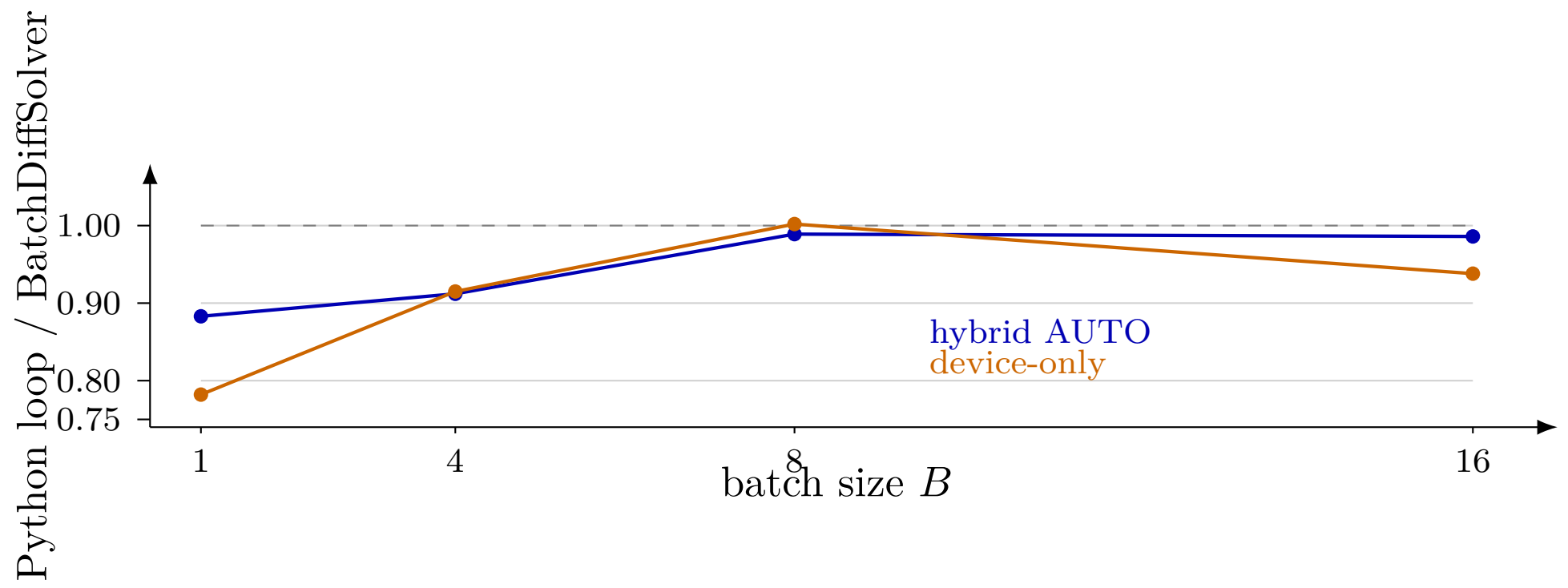
\begin{figure}[t]
\centering
\input{figs/fig_batch_ratio}
\caption{Ratio of median fwd+bwd step times
(Python loop)$/($\batchdiff{}) from \Cref{tab:batch}.
Values below $1$ mean the loop is faster; hybrid AUTO is
2--13\% slower on this idle host.}
\label{fig:exp-batch}
\end{figure}

On this QP, a Python loop of warm \diffsolver{} objects remains
slightly faster than \batchdiff{} at every $B\le 16$
(15.8\,ms $\to$ 260\,ms, a $16.5\times$ factor at $B=16$).
\Cref{fig:exp-batch} plots the ratio (loop)/(batch): hybrid AUTO
lies in $0.88$--$0.99$, i.e.\ 1--12\% slower, which the overlapping
wall-clock curves of an earlier draft concealed.
Each problem already drops the GIL and solves in $\sim$16\,ms; the
uniform-batch scheduler does not amortize below that floor on an
idle host.
We did not time a GIL-contending neighbour (DataLoader or a
concurrent autograd region), and we do not claim a speedup in
that regime.
\batchdiff{} is therefore an architecture
(\Cref{sec:batch})---stable slots, independent iteration counts,
and cuDSS uniform batches---not a measured throughput
improvement, and $B>16$ is not in the protocol.

\subsection{Ablations}
\label{sec:exp-ablate}

Same $n{=}5000$, $k{=}50$ QP as \Cref{tab:portfolio-fwd}.
Equilibration is device-side when enabled, so that row is GPU Ruiz
off.
The elimination-tree cache is not independently switched in this
protocol, so no row is reported for it.

\begin{table}[t]
\centering
\caption{Ablations on the $n{=}5000$, $k{=}50$ QP.
Median of 20 steady solves except the prefactor row (first
solve after construction).}
\label{tab:ablate}
\small
\begin{tabular}{@{}lc@{}}
\toprule
Configuration & Time (ms) \\
\midrule
\gpu{} default & 24.6 \\
Ruiz equilibration off & 24.1 \\
Sequential RHS (no $\mathrm{nrhs}{=}2$) & 27.2 \\
No prefactor at setup (first solve) & 102 \\
\cpu{} default & 132 \\
\cpu{} scalar QDLDL & 122 \\
\bottomrule
\end{tabular}
\end{table}

Prefactorization moves the first-solve stall into setup
(setup 96\,ms with prefactor versus 12\,ms without; the first solve
then drops from 102\,ms to 25\,ms, and subsequent solves stay at
$\sim$24.5\,ms either way).
Disabling Ruiz leaves the iteration count unchanged (13) and is
within noise on this instance.
Forcing scalar QDLDL was \emph{faster} than the default
dispatch on this instance (122\,ms versus 132\,ms): the SIMD
gather/dispatch prologue is a net loss when QDLDL dominates; see
\Cref{sec:cpu-opts}.

\paragraph{Consistency of the $n{=}5000$ \cpu{} numbers.}
Every row in this section instantiates the same $n{=}5000$,
$k{=}50$, $T{=}1$ specialization of the multi-period factor model
of \Cref{sec:exp-portfolio}.
The ranking number is the isolated subsequent mean of
\Cref{tab:portfolio-fwd}: $150$\,ms at seed~42 (11 iterations).
Other protocols on the same instance sit nearby---$132$\,ms as
the median of 20 steady solves in \Cref{tab:ablate},
$131$\,ms first-solve and $133$\,ms steady $q$-update in
\Cref{tab:portfolio-dpp}, $135$\,ms forward in
\Cref{tab:diff}---but are not substituted for that ranking cell.
A three-seed repeat of the Table~7 isolation protocol is in
\Cref{sec:exp-rules}.
An earlier isolation of seed~42 at $133$\,ms is campaign
scatter, not a second official point estimate.
That a $q$-update is no faster than a cold re-solve
is expected, because updating $\dataq$ invalidates the numeric
factor and reruns the full IPM---only the symbolic analysis is
reused.

\subsection{Maros--M\'esz\'aros subset}
\label{sec:exp-maros-pointer}

As a check outside the factor-model family we also timed an
18-problem Maros--M\'esz\'aros QP subset, including CuClarabel on
the same machine; the protocol, subsequent times, and a documented
failure of Clarabel.rs, \cpu{}, and \gpu{} on YAO are reported in
\Cref{sec:exp-maros}.
The subset does not overturn the portfolio ranking and is not a
corpus-wide claim.

%% file: tables/tab_portfolio_first.tex
\begin{table}[t]
\centering
\caption{Factor-model QP/SOCP of \Cref{lst:portfolio-qp,lst:portfolio-socp}: first solver solve time in milliseconds (CVXPY canon, data updates, and solver construction excluded), iterations in parentheses. Every solver sees seed-42 data. Mosek/Gurobi are shown at 1 and 8 threads. Gap/feasibility $10^{-8}$. QP IR off; SOCP IR on (Clarabel.rs/CuClarabel True, \cpu{} AUTO, \gpu{} AUTO2; protocol budget two refinement steps). Four AlmostSolved SOCP cells that returned at a one-step budget were re-timed with two refinement steps and solved (\gpu{}/CuClarabel at $n{=}1000$, $T{=}4$; \cpu{} at $n{=}5000$, $T{=}3$ and $T{=}4$). QP uses a native quadratic objective. A dagger marks a status other than solved/optimal. CuClarabel is post-JIT (first timed solve after a discarded setup+solve).
The \cpu{} column is a later isolated \cpu{} process; the other columns are the original multi-solver campaign.}
\label{tab:portfolio-first}
\footnotesize
\setlength{\tabcolsep}{1.4pt}
\resizebox{\textwidth}{!}{%
\begin{tabular}{@{}l cccccccc@{}}
\toprule
Instance & \makecell{Clarabel\\rs} & \makecell{Clarabel\\CPU} & \makecell{Clarabel\\GPU} & \makecell{CuCla-\\rabel} & \makecell{Mosek\\{\scriptsize 1-thr}} & \makecell{Mosek\\{\scriptsize 8-thr}} & \makecell{Gurobi\\{\scriptsize 1-thr}} & \makecell{Gurobi\\{\scriptsize 8-thr}} \\
\midrule
QP $n{=}200$, $k{=}50$ & 3.22 (9) & 4.78 (9) & 10.2 (9) & 39.1 (9) & 4.96 (10) & 5.19 (10) & 5.03 (13) & 5.43 (13) \\
QP $n{=}500$, $k{=}50$ & 7.33 (10) & 7.07 (10) & 15.2 (10) & 47.9 (10) & 14.9 (14) & 16.0 (14) & 11.2 (14) & 11.7 (14) \\
QP $n{=}500$, $k{=}50$, $T{=}2$ & 26.2 (11) & 20.7 (10) & 19.5 (11) & 57.5 (11) & 44.9 (12) & 42.8 (12) & 44.9 (15) & 35.6 (16) \\
QP $n{=}500$, $k{=}50$, $T{=}3$ & 52.2 (11) & 39.5 (10) & 22.8 (11) & 67.5 (11) & 157 (12) & 156 (12) & 104 (15) & 66.5 (17) \\
QP $n{=}500$, $k{=}50$, $T{=}4$ & 96.1 (12) & 71.2 (11) & 30.4 (12) & 73.3 (12) & 245 (12) & 237 (12) & 172 (16) & 95.4 (17) \\
QP $n{=}1000$, $k{=}50$ & 16.3 (11) & 12.4 (10) & 14.4 (11) & 52.8 (11) & 28.0 (14) & 29.4 (14) & 20.7 (15) & 21.7 (15) \\
QP $n{=}1000$, $k{=}50$, $T{=}2$ & 50.6 (11) & 38.1 (10) & 20.2 (11) & 63.0 (11) & 78.0 (13) & 75.5 (13) & 85.4 (16) & 62.6 (16) \\
QP $n{=}1000$, $k{=}50$, $T{=}3$ & 92.5 (11) & 84.4 (11) & 27.7 (11) & 68.6 (11) & 119 (12) & 118 (12) & 418 (15) & 111 (17) \\
QP $n{=}1000$, $k{=}50$, $T{=}4$ & 166 (12) & 142 (11) & 36.8 (12) & 79.9 (12) & 190 (13) & 179 (13) & 327 (15) & 188 (17) \\
QP $n{=}5000$, $k{=}50$ & 98.0 (13) & 139 (11) & 25.2 (13) & 76.4 (13) & 136 (15) & 139 (15) & 124 (19) & 134 (20) \\
QP $n{=}5000$, $k{=}50$, $T{=}2$ & 379 (16) & 259 (12) & 59.3 (16) & 129 (16) & 404 (14) & 376 (14) & 537 (21) & 407 (22) \\
QP $n{=}5000$, $k{=}50$, $T{=}3$ & 716 (16) & 501 (13) & 89.6 (16) & 166 (16) & 696 (14) & 603 (14) & 1219 (21) & 761 (23) \\
QP $n{=}5000$, $k{=}50$, $T{=}4$ & 1509 (16) & 1236 (13) & 118 (16) & 201 (16) & 1072 (14) & 895 (14) & 2120 (21) & 1655 (23) \\
QP $n{=}2000$, $k{=}50$ & 36.2 (12) & 24.9 (10) & 16.8 (12) & 59.5 (12) & 52.4 (14) & 53.2 (14) & 44.9 (16) & 49.1 (16) \\
QP $n{=}10000$, $k{=}100$ & 1596 (15) & 885 (12) & 59.2 (15) & 131 (15) & 440 (13) & 409 (13) & 697 (19) & 484 (22) \\
\midrule
SOCP $n{=}200$, $k{=}50$ & 12.4 (19) & 8.47 (19) & 29.4 (19) & 140 (19) & 13.5 (11) & 13.2 (11) & 61.8 (12) & 51.8 (14) \\
SOCP $n{=}500$, $k{=}50$ & 25.8 (19) & 15.3 (18) & 27.8 (19) & 143 (19) & 15.4 (15) & 16.3 (15) & 92.3 (17) & 82.5 (17) \\
SOCP $n{=}500$, $k{=}50$, $T{=}2$ & 79.8 (21) & 56.2 (22) & 43.8 (22) & 179 (22) & 54.4 (13) & 56.1 (13) & 253 (15) & 198 (17) \\
SOCP $n{=}500$, $k{=}50$, $T{=}3$ & 138 (19) & 119 (25) & 50.3 (20) & 184 (19) & 170 (14) & 175 (14) & 548 (19) & 326 (19) \\
SOCP $n{=}500$, $k{=}50$, $T{=}4$ & 217 (19) & 189 (24) & 61.7 (20) & 197 (19) & 268 (15) & 266 (15) & 775 (16) & 418 (17) \\
SOCP $n{=}1000$, $k{=}50$ & 52.5 (20) & 32.6 (21) & 31.8 (20) & 162 (20) & 31.2 (16) & 33.3 (16) & 165 (16) & 160 (22) \\
SOCP $n{=}1000$, $k{=}50$, $T{=}2$ & 144 (20) & 120 (26) & 47.5 (20) & 500 (20) & 76.6 (13) & 76.2 (13) & 527 (19) & 419 (22) \\
SOCP $n{=}1000$, $k{=}50$, $T{=}3$ & 261 (20) & 227 (26) & 65.5 (21) & 205 (20) & 147 (15) & 147 (15) & 1090 (20) & 612 (19) \\
SOCP $n{=}1000$, $k{=}50$, $T{=}4$ & 440 (21) & 400 (27) & 95.8 (22) & 262 (23) & 227 (15) & 211 (15) & 1540 (18) & 911 (17) \\
SOCP $n{=}5000$, $k{=}50$ & 251 (19) & 193 (25) & 48.5 (19) & 191 (19) & 187 (23) & 189 (23) & 951 (20) & 913 (23) \\
SOCP $n{=}5000$, $k{=}50$, $T{=}2$ & 912 (24) & 905 (35) & 132 (25) & 364 (24) & 456 (15) & 412 (15) & 3509 (24) & 2668 (30) \\
SOCP $n{=}5000$, $k{=}50$, $T{=}3$ & 1829 (25) & 1630 (33) & 189 (28) & 459 (25) & 827 (16) & 688 (16) & 7278 (24) & 5093 (27) \\
SOCP $n{=}5000$, $k{=}50$, $T{=}4$ & 3610 (25) & 3642 (32) & 242 (25) & 548 (25) & 1210 (15) & 991 (15) & 11695 (26) & 9450 (31) \\
SOCP $n{=}2000$, $k{=}50$ & 90.6 (17) & 77.2 (26) & 28.7 (17) & 142 (17) & 60.6 (17) & 60.9 (17) & 357 (19) & 347 (23) \\
SOCP $n{=}10000$, $k{=}100$ & 978 (21) & 795 (26) & 117 (21) & 336 (21) & 615 (19) & 556 (19) & 4729 (21) & 4497 (45) \\
\bottomrule
\end{tabular}%
}
\end{table}

%% file: tables/tab_portfolio_fwd.tex
\begin{table}[t]
\centering
\caption{Same instances and seed-42 data: subsequent mean of ten sequential solver re-solves on a reused object (CVXPY warm start on the first solve and the ten re-solves; \cpu{}/\gpu{} run each re-solve from the default start, no warm start). Times are solver solve time only. Mosek/Gurobi columns are 1-thread and 8-thread. The four AlmostSolved SOCP cells named in \Cref{tab:portfolio-first} use two refinement steps. A dagger marks a non-optimal status.
The \cpu{} column is a later isolated \cpu{} process (no RMM, no commercial solver); the other columns are the original multi-solver campaign.
Column-to-column rankings at small $n$ are therefore indicative.}
\label{tab:portfolio-fwd}
\footnotesize
\setlength{\tabcolsep}{1.4pt}
\resizebox{\textwidth}{!}{%
\begin{tabular}{@{}l cccccccc@{}}
\toprule
Instance & \makecell{Clarabel\\rs} & \makecell{Clarabel\\CPU} & \makecell{Clarabel\\GPU} & \makecell{CuCla-\\rabel} & \makecell{Mosek\\{\scriptsize 1-thr}} & \makecell{Mosek\\{\scriptsize 8-thr}} & \makecell{Gurobi\\{\scriptsize 1-thr}} & \makecell{Gurobi\\{\scriptsize 8-thr}} \\
\midrule
QP $n{=}200$, $k{=}50$ & 3.21 (9) & 2.81 (9) & 10.2 (9) & 42.7 (9) & 4.43 (10) & 4.64 (10) & 5.16 (13) & 5.32 (13) \\
QP $n{=}500$, $k{=}50$ & 7.30 (10) & 6.90 (10) & 12.4 (10) & 48.0 (10) & 13.9 (14) & 15.6 (14) & 11.0 (14) & 11.6 (14) \\
QP $n{=}500$, $k{=}50$, $T{=}2$ & 25.9 (11) & 21.3 (10) & 17.5 (11) & 56.3 (11) & 42.6 (12) & 43.2 (12) & 43.3 (15) & 36.5 (16) \\
QP $n{=}500$, $k{=}50$, $T{=}3$ & 52.3 (11) & 38.1 (10) & 22.8 (11) & 75.5 (11) & 158 (12) & 155 (12) & 103 (15) & 69.4 (17) \\
QP $n{=}500$, $k{=}50$, $T{=}4$ & 91.4 (12) & 70.3 (11) & 30.2 (12) & 72.8 (12) & 238 (12) & 237 (12) & 171 (16) & 90.1 (17) \\
QP $n{=}1000$, $k{=}50$ & 16.2 (11) & 12.3 (10) & 14.2 (11) & 51.2 (11) & 25.9 (14) & 29.2 (14) & 21.0 (15) & 22.0 (15) \\
QP $n{=}1000$, $k{=}50$, $T{=}2$ & 49.6 (11) & 39.1 (10) & 19.7 (11) & 60.1 (11) & 75.5 (13) & 74.0 (13) & 83.4 (16) & 62.8 (16) \\
QP $n{=}1000$, $k{=}50$, $T{=}3$ & 93.3 (11) & 82.4 (11) & 26.7 (11) & 84.1 (11) & 119 (12) & 117 (12) & 180 (15) & 110 (17) \\
QP $n{=}1000$, $k{=}50$, $T{=}4$ & 166 (12) & 137 (11) & 36.0 (12) & 80.4 (12) & 191 (13) & 177 (13) & 322 (15) & 185 (17) \\
QP $n{=}5000$, $k{=}50$ & 98.2 (13) & 150 (11) & 24.2 (13) & 76.2 (13) & 135 (15) & 135 (15) & 122 (19) & 130 (20) \\
QP $n{=}5000$, $k{=}50$, $T{=}2$ & 381 (16) & 263 (12) & 58.9 (16) & 130 (16) & 400 (14) & 371 (14) & 534 (21) & 405 (22) \\
QP $n{=}5000$, $k{=}50$, $T{=}3$ & 713 (16) & 495 (13) & 89.3 (16) & 165 (16) & 692 (14) & 597 (14) & 1223 (21) & 767 (23) \\
QP $n{=}5000$, $k{=}50$, $T{=}4$ & 1461 (16) & 977 (13) & 118 (16) & 204 (16) & 1065 (14) & 872 (14) & 2083 (21) & 1609 (23) \\
QP $n{=}2000$, $k{=}50$ & 36.6 (12) & 24.9 (10) & 16.4 (12) & 65.8 (12) & 50.7 (14) & 51.4 (14) & 45.7 (16) & 48.7 (16) \\
QP $n{=}10000$, $k{=}100$ & 1585 (15) & 1091 (12) & 58.8 (15) & 136 (15) & 439 (13) & 406 (13) & 698 (19) & 457 (22) \\
\midrule
SOCP $n{=}200$, $k{=}50$ & 12.4 (19) & 8.37 (19) & 28.4 (19) & 147 (19) & 12.5 (11) & 12.5 (11) & 61.7 (12) & 50.8 (14) \\
SOCP $n{=}500$, $k{=}50$ & 25.8 (19) & 15.2 (18) & 26.6 (19) & 157 (19) & 14.7 (15) & 16.0 (15) & 91.9 (17) & 82.5 (17) \\
SOCP $n{=}500$, $k{=}50$, $T{=}2$ & 79.6 (21) & 56.1 (22) & 42.3 (22) & 183 (22) & 54.1 (13) & 55.6 (13) & 254 (15) & 200 (17) \\
SOCP $n{=}500$, $k{=}50$, $T{=}3$ & 133 (19) & 118 (25) & 50.8 (20) & 206 (19) & 173 (14) & 173 (14) & 536 (19) & 326 (19) \\
SOCP $n{=}500$, $k{=}50$, $T{=}4$ & 217 (19) & 191 (24) & 65.4 (20) & 196 (19) & 269 (15) & 270 (15) & 760 (16) & 418 (17) \\
SOCP $n{=}1000$, $k{=}50$ & 52.1 (20) & 32.4 (21) & 31.8 (20) & 156 (20) & 29.5 (16) & 32.2 (16) & 166 (16) & 162 (22) \\
SOCP $n{=}1000$, $k{=}50$, $T{=}2$ & 143 (20) & 126 (26) & 46.8 (20) & 188 (20) & 76.3 (13) & 76.4 (13) & 542 (19) & 419 (22) \\
SOCP $n{=}1000$, $k{=}50$, $T{=}3$ & 262 (20) & 229 (26) & 66.0 (21) & 204 (20) & 147 (15) & 146 (15) & 1062 (20) & 619 (19) \\
SOCP $n{=}1000$, $k{=}50$, $T{=}4$ & 438 (21) & 425 (27) & 86.9 (22) & 256 (23) & 225 (15) & 212 (15) & 1556 (18) & 906 (17) \\
SOCP $n{=}5000$, $k{=}50$ & 251 (19) & 197 (25) & 47.4 (19) & 190 (19) & 184 (23) & 187 (23) & 959 (20) & 905 (23) \\
SOCP $n{=}5000$, $k{=}50$, $T{=}2$ & 911 (24) & 866 (35) & 132 (25) & 364 (24) & 455 (15) & 412 (15) & 3522 (24) & 2672 (30) \\
SOCP $n{=}5000$, $k{=}50$, $T{=}3$ & 1829 (25) & 1695 (33) & 201 (28) & 464 (25) & 814 (16) & 695 (16) & 7278 (24) & 5115 (27) \\
SOCP $n{=}5000$, $k{=}50$, $T{=}4$ & 3617 (25) & 3141 (32) & 244 (25) & 556 (25) & 1211 (15) & 1000 (15) & 11631 (26) & 9614 (31) \\
SOCP $n{=}2000$, $k{=}50$ & 90.1 (17) & 78.4 (26) & 28.3 (17) & 143 (17) & 60.1 (17) & 61.3 (17) & 357 (19) & 350 (23) \\
SOCP $n{=}10000$, $k{=}100$ & 981 (21) & 817 (26) & 117 (21) & 346 (21) & 610 (19) & 562 (19) & 4734 (21) & 4552 (45) \\
\bottomrule
\end{tabular}%
}
\end{table}

%% file: tables/tab_portfolio_gpu_obj.tex
\begin{table}[t]
\centering
\caption{Ten sequential \gpu{} re-solves (each from the default start, no warm start) on the seed-42 factor instances. Objectives are Clarabel \texttt{obj\_val}. A dagger marks a ten-point sample that is not a subset of the closed interval between Mosek and Gurobi; on every such row the overshoot is smaller than $|f_{\mathrm{Mosek}}-f_{\mathrm{Gurobi}}|$, and the sample width itself is also smaller than that gap. The $n{=}100$, $k{=}50$ instances are omitted.}
\label{tab:portfolio-gpu-obj}
\footnotesize
\setlength{\tabcolsep}{2.4pt}
\resizebox{\textwidth}{!}{%
\begin{tabular}{@{}l cccccc@{}}
\toprule
Instance & \gpu{} mean & \gpu{} std & \gpu{} range & Mosek & Gurobi & $|f_{\mathrm{M}}{-}f_{\mathrm{G}}|$ \\
\midrule
QP $n{=}200$, $k{=}50$ & -0.06511696 & $1.46\times 10^{-17}$ & $2.78\times 10^{-17}$ & -0.06511695 & -0.06511696 & $4.30\times 10^{-9}$ \\
QP $n{=}500$, $k{=}50$$^{\dagger}$ & -0.08753556 & $2.56\times 10^{-15}$ & $8.15\times 10^{-15}$ & -0.08753549 & -0.08753556 & $6.97\times 10^{-8}$ \\
QP $n{=}500$, $k{=}50$, $T{=}2$$^{\dagger}$ & -0.17644199 & $9.00\times 10^{-15}$ & $3.04\times 10^{-14}$ & -0.17644191 & -0.17644199 & $8.39\times 10^{-8}$ \\
QP $n{=}500$, $k{=}50$, $T{=}3$ & -0.26544450 & $4.85\times 10^{-14}$ & $1.39\times 10^{-13}$ & -0.26544392 & -0.26544450 & $5.77\times 10^{-7}$ \\
QP $n{=}500$, $k{=}50$, $T{=}4$ & -0.36017053 & $5.27\times 10^{-14}$ & $1.61\times 10^{-13}$ & -0.36017030 & -0.36017053 & $2.28\times 10^{-7}$ \\
QP $n{=}1000$, $k{=}50$$^{\dagger}$ & -0.09891326 & $2.85\times 10^{-15}$ & $9.69\times 10^{-15}$ & -0.09891308 & -0.09891326 & $1.78\times 10^{-7}$ \\
QP $n{=}1000$, $k{=}50$, $T{=}2$ & -0.19627644 & $2.07\times 10^{-14}$ & $6.35\times 10^{-14}$ & -0.19627580 & -0.19627644 & $6.41\times 10^{-7}$ \\
QP $n{=}1000$, $k{=}50$, $T{=}3$ & -0.29481191 & $3.02\times 10^{-14}$ & $9.21\times 10^{-14}$ & -0.29481021 & -0.29481191 & $1.71\times 10^{-6}$ \\
QP $n{=}1000$, $k{=}50$, $T{=}4$ & -0.38923061 & $8.01\times 10^{-15}$ & $2.41\times 10^{-14}$ & -0.38923022 & -0.38923061 & $3.94\times 10^{-7}$ \\
QP $n{=}2000$, $k{=}50$$^{\dagger}$ & -0.10491352 & $1.79\times 10^{-15}$ & $5.87\times 10^{-15}$ & -0.10491264 & -0.10491352 & $8.80\times 10^{-7}$ \\
QP $n{=}5000$, $k{=}50$ & -0.11082304 & $4.74\times 10^{-16}$ & $1.29\times 10^{-15}$ & -0.11082196 & -0.11082304 & $1.09\times 10^{-6}$ \\
QP $n{=}5000$, $k{=}50$, $T{=}2$$^{\dagger}$ & -0.21983046 & $2.95\times 10^{-15}$ & $1.04\times 10^{-14}$ & -0.21982777 & -0.21983045 & $2.68\times 10^{-6}$ \\
QP $n{=}5000$, $k{=}50$, $T{=}3$ & -0.32753716 & $2.49\times 10^{-15}$ & $8.27\times 10^{-15}$ & -0.32753303 & -0.32753717 & $4.13\times 10^{-6}$ \\
QP $n{=}5000$, $k{=}50$, $T{=}4$ & -0.43681468 & $3.30\times 10^{-15}$ & $1.26\times 10^{-14}$ & -0.43680851 & -0.43681468 & $6.17\times 10^{-6}$ \\
QP $n{=}10000$, $k{=}100$ & -0.11150290 & $1.48\times 10^{-16}$ & $4.58\times 10^{-16}$ & -0.11149824 & -0.11150290 & $4.66\times 10^{-6}$ \\
\midrule
SOCP $n{=}200$, $k{=}50$$^{\dagger}$ & 0.03982424 & $7.31\times 10^{-18}$ & $0$ & 0.03982425 & 0.03982424 & $7.19\times 10^{-9}$ \\
SOCP $n{=}500$, $k{=}50$ & -0.01892099 & $3.67\times 10^{-12}$ & $1.20\times 10^{-11}$ & -0.01892095 & -0.01892099 & $4.14\times 10^{-8}$ \\
SOCP $n{=}500$, $k{=}50$, $T{=}2$$^{\dagger}$ & -0.03953206 & $2.43\times 10^{-12}$ & $6.70\times 10^{-12}$ & -0.03953205 & -0.03953206 & $1.40\times 10^{-8}$ \\
SOCP $n{=}500$, $k{=}50$, $T{=}3$$^{\dagger}$ & -0.06038883 & $7.08\times 10^{-10}$ & $2.04\times 10^{-9}$ & -0.06038878 & -0.06038883 & $4.55\times 10^{-8}$ \\
SOCP $n{=}500$, $k{=}50$, $T{=}4$$^{\dagger}$ & -0.09670704 & $9.17\times 10^{-10}$ & $1.80\times 10^{-9}$ & -0.09670701 & -0.09670704 & $3.68\times 10^{-8}$ \\
SOCP $n{=}1000$, $k{=}50$ & -0.04184499 & $2.76\times 10^{-13}$ & $7.98\times 10^{-13}$ & -0.04184492 & -0.04184499 & $7.35\times 10^{-8}$ \\
SOCP $n{=}1000$, $k{=}50$, $T{=}2$$^{\dagger}$ & -0.08264458 & $6.14\times 10^{-10}$ & $1.36\times 10^{-9}$ & -0.08264443 & -0.08264458 & $1.57\times 10^{-7}$ \\
SOCP $n{=}1000$, $k{=}50$, $T{=}3$$^{\dagger}$ & -0.12608273 & $3.72\times 10^{-10}$ & $1.21\times 10^{-9}$ & -0.12608264 & -0.12608273 & $8.99\times 10^{-8}$ \\
SOCP $n{=}1000$, $k{=}50$, $T{=}4$$^{\dagger}$ & -0.16553371 & $5.59\times 10^{-10}$ & $1.28\times 10^{-9}$ & -0.16553353 & -0.16553371 & $1.74\times 10^{-7}$ \\
SOCP $n{=}5000$, $k{=}50$ & -0.07642156 & $4.83\times 10^{-13}$ & $1.69\times 10^{-12}$ & -0.07642044 & -0.07642156 & $1.12\times 10^{-6}$ \\
SOCP $n{=}5000$, $k{=}50$, $T{=}2$$^{\dagger}$ & -0.15398216 & $1.72\times 10^{-10}$ & $5.44\times 10^{-10}$ & -0.15398125 & -0.15398216 & $9.10\times 10^{-7}$ \\
SOCP $n{=}5000$, $k{=}50$, $T{=}3$$^{\dagger}$ & -0.22650769 & $2.72\times 10^{-11}$ & $6.83\times 10^{-11}$ & -0.22650510 & -0.22650769 & $2.59\times 10^{-6}$ \\
\bottomrule
\end{tabular}%
}
\end{table}

%% file: tables/tab_n20000.tex
\begin{table}[t]
\centering
\caption{Companion $n{=}20000$, $k{=}50$, $T{=}1$ campaign
on the listing generator (seed~42), \emph{not} part of
\Cref{tab:portfolio-first,tab:portfolio-fwd}.
Subsequent mean of ten cold re-solves (solver solve only;
construction excluded).
\cpu{} is isolated (no GPU/RMM in-process).
CuClarabel and 8-thread commercial columns were not run at
this size.
Iterations in parentheses.
The same instances supply the $n{=}20000$ rows of
\Cref{tab:portfolio-acc,tab:residuals}.}
\label{tab:n20000}
\small
\begin{tabular}{@{}lccccc@{}}
\toprule
Instance & Clarabel.rs & \cpu{} & \gpu{} & Mosek & Gurobi \\
\midrule
QP subsequent (ms)
  & $2528$ (18) & $1085$ (13) & $92.7$ (18) & $661$ (16) & $721$ (20) \\
SOCP subsequent (ms)
  & $1593$ (25) & $975$ (29) & $191$ (25) & $887$ (24) & $5180$ (23) \\
\bottomrule
\end{tabular}
\end{table}

%% file: tables/tab_residuals.tex
\begin{table}[t]
\centering
\caption{Clarabel-form residuals of \eqref{eq:kkt} on the
$T{=}1$ factor-model QP of the listings (seed~42, $k{=}50$),
$\mathrm{pres}=\|Ax+s-b\|_2/(1+\|b\|_2)$,
$\mathrm{dres}=\|Px+q+A^\top z\|_2/(1+\|q\|_2)$,
$\mathrm{gap}=|s^\top z|/(1+|f|)$.
Clarabel.rs, \cpu{}, and \gpu{} meet the $10^{-8}$ protocol
by one to five orders of magnitude.
A Mosek Task on the same $(P,q,A,b)$ reports its own
relative gap and \texttt{solsta} (native merit);
we do not reconstruct Clarabel-form $(s,z)$ from Mosek duals.
\Cref{tab:portfolio-acc} compares objectives only.}
\label{tab:residuals}
\small
\begin{tabular}{@{}llccc@{}}
\toprule
$n$ & Solver & pres & dres & gap \\
\midrule
$5000$ & Clarabel.rs
  & $1.2{\times}10^{-10}$ & $8.6{\times}10^{-12}$ & $5.5{\times}10^{-9}$ \\
$5000$ & \cpu{}
  & $1.6{\times}10^{-9}$ & $8.5{\times}10^{-13}$ & $1.1{\times}10^{-9}$ \\
$5000$ & \gpu{}
  & $1.2{\times}10^{-10}$ & $8.6{\times}10^{-12}$ & $5.5{\times}10^{-9}$ \\
$20000$ & Clarabel.rs
  & $5.3{\times}10^{-11}$ & $9.7{\times}10^{-12}$ & $7.2{\times}10^{-9}$ \\
$20000$ & \cpu{}
  & $1.7{\times}10^{-9}$ & $6.9{\times}10^{-12}$ & $5.4{\times}10^{-9}$ \\
$20000$ & \gpu{}
  & $5.2{\times}10^{-11}$ & $9.7{\times}10^{-12}$ & $7.2{\times}10^{-9}$ \\
$5000$ & Mosek (native)
  & $1.3{\times}10^{-9\,\ast}$ & --- & $7.5{\times}10^{-10}$ \\
$20000$ & Mosek (native)
  & $3.0{\times}10^{-9\,\ast}$ & --- & $2.0{\times}10^{-8}$ \\
\bottomrule
\end{tabular}
\\[0.4em]
{\footnotesize $^\ast$Mosek \texttt{getpviolcon} infinity-norm
(native constraint violation), not Clarabel $\mathrm{pres}$.
Both Tasks report \texttt{solsta.optimal}.}
\end{table}

%% file: tables/tab_warmstart_equil.tex
\begin{tabular}{@{}llllccccc@{}}
\toprule
 & & & & \multicolumn{5}{c}{relative $\mu$ perturbation} \\
\cmidrule(l){5-9}
$n$ & Family & Equil. & Method & $10^{-4}$ & $10^{-3}$ & $10^{-2}$ & $10^{-1}$ & $1$ \\
\midrule
1000 & QP & on & \cpu{} & 10\,$\to$\,5 & 10\,$\to$\,7 & 10\,$\to$\,8 & 10\,$\to$\,9 & 11\,$\to$\,11 \\
 &  &  & \gpu{} & 11\,$\to$\,6 & 11\,$\to$\,8 & 11\,$\to$\,9 & 11\,$\to$\,10 & 13\,$\to$\,13 \\
 &  & off & \cpu{} & 11\,$\to$\,7 & 11\,$\to$\,9 & 11\,$\to$\,10 & 11\,$\to$\,11 & 11\,$\to$\,11 \\
 &  &  & \gpu{} & 11\,$\to$\,7 & 11\,$\to$\,9 & 12\,$\to$\,10 & 12\,$\to$\,11 & 11\,$\to$\,11 \\
\addlinespace
 & SOCP & on & \cpu{} & 21\,$\to$\,13 & 21\,$\to$\,14 & 20\,$\to$\,16 & 20\,$\to$\,20 & 21\,$\to$\,24 \\
 &  &  & \gpu{} & 20\,$\to$\,15 & 20\,$\to$\,15 & 19\,$\to$\,22 & 17\,$\to$\,22 & 18\,$\to$\,19 \\
 &  & off & \cpu{} & 14\,$\to$\,20 & 14\,$\to$\,24 & 14\,$\to$\,24 & 14\,$\to$\,19 & 14\,$\to$\,20 \\
 &  &  & \gpu{} & 16\,$\to$\,17 & 16\,$\to$\,22 & 16\,$\to$\,20 & 16\,$\to$\,16 & 17\,$\to$\,20 \\
\midrule
5000 & QP & on & \cpu{} & 11\,$\to$\,6 & 11\,$\to$\,7 & 11\,$\to$\,9 & 11\,$\to$\,11 & 11\,$\to$\,15 \\
 &  &  & \gpu{} & 13\,$\to$\,7 & 13\,$\to$\,8 & 14\,$\to$\,10 & 14\,$\to$\,13 & 18\,$\to$\,18 \\
 &  & off & \cpu{} & 13\,$\to$\,8 & 13\,$\to$\,10 & 13\,$\to$\,11 & 13\,$\to$\,13 & 13\,$\to$\,14 \\
 &  &  & \gpu{} & 13\,$\to$\,8 & 13\,$\to$\,10 & 13\,$\to$\,11 & 13\,$\to$\,13 & 13\,$\to$\,14 \\
\addlinespace
 & SOCP & on & \cpu{} & 25\,$\to$\,14 & 25\,$\to$\,18 & 25\,$\to$\,20 & 25\,$\to$\,25 & 25\,$\to$\,30 \\
 &  &  & \gpu{} & 19\,$\to$\,16 & 19\,$\to$\,24 & 19\,$\to$\,23 & 20\,$\to$\,21 & 20\,$\to$\,25 \\
 &  & off & \cpu{} & 18\,$\to$\,17 & 18\,$\to$\,22 & 18\,$\to$\,23 & 18\,$\to$\,24 & 18\,$\to$\,23 \\
 &  &  & \gpu{} & 19\,$\to$\,16 & 19\,$\to$\,20 & 19\,$\to$\,22 & 18\,$\to$\,20 & 23\,$\to$\,26 \\
\bottomrule
\end{tabular}

%% file: tables/tab_eps_sweep.tex
\begin{table}[t]
\centering
\caption{$\varepsilon$-convergence sweep for the VJP check of
\Cref{tab:diff}: relative error
$\|\hat g-g_{\mathrm{FD}}\|_2/\|g_{\mathrm{FD}}\|_2$ of the layer VJP
against a central finite difference of $\mathbf{1}^\top x$ (eight
smooth interior coordinates, kink filter on), on the $n{=}500$ and
$n{=}5000$ QP rows for both backends.
The residual is essentially flat over the four decades swept on all
four configurations, i.e.\ it is an $O(\delta)$ solver-accuracy floor
set by the IPM tolerance, not $O(\varepsilon^2)$ truncation; no
truncation rise at $\varepsilon{=}10^{-4}$ and no
$O(\delta/\varepsilon)$ noise rise at $\varepsilon{=}10^{-7}$ is
resolved above run-to-run noise.
Cosine is $\ge 0.9997$ everywhere.
The $\varepsilon{=}10^{-6}$ used in \Cref{tab:diff} sits inside the
flat regime in all four configurations.}
\label{tab:eps-sweep}
\small
\setlength{\tabcolsep}{6pt}
\begin{tabular}{@{}lrrrr@{}}
\toprule
 & \multicolumn{2}{c}{$n{=}500$ QP} & \multicolumn{2}{c}{$n{=}5000$ QP} \\
\cmidrule(lr){2-3}\cmidrule(lr){4-5}
$\varepsilon$ & GPU & CPU & GPU & CPU \\
\midrule
$10^{-4}$ & $7.85{\times}10^{-4}$ & $9.24{\times}10^{-4}$ & $5.16{\times}10^{-2}$ & $3.00{\times}10^{-2}$ \\
$10^{-5}$ & $7.81{\times}10^{-4}$ & $8.85{\times}10^{-4}$ & $5.16{\times}10^{-2}$ & $3.00{\times}10^{-2}$ \\
$10^{-6}$ & $7.81{\times}10^{-4}$ & $8.85{\times}10^{-4}$ & $5.16{\times}10^{-2}$ & $3.00{\times}10^{-2}$ \\
$10^{-7}$ & $7.80{\times}10^{-4}$ & $8.85{\times}10^{-4}$ & $5.16{\times}10^{-2}$ & $2.95{\times}10^{-2}$ \\
\bottomrule
\end{tabular}
\end{table}

%% file: tables/tab_vjp_ratios.tex
\begin{table}[t]
\centering
\caption{Per-coordinate VJP/FD ratios
$\hat g_i/g_{\mathrm{FD},i}$ on interior QP coordinates
($\varepsilon=10^{-6}$, $L=\mathbf{1}^\top x$, GPU, 24 samples)
and finite-difference checks of $\partial L/\partial P$,
$\partial L/\partial A$ on a $n{=}200$ QP (20 nonzero entries).
The median ratio is near $1$; the interquartile range is a few
percent, but the tails reach $O(1)$
($n{=}5000$ min $0.017$).
$P$/$A$ pass rates are among entries with $|g_{\mathrm{FD}}|>10^{-10}$.}
\label{tab:vjp-ratios}
\small
\begin{tabular}{@{}lcccc@{}}
\toprule
Check & $n$ & median [q25, q75] & within 10\% & notes \\
\midrule
$\partial(\mathbf{1}^\top x)/\partial q_i$
  & $500$ & $1.012$ $[0.983,1.039]$ & $21/24$ & min $0.81$, max $2.33$ \\
$\partial(\mathbf{1}^\top x)/\partial q_i$
  & $5000$ & $1.008$ $[0.987,1.054]$ & $18/24$ & min $0.017$, max $1.22$ \\
$\partial L/\partial P_{rc}$
  & $200$ & rel.\ med.\ $4.9{\times}10^{-2}$ & $9/20$ at $10^{-2}$ & $11/20$ at $10^{-1}$ \\
$\partial L/\partial A_{rc}$
  & $200$ & rel.\ med.\ $1.6{\times}10^{-3}$ & $12/20$ at $10^{-2}$ & $12/20$ at $10^{-1}$ \\
\bottomrule
\end{tabular}
\end{table}

%% file: tables/tab_cvxlayer_acc.tex
\begin{table}[t]
\centering
\caption{Layer accuracy with CVXPYLayers (driven by Clarabel)
as the benchmark, on the $T{=}1$ specialization of the multi-period
factor QP/SOCP of \Cref{sec:exp-cvxlayer} (seed~42).
Relative error and cosine are on asset weights $x\in\R^n$ and on
$d\mu=-\partial L/\partial q$ restricted to those coordinates.
Objectives agree to about $10^{-8}$ and are omitted.}
\label{tab:cvxlayer-acc}
\footnotesize
\setlength{\tabcolsep}{4pt}
\begin{tabular}{@{}llcccc@{}}
\toprule
Problem & Layer & $x$ rel.\ & $\cos x$ & $d\mu$ rel.\ & $\cos d\mu$ \\
\midrule
QP $n{=}1000$ & \cpu{} \diffsolver{}
  & $4.7{\times}10^{-5}$ & $1.000$ & $9.4{\times}10^{-5}$ & $1.000$ \\
QP $n{=}1000$ & \gpu{} \diffsolver{}
  & $3.2{\times}10^{-5}$ & $1.000$ & $9.5{\times}10^{-5}$ & $1.000$ \\
QP $n{=}5000$ & \cpu{} \diffsolver{}
  & $6.2{\times}10^{-5}$ & $1.000$ & $8.3{\times}10^{-3}$ & $1.000$ \\
QP $n{=}5000$ & \gpu{} \diffsolver{}
  & $1.2{\times}10^{-4}$ & $1.000$ & $8.3{\times}10^{-3}$ & $1.000$ \\
SOCP $n{=}1000$ & \cpu{} \diffsolver{}
  & $5.7{\times}10^{-6}$ & $1.000$ & $7.9{\times}10^{-5}$ & $1.000$ \\
SOCP $n{=}1000$ & \gpu{} \diffsolver{}
  & $6{\times}10^{-9}$ & $1.000$ & $7.9{\times}10^{-5}$ & $1.000$ \\
SOCP $n{=}5000$ & \cpu{} \diffsolver{}
  & $2.5{\times}10^{-5}$ & $1.000$ & $5.6{\times}10^{-4}$ & $1.000$ \\
SOCP $n{=}5000$ & \gpu{} \diffsolver{}
  & $8{\times}10^{-10}$ & $1.000$ & $5.0{\times}10^{-4}$ & $1.000$ \\
\bottomrule
\end{tabular}
\end{table}

%% file: tables/tab_cvxlayer_perf.tex
\begin{table}[t]
\centering
\caption{Subsequent-mean layer time versus CVXPYLayers
(driven by Clarabel) on the instances of
\Cref{tab:cvxlayer-acc}.
Times in milliseconds: arithmetic mean of ten isolated repeats for
CVXPYLayers; five repeats for the \diffsolver{} rows.
CVXPYLayers includes per-call CVXPY canonicalization and a
\texttt{diffcp} backward; \diffsolver{} is already constructed.
The last column is (CVXPYLayers total)$/$(\diffsolver{} total).}
\label{tab:cvxlayer-perf}
\footnotesize
\setlength{\tabcolsep}{3.5pt}
\begin{tabular}{@{}llcccc@{}}
\toprule
Problem & Layer & Fwd & Bwd & Total & vs.\ CVX \\
\midrule
QP $n{=}1000$ & CVXPYLayers & 160 & 314 & 474 & --- \\
QP $n{=}1000$ & \cpu{} \diffsolver{} & 12.9 & 1.95 & 14.8 & $32{\times}$ \\
QP $n{=}1000$ & \gpu{} \diffsolver{} & 15.4 & 0.66 & 16.0 & $30{\times}$ \\
QP $n{=}5000$ & CVXPYLayers & 1352 & 2197 & 3549 & --- \\
QP $n{=}5000$ & \cpu{} \diffsolver{} & 141 & 16.0 & 157 & $23{\times}$ \\
QP $n{=}5000$ & \gpu{} \diffsolver{} & 24.8 & 1.03 & 25.8 & $137{\times}$ \\
SOCP $n{=}1000$ & CVXPYLayers & 113 & 481 & 594 & --- \\
SOCP $n{=}1000$ & \cpu{} \diffsolver{} & 33.2 & 2.57 & 35.7 & $17{\times}$ \\
SOCP $n{=}1000$ & \gpu{} \diffsolver{} & 34.1 & 1.95 & 36.0 & $17{\times}$ \\
SOCP $n{=}5000$ & CVXPYLayers & 1021 & 4100 & 5121 & --- \\
SOCP $n{=}5000$ & \cpu{} \diffsolver{} & 184 & 12.6 & 197 & $26{\times}$ \\
SOCP $n{=}5000$ & \gpu{} \diffsolver{} & 46.2 & 4.39 & 50.6 & $101{\times}$ \\
\bottomrule
\end{tabular}
\end{table}

%% file: figs/fig_cvxlayer_speedup.tex
\begin{tikzpicture}[x=2.15cm,y=1.55cm]
  \def\xmin{0.25}
  \def\xmax{4.85}
  \def\ymin{0.85}
  \def\ymax{4.05}
  \draw[gray!30] (\xmin,1) -- (\xmax,1);
  \draw[gray!30] (\xmin,2) -- (\xmax,2);
  \draw[gray!30] (\xmin,3) -- (\xmax,3);
  \draw[gray!30] (\xmin,4) -- (\xmax,4);
  \draw[-{Latex}] (\xmin,\ymin) -- (\xmax,\ymin);
  \draw[-{Latex}] (\xmin,\ymin) -- (\xmin,\ymax);
  \foreach \y/\lab in {1/10,2/100,3/1000,4/10000}
    \draw (\xmin,\y) -- ++(-0.08,0) node[left,font=\scriptsize] {\lab};
  \foreach \x/\lab in {1/{QP $n{=}1000$},2/{QP $n{=}5000$},3/{SOCP $n{=}1000$},4/{SOCP $n{=}5000$}}
    \node[below,font=\scriptsize,align=center] at (\x,\ymin-0.02) {\lab};
  \foreach \x/\hcvx/\hcpu/\hgpu in {
      1/2.676/1.170/1.204,
      2/3.550/2.195/1.412,
      3/2.774/1.553/1.556,
      4/3.709/2.293/1.704} {
    \fill[black!35] (\x-0.35,\ymin) rectangle (\x-0.13,\hcvx);
    \fill[blue!70!black] (\x-0.11,\ymin) rectangle (\x+0.11,\hcpu);
    \fill[orange!80!black] (\x+0.13,\ymin) rectangle (\x+0.35,\hgpu);
  }
  \node[font=\small,rotate=90] at (\xmin-0.72,2.45) {subsequent fwd+bwd (ms)};
  \fill[black!35] (0.55,3.72) rectangle (0.75,3.88);
  \node[anchor=west,font=\scriptsize] at (0.80,3.80) {CVXPYLayers};
  \fill[blue!70!black] (2.15,3.72) rectangle (2.35,3.88);
  \node[anchor=west,font=\scriptsize] at (2.40,3.80) {\cpu{}};
  \fill[orange!80!black] (3.25,3.72) rectangle (3.45,3.88);
  \node[anchor=west,font=\scriptsize] at (3.50,3.80) {\gpu{}};
\end{tikzpicture}

%% file: figs/fig_batch_ratio.tex
\begin{tikzpicture}[x=0.70cm,y=6.4cm]
  \def\xmin{0.4}
  \def\xmax{17.0}
  \def\ymin{0.74}
  \draw[gray!35] (1,0.80) -- (16,0.80);
  \draw[gray!35] (1,0.90) -- (16,0.90);
  \draw[gray!35] (1,1.00) -- (16,1.00);
  \draw[-{Latex}] (\xmin,\ymin) -- (\xmax,\ymin);
  \draw[-{Latex}] (\xmin,\ymin) -- (\xmin,1.08);
  \foreach \x in {1,4,8,16}
    \draw (\x,\ymin) -- ++(0,-0.008) node[below,font=\scriptsize] {\x};
  \foreach \y/\lab in {0.75/0.75,0.80/0.80,0.90/0.90,1.00/1.00}
    \draw (\xmin,\y) -- ++(-0.15,0) node[left,font=\scriptsize] {\lab};
  \draw[dashed,gray] (1,1.0) -- (16,1.0);
  \draw[blue!70!black,thick]
    (1,0.883) -- (4,0.912) -- (8,0.989) -- (16,0.986);
  \fill[blue!70!black]
    (1,0.883) circle (1.7pt)
    (4,0.912) circle (1.7pt)
    (8,0.989) circle (1.7pt)
    (16,0.986) circle (1.7pt);
  \draw[orange!80!black,thick]
    (1,0.782) -- (4,0.915) -- (8,1.002) -- (16,0.938);
  \fill[orange!80!black]
    (1,0.782) circle (1.7pt)
    (4,0.915) circle (1.7pt)
    (8,1.002) circle (1.7pt)
    (16,0.938) circle (1.7pt);
  \node[font=\small] at (8.5,\ymin-0.07) {batch size $B$};
  \node[font=\small,rotate=90] at (\xmin-1.35,0.91) {Python loop $/$ BatchDiffSolver};
  \node[blue!70!black,font=\scriptsize,anchor=west] at (9.4,0.86) {hybrid AUTO};
  \node[orange!80!black,font=\scriptsize,anchor=west] at (9.4,0.82) {device-only};
\end{tikzpicture}

%% file: sections/11_discussion.tex
\section{Discussion}
\label{sec:discussion}

The significance of the stack is not a new embedding, but a single
IPM that meets the three requirements of \Cref{sec:intro} on
parametric high-accuracy QCPs: tight tolerances, reused symbolic
factorization, and an in-solver VJP.
The measurements below locate the regimes in which that claim is
supported, and those in which it is not.

\paragraph{Measured crossover.}
The crossover structure of \Cref{tab:portfolio-fwd} is consistent:
\cpu{} and Clarabel.rs lead at small $n$, \gpu{} leads the
whole field from $n{=}2000$ on this factor-model QP family, and
its advantage widens with $n$.
The $n{=}1000$ $T{=}1$ QP cell (isolated \cpu{} 12.3\,ms versus
original-campaign \gpu{} 14.2\,ms) is too close, and too mixed
across campaigns, to treat as a crossover.
Raising $T$ at fixed $n$ moves the GPU lead earlier
(\gpu{} already leads at $n{=}500$, $T{=}2$).
The \cpu{} column of
\Cref{tab:portfolio-first,tab:portfolio-fwd} is a later isolated
\cpu{} process (no commercial solver, no RMM); the other
columns are the original multi-solver campaign.
First execution reverses some rankings: \gpu{}'s tabled first solve
is the first solver call after construction, while Clarabel.rs
incurs a symbolic-setup cost that is non-monotone in $n$
(\Cref{tab:portfolio-first}).
\batchdiff{} did not outperform a Python loop of warm \diffsolver{}
objects at $B\le 16$ on that QP: each slot already solves in
$\sim$15\,ms, so the uniform-batch scheduler has little residual
work to amortize.
It is an architecture (stable slots and cuDSS uniform batches),
not an additional performance claim.

\paragraph{When to use which backend.}
\gpu{} is preferable when the KKT is large and sparse enough that
cuDSS occupies most of the iteration, which is the factor-model
portfolio regime of \Cref{sec:portfolio-layer} at a few thousand
assets.
\cpu{} is preferable when the model canonicalizes to many tiny
nonsymmetric cones and a small KKT---in particular a
high-dimensional power cone reduced to $O(\cdim_1)$
three-dimensional power cones.
In that shape, QDLDL on the host is faster than cuDSS plus
per-cone GPU launches.
A practical heuristic is: when the model has at least eight 3D
power cones and the KKT is small, use the CPU path; mixed cone
families with large unsplit SOCs or mixed-dimension PSD blocks
favour the GPU four-tier kernels and the chordal path.

\paragraph{When to warm-start.}
Smoothing should be left on for parametric QPs: on the
$n{=}1000$ family the one-shot saving shrinks smoothly with the jump
($10\to 5$ at $10^{-4}$, gone at a jump of size $1$) and a
twenty-step rebalance at $10^{-2}$ is a stable 3--4 iteration
reduction (\Cref{sec:exp-warm}).
On SOCPs the same control is useful only for small jumps.
On the production maps it reduces the iteration count at
$10^{-4}$--$10^{-3}$ on both backends and in the CPU sequential
$10^{-2}$ walk (the GPU floor, cap off, is already slightly
above cold at one-shot $10^{-4}$, \Cref{tab:soc-cap}); it does
not reduce it for $O(1)$ $\mu$ jumps, and the GPU sequential
$10^{-2}$ SOCP chain is slightly slower than cold.
That GPU cell is a path effect of Theorem~4.1's lift of
$s_0$ with four-tier SOC scaling (\Cref{tab:soc-cap}):
turning the cap off raises sequential $R$ from $1.20$ to
$1.34$, so the cap is neither a sufficient nor a necessary
cause.
Smoothing is not appropriate for a single isolated solve.

\paragraph{End-to-end training.}
A two-stage predict-then-allocate workflow never includes the QCP
in the training graph; an end-to-end loop does, but only if the
optimizer returns $\partial L/\partial\mu$
(\Cref{sec:intro,sec:portfolio-layer}).
\Cref{tab:cvxlayer-acc,tab:cvxlayer-perf} provide the layer-level
evidence: on a 5000-name instance the in-solver GPU layer
matches the CVXPYLayers (Clarabel) $x$ and $d\mu$ and reduces a
subsequent step from seconds to tens of milliseconds
($137{\times}$ / $101{\times}$).
That wall-clock reduction is the reason the stack can support
production fine-tuning; it is not a measured Sharpe or
transaction-cost improvement, and we do not treat it as one.
The remaining work is the learning curve itself---a network, a
decision loss, and a backtest---which this paper does not run.
The forward solver is used in private deployments; we do not
treat that as a measured claim of this paper.

\paragraph{Batching.}
\batchdiff{} binds same-pattern QCPs to stable slots and gathers
cuDSS phases into uniform batches.
An isolated Python loop of \diffsolver{} objects already drops the
GIL per call and is slightly faster on the measured QP
(\Cref{sec:exp-batch}); a GIL-contending neighbour was not timed,
and we do not claim a speedup in that regime.
The design requires memory linear in the number of stable slots and
requires a shared sparsity pattern; the cuDSS shard width is a
tunable independent of the problem and worker counts, and oversized
shards can saturate the device and lose to smaller ones.
It is not appropriate for one sequential backtest trajectory: a
single \diffsolver{} with $q$-updates and smoothing should be used
instead.

\paragraph{Precision.}
As recorded in \Cref{sec:gpu-delta}, we keep binary64 where
CuClarabel offers mixed precision \cite{chen2025cuclarabel}: the IPM
Jacobian becomes increasingly ill-conditioned, iterative refinement
then dominates, and a differentiable layer inherits any forward
inconsistency.
A future mixed-precision path, if any, should be confined to
refinement and gated like the latched AUTO policy of
\Cref{sec:gpu-cudss}.

\paragraph{Limitations.}
The portfolio tables are one factor-model family and one seed;
the CVXPYLayers comparison of \Cref{sec:exp-cvxlayer} uses the
same family at $T{=}1$, still one seed.
That remains the primary empirical limitation.
\Cref{sec:exp-maros} adds an 18-problem Maros--M\'esz\'aros subset on the same
machine, not the 138-problem corpus; CuClarabel is timed only on
that subset.
The Maros subset covers 18 of the 138 problems; the subset does not
replace a full corpus or a CBLIB pass.
A second limitation is statistical: every portfolio and layer number
is one seed (42) except the three-seed isolated $n{=}5000$ QP of
\Cref{sec:exp-rules} (seeds $42,7,99$).
The measured run-to-run scatter of about 7\%
is of the same order as the 1--13\% batch gap
and some ablation rows, so percentage-level conclusions on one seed
should be read as indicative, not decisive.
A third is baseline coverage: \texttt{diffqcp} is not timed
(\Cref{sec:exp-rules}).
\Cref{tab:residuals} now records Clarabel-form residuals on the
same standard form; \Cref{tab:portfolio-acc} remains an objective
comparison, not a residual oracle.
\cpu{} has no multi-problem batch.
Device-side updates of $P$ and $A$ from device tensors without a host
permutation are not complete.
Large generalized power cones with nonuniform exponents remain more
reliable when split into 3D power cones than when the native
high-dimensional IPM is used---a limitation shared with Clarabel.rs,
which also omits primal-dual scaling on those cones.
Chordal clique merge can refuse decomposition rather than launch an
incorrect batch.
Hybrid cuDSS execute applies to shard sizes of at most eight;
larger shards run device-only.

%% file: sections/12_conclusion.tex
\section{Conclusion}
\label{sec:conclusion}

\paragraph{Results.}
We described Sequential Portfolio Optimizer (\spo{}):
\cpu{} and \gpu{}, CPU and CUDA interior-point solvers
for quadratic cone programs in the lineage of Clarabel and
CuClarabel, with an in-solver differentiable layer and, on the GPU,
a multi-problem batch architecture.
The Mehrotra homogeneous embedding is unchanged
(\Cref{ass:standing}).
The original mathematical content is the reduced two-block VJP of
\Cref{prop:vjp}, the documented smoothing variants of
\Cref{def:smoothing-impl}, and the problem distance of
\Cref{def:dkkt,prop:dkkt-invar}.
The systems work is a SIMD-aware QDLDL plus a KLU adjoint and pooled
workspaces on the CPU; a device-resident setup path, mixed-dimension
PSD cones, four-tier SOC kernels, and a cuDSS reuse protocol on the
GPU; and \batchdiff{}, which binds same-pattern QCPs to stable slots
and gathers cuDSS phases into uniform batches.
A forward optimizer and a \diffsolver{} share one IPM, so modelling
rebalancing and training do not share a binding.
On the machine of \Cref{sec:experiments} the crossover is that of
\Cref{tab:portfolio-fwd}: \cpu{} and Clarabel.rs are fastest at
small $n$, \gpu{} is fastest from $n{=}2000$ on this factor-model
QP family, and the gap widens with $n$.
The $n{=}1000$ $T{=}1$ comparison mixes an isolated \cpu{}
column with the original multi-solver campaign and is not a
hard crossover.
On sequential $\mu$-rebalancing at relative step $10^{-2}$ on the
$n{=}1000$ family, smoothing reduces mean QP iterations from
$10.0$ to $8.0$ on \cpu{} ($R=0.80$) and from $11.0$ to $9.0$ on
\gpu{} ($R=0.82$), and SOCP iterations from $20.0$ to $15.3$ on
\cpu{} ($R=0.77$).
The in-solver VJP matches finite differences on interior QP
coordinates at cosine $1.000$ with relative error
$8{\times}10^{-4}$ at $n{=}500$ and $3$--$5{\times}10^{-2}$ at
$n{=}5000$ (\Cref{tab:diff});
the $\partial L/\partial P$ check is weaker
($9/20$ entries within $10^{-2}$ on a small QP,
\Cref{tab:vjp-ratios}) and is not claimed at the same
level.
Against CVXPYLayers with Clarabel as the training-stack oracle,
the same GPU layer matches $x$ and $d\mu$ (cosine $1.000$) and
is $137{\times}$ / $101{\times}$ faster on subsequent
forward+backward at $n{=}5000$ (\Cref{tab:cvxlayer-perf}).

\paragraph{Limitations.}
Those rankings are not claimed as corpus-wide;
\Cref{sec:exp-maros} times an 18-problem Maros--M\'esz\'aros QP
subset on the same machine, and CBLIB \cite{cbf2016} and Netlib
remain outside the protocol.
Portfolio timings use one factor-model family and one seed.
The documented failure mode is the GPU SOCP warm-start cell: the
same $10^{-2}$ chain lies above cold ($18.0\to 21.7$, $R=1.20$).
That cell is a GPU-path effect of the lift of $s_0$ in
\Cref{thm:cgj} with four-tier sparse SOC scaling.
The $2{\times}2$ of \Cref{tab:soc-cap} keeps CPU at $R<1$
(floor $0.77$, cap $0.85$) and GPU at $R>1$ (floor $1.34$,
cap $1.20$), so the cap is neither sufficient nor necessary.
That cell is a limitation, not a claimed speedup.
The GPU SOCP VJP check is weaker (cosine $0.999$, relative error
$5.4{\times}10^{-2}$ on the $n{=}500$ row of \Cref{tab:diff}), and the
layer returns a single Clarke element at an active bound rather than
the full subdifferential.
Relative to CuClarabel the increment is feature-level
(\Cref{tab:gpu-vs-cu}), not a claimed wall-clock advantage.
\batchdiff{} is at parity with a warm \diffsolver{} loop up to
$B{=}16$ on the measured instances (\Cref{sec:exp-batch}), so it is
a GIL-free architecture rather than a performance contribution.

\paragraph{Open problems.}
Completing device updates of $P$ and $A$ without a host permutation;
an optional mixed-precision refinement path that does not touch the
stored adjoint; multi-GPU sharding of \batchdiff{} for trajectory
counts beyond one device's memory; and timings on Netlib LP, CBLIB,
and entropy- or SDP-dominated models.
A theory of the Chen--Goulart--Jones operator in equilibrated
coordinates, including its interaction with four-tier sparse SOC
scaling, is also open (\Cref{rem:cgj-scope}).
The GPU SOC cap is an optional variant, neither a sufficient nor
a necessary cause of the GPU sequential $R>1$ cell.

%% file: sections/appendix.tex
\section{Mathematical derivations}
\label{app:barriers}

This appendix collects the derivations that the main text cites but
does not interrupt for: the barrier functions behind the smoothing
operator, the proof of \Cref{prop:vjp} (the linearization from $F$
to $\Ft$ with the assembly of $g$ and $(w_1,w_2)$), and the
rank-two SOC projection derivative used in the adjoint.

\subsection{Barriers used by the smoother}
\label{app:lhscb}

The smoothing operator \eqref{eq:smoothing-dual} uses the standard
logarithmically homogeneous barriers \cite{nesterov1997selfscaled,chen2025warmstart}:
\begin{align*}
f_{\R_+}(s) &= -\sum_i \log s_i, \\
f_{\soc}(t,u) &= -\log(t^2-\|u\|_2^2), \\
f_{\exp}(x,y,z) &= -\log\bigl(y\log(z/y)-x\bigr) - \log y - \log z, \\
f_{\mathrm{pow}}^\alpha(x,y,z) &= -\log(x^{2\alpha} y^{2-2\alpha}-z^2) - (1-\alpha)\log x - \alpha\log y, \\
f_{\Splus}(S) &= -\log\det S.
\end{align*}
The zero cone is handled by the Moreau rule of \Cref{sec:smoothing-impl}
rather than a barrier.

\subsection{From $F$ to $\Ft$: linearization, $g$, and $(w_1,w_2)$}
\label{app:reduced-jacobian}

This subsection is the proof of \Cref{prop:vjp}.
The Minty form of \eqref{eq:kkt} is $G(x,z,s;\theta)=0$ with
\begin{equation}
\label{eq:G-minty}
G
=
\begin{bmatrix}
\dataP x + \dataA\transp z + \dataq \\
\dataA x + s - \datab \\
z - \proj_{\Kstar}(z-s)
\end{bmatrix}.
\end{equation}
Write $v=z-s$ and $u=dv=dz-ds$.
Euclidean projection gives $dz=\Dpi u$ and
$ds=-(I-\Dpi)u$, using $\Dpi^\top=\Dpi$.
Linearizing the first two blocks at a computed solution,
\begin{align}
\label{eq:dG}
dG_1
&=
\dataP\,dx + \dataA\transp\Dpi\,u
+ d\dataP\,x + d\dataq + d\dataA\transp z,
\\
dG_2
&=
\dataA\,dx - (I-\Dpi)\,u
- d\datab + d\dataA\,x.
\end{align}
The reduced residual used in the code is $(G_1,-G_2)$, so the
map on $(dx,u)$ is exactly \eqref{eq:F}:
the second row is $-\dataA$ and $I-\Dpi$, not $+\dataA$ and
$-(I-\Dpi)$.
(The overall sign of $G_2$ is a convention; it is absorbed by
the sign of $w_2$ and of $\partial L/\partial\datab$.)

Block transpose of the unsymmetric $F$ is then mechanical:
\begin{equation}
\label{eq:Ft-blocks}
F^\top
=
\begin{bmatrix}
F_{11}^\top & F_{21}^\top \\
F_{12}^\top & F_{22}^\top
\end{bmatrix}
=
\begin{bmatrix}
\dataP & (-\dataA)^\top \\
(\dataA\transp\Dpi)^\top & (I-\Dpi)^\top
\end{bmatrix}
=
\begin{bmatrix}
\dataP & -\dataA\transp \\
\Dpi\transp \dataA & I-\Dpi\transp
\end{bmatrix},
\end{equation}
which is \eqref{eq:Ft}.
The $(1,2)$ entry $-\dataA\transp$ is $F_{21}^\top$, not
$F_{12}^\top$.
Comparing it to $F_{12}=\dataA\transp\Dpi$ is a block-index
error: those two blocks are not transposes of each other.

The upstream loss differentiates as
$dL=g_x^\top dx+g_z^\top dz+g_s^\top ds$.
Substituting $dz$ and $ds$,
\begin{equation}
\label{eq:dL-u}
dL
=
g_x^\top dx
+
\bigl(\Dpi\transp g_z-(I-\Dpi)^\top g_s\bigr)^\top u.
\end{equation}
Symmetry of $\Dpi$ collapses the second factor to
$\Dpi(g_z+g_s)-g_s$, which is the lower block of \eqref{eq:g}.
When only $g_x$ is supplied, that block is zero.
Write $g=(g_x;\,g_u)$ for this pair.
The implicit-function identity
$F\begin{bmatrix}dx\\u\end{bmatrix}
+(\partial G_{\mathrm{red}}/\partial\theta)\,d\theta=0$
then gives
$dL=-g^\top F^{-1}(\partial G_{\mathrm{red}}/\partial\theta)\,d\theta$.
The convention $\Ft w=-g$ means $w=-F^{-\top}g$, so
$dL=w^\top(\partial G_{\mathrm{red}}/\partial\theta)\,d\theta$.
Here $G_{\mathrm{red}}=(G_1,-G_2)$, hence
$\partial G_{\mathrm{red}}/\partial\dataq=(I;0)$ and
$\partial G_{\mathrm{red}}/\partial\datab=(0;I)$, which is
\eqref{eq:grads}: $w_1=\partial L/\partial\dataq$ and
$w_2=\partial L/\partial\datab$.
The formulae for $\partial L/\partial\dataP$ and
$\partial L/\partial\dataA$ are the same identity applied to
the $d\dataP\,x$ and $d\dataA$ terms of \eqref{eq:dG}.
A singular face direction of a boundary SOC is removed from
$g$ before the solve; it is not a fourth Woodbury column
(\Cref{app:soc-dpi}).
Invertibility of $F$ gives existence and uniqueness of $w$, which
completes the proof of \Cref{prop:vjp}.

\subsection{Rank-two SOC projection derivative}
\label{app:soc-dpi}

On the boundary of a second-order cone, Euclidean projection is
differentiable in the tangent space orthogonal to the generator of
the face.
Writing $v=(v_0,v_1)$, $u=v_1/\|v_1\|_2$,
$p=(1,u)$, $q=(1,-u)$, $\hat u=(0,u)$, and
$\tau=v_0/\eta$, where $\eta=\sqrt{v_0^2+\|v_1\|_2^2}$ is the
Nesterov--Todd scaling norm of $v$ stored by the solver,
the code uses
\begin{equation}
\label{eq:soc-dpi}
\Dpi
=
D_d
+
V D_a^{-1} V^\top,
\qquad
V
=
\bigl[p/2,\; q/2,\; \sqrt{|\tau|/2}\,\hat u\bigr],
\qquad
D_a
=
\mathrm{diag}\bigl(+1,-1,-\mathrm{sign}\,\tau\bigr),
\end{equation}
with $D_d=\mathrm{diag}\bigl(\tfrac12,\,\tfrac{1+\tau}{2},\ldots\bigr)$.
Expanding the Woodbury product gives three signed rank-one terms
\begin{equation}
\label{eq:soc-dpi-expand}
V D_a^{-1} V^\top
=
\tfrac14 pp^\top
-
\tfrac14 qq^\top
-
\tfrac{\tau}{2}\,\hat u\hat u^\top.
\end{equation}
The column space is two-dimensional, not three:
$\hat u=\tfrac12(p-q)$, so the third generator is linearly
dependent on the first two.
Three columns are kept anyway because (i)~$D_a$ stays diagonal
(a thin two-column SVD would replace it by a dense $2\times 2$
block and a $\tau$-dependent pattern), and (ii)~the expanded
KKT always reserves three auxiliary indices, so Identity / Zero
/ RankTwo share one sparsity pattern (three slots are always
reserved).
When $|\tau|<10^{-16}$ the third column is zeroed and only two
are live.
This is distinct from the rank-two \emph{forward} sparse-SOC
expansion of \Cref{sec:ipm}.
\cpu{} stores the expansion in the sparse pattern of $F_e$;
\gpu{} applies it matrix-free except for the sparse KKT fill.
A separate nullspace step removes the singular face direction
from $g$ before the adjoint solve; that direction is not a
fourth column.
Identity and zero projections (strict dual interior / polar
interior) remain purely diagonal.

\section{Maros--M\'esz\'aros QP subset}
\label{sec:exp-maros}

Eighteen problems from the Maros--M\'esz\'aros convex QP collection
\cite{maros1999}, loaded as YimingYAN MAT box QPs
($\tfrac12 x^\top Qx+c^\top x$, $r_\ell\le Ax\le r_u$,
$\ell\le x\le u$) and converted to Clarabel $(P,q,A,b,\mathcal{K})$
with a zero cone for equalities and a nonnegative cone for finite
inequalities and bounds.
The subset mixes tiny Hock--Schittkowski QPs with medium sparse
instances (MOSARQP1, LASER, YAO, CONT-050, CVXQP1\_M); the
138-problem corpus is not claimed.
Tolerances match \Cref{sec:experiments} ($10^{-8}$, 200 IPM
iterations, IR off).
\gpu{} is timed through the native forward solver, not
\diffsolver{}.
Mosek and Gurobi use their native quadratic models.
\Cref{tab:maros-sub} reports the mean of ten independent
re-solves on a reused object, each from the default start with no
warm start.
CuClarabel discards a session-level dummy QP and the first
setup+solve of each problem (Julia JIT / method specialization);
the ten subsequent solves follow that.

\input{tables/tab_maros_sub}

Clarabel.rs is the fastest of the three on the tiny and
small instances (HS35--CVXQP3\_S, typically $0.01$--$2$\,ms).
\gpu{} subsequent time on those problems is $0.6$--$8$\,ms: launch
and host--device overhead, not KKT work.
CuClarabel's post-JIT subsequent times on the same tiny QPs remain
$28$--$109$\,ms (and $64$--$72$\,ms on the AlmostSolved HS35/HS51
cells): that floor is Julia runtime plus a failed or
high-iteration solve, not a KKT comparison with \gpu{}.
The ranking changes with size.
On CONT-050 ($n{=}2597$) the subsequent means are
36\,ms Clarabel.rs, 42\,ms \cpu{}, 20\,ms \gpu{}, 58\,ms
CuClarabel, 48\,ms Mosek, 17\,ms Gurobi (all solved, nine IPM
iterations on Clarabel.rs, \cpu{}, and \gpu{}).
On CVXQP1\_M ($n{=}1000$) they are 80\,ms, 40\,ms, 24\,ms,
311\,ms, 217\,ms, 57\,ms.

Two instances are not a clean six-way solve.
On QPCSTAIR, Mosek returns a primal-infeasibility certificate
after 31 iterations; Clarabel.rs, \cpu{}, \gpu{}, CuClarabel, and
Gurobi agree on the objective $6.204387{\times}10^{6}$ to relative
$10^{-9}$.
On YAO, Mosek and Gurobi agree at $-75.42$.
The Clarabel.rs, \cpu{}, and \gpu{} cold solves do \emph{not} hit the
200-iteration cap: Clarabel.rs and CuClarabel return
AlmostSolved in 33--34 iterations at objective $-183$,
\cpu{} returns AlmostSolved (status~4) in 79 iterations at
$-149$, and \gpu{} stops with insufficient progress in 69
iterations at $-150$.
We did not record $\kappa$ of the KKT matrix; the evidence we
have is that the homogeneous IPM fails to meet the $10^{-8}$
gap/feasibility tolerances on this instance and reports a
different, more negative objective than the commercial solvers.
We treat YAO as a robustness failure of Clarabel.rs, \cpu{}, and \gpu{} on this
subset, not as a timing point, and omit it from
\Cref{tab:maros-sub}.
On CVXQP1\_M, Clarabel.rs, \cpu{}, \gpu{}, and Gurobi agree with
each other to $10^{-8}$ relative, while Mosek's objective lies
$7{\times}10^{-5}$ away---the same four-versus-one pattern as the
$n{=}20000$ portfolio QP (\Cref{tab:portfolio-acc}), and again
without recorded KKT residuals.
CuClarabel is AlmostSolved there at $7.55{\times}10^{5}$ and is
not part of that agreement.
Several tiny QPs also carry a dagger in the CuClarabel column
(HS35, HS51, DUAL1, DUAL4).

This subset does not overturn the portfolio ranking, and it does
not replace a full Maros--M\'esz\'aros campaign.
It shows that on these 17 timed QPs the small-$n$ subsequent leader is
Clarabel.rs, that \gpu{} subsequent time becomes competitive once
$n$ is a few thousand, that CuClarabel's post-JIT subsequent
means on this machine sit at tens to a few hundred milliseconds
on the same subset---subject to the cuDSS 0.8.0 versus 0.7.1
caveat of \Cref{sec:exp-rules}---and that YAO is a documented
failure of the Clarabel IPM at the stated tolerance.

\subsection{SOC-cap isolation}
\label{app:soc-cap}
\Cref{tab:soc-cap} is the SOC-cap $2{\times}2$:
CPU floor / CPU+cap / GPU floor / GPU cap.
CPU stays at $R<1$ in both columns; GPU stays at $R>1$, and
the GPU floor is worse than the production cap.

\input{tables/tab_soc_cap}

\section{Reproducibility notes}
\label{app:repro}
\label{app:settings}

Clarabel.rs, \cpu{}, and \gpu{} share the Settings field names of
the official solver
\cite{clarabelsettings}.
Those defaults ``are configured to give good general performance
for most problems'' on binary64 data
\cite{clarabelsettings};
\Cref{sec:exp-rules} keeps the recommended tolerances and
\texttt{max\_iter} and overrides only the rows marked below.
Meanings in the first block are those of the official table.

\begin{table}[h]
\centering
\caption{Optimizer settings that the tables depend on.
Official names, recommended defaults, and descriptions follow
\cite{clarabelsettings}.
The protocol column is \Cref{sec:exp-rules}.
The second block has no official Clarabel.rs analogue.}
\label{tab:settings}
\small
\setlength{\tabcolsep}{3.2pt}
\resizebox{\textwidth}{!}{%
\begin{tabular}{@{}>{\raggedright\arraybackslash}p{3.55cm}ccp{5.55cm}@{}}
\toprule
Setting & Official & Protocol & Meaning \\
\midrule
\multicolumn{4}{@{}l}{\emph{Official Clarabel Settings}
\cite{clarabelsettings}} \\
\texttt{max\_iter} & 200 & 200
  & maximum number of iterations \\
\texttt{tol\_gap\_abs} & $10^{-8}$ & $10^{-8}$
  & absolute duality gap tolerance \\
\texttt{tol\_gap\_rel} & $10^{-8}$ & $10^{-8}$
  & relative duality gap tolerance \\
\texttt{tol\_feas} & $10^{-8}$ & $10^{-8}$
  & primal and dual feasibility tolerance \\
\texttt{equilibrate\_enable} & true & true$^{\dagger}$
  & data equilibration pre-scaling \\
\makecell[l]{\texttt{iterative\_refinement}\\[-0.15em]\texttt{\_enable}}
  & true & QP off; SOCP on
  & KKT direct solve with iterative refinement \\
\makecell[l]{\texttt{iterative\_refinement}\\[-0.15em]\texttt{\_max\_iter}}
  & 10 & 2
  & iterative refinement maximum iterations \\
\makecell[l]{\texttt{chordal\_decomposition}\\[-0.15em]\texttt{\_enable}}
  & true & false
  & enable chordal decomposition \\
\texttt{max\_threads} & 0 (auto) & 1
  & host threads for multithreaded KKT solves \\
\midrule
\multicolumn{4}{@{}l}{\emph{CPU/GPU-only} (no official analogue)} \\
IR latch & --- & CPU AUTO; GPU AUTO2
  & latches near tolerance (\Cref{sec:gpu-cudss}) \\
Smoothing warm start & --- & off unless named
  & not an instance-cache flag \\
SOC $\mubar$ cap & --- & CPU off; GPU on
  & optional; \Cref{tab:soc-cap} \\
Two-RHS KKT solve & --- & on
  & affine + constant in one factor \\
Hybrid cuDSS execute & --- & AUTO
  & CPU root supernodes; device-only if shard $>8$ \\
Prefactor at setup & --- & on
  & fused analysis+factor at construction \\
CPU adjoint & --- & KLU
  & unsymmetric $F^\top$; optional QDLDL lift \\
\bottomrule
\end{tabular}%
}
\\[0.35em]
{\footnotesize $^{\dagger}$Off in one row of \Cref{tab:ablate}.
SOCP IR: Clarabel.rs \texttt{true}, \cpu{} AUTO, \gpu{} AUTO2.}
\end{table}

%% file: tables/tab_maros_sub.tex
\begin{table}[t]
\centering
\caption{Maros--M\'esz\'aros subset, subsequent mean of ten independent re-solves on a reused object, each from the default start (no warm start). YimingYAN MAT box QPs; gap/feasibility $10^{-8}$. YAO is omitted (failure of Clarabel.rs, \cpu{}, and \gpu{}; see the text). \gpu{} is the native forward solver (not \diffsolver{}). Mosek and Gurobi use native quadratic models. Mosek deletes the interior-point solution between solves; Gurobi resets to the default start. A dagger marks a non-optimal status. CuClarabel times are post-JIT (session dummy QP and the first run of each problem discarded) and use CUDSS.jl~0.6 / \texttt{libcudss}~0.7.1; \gpu{} uses conda \texttt{libcudss}~0.8.0.}
\label{tab:maros-sub}
\small

\setlength{\tabcolsep}{2.4pt}
\begin{tabular}{@{}l r r ccccc@{}}
\toprule
Problem & $n$ & \makecell{Clarabel\\rs} & \makecell{Clarabel\\CPU} & \makecell{Clarabel\\GPU} & \makecell{CuCla-\\rabel} & Mosek & Gurobi \\
\midrule
HS35 & 3 & 0.012 (7) & 0.092 (7) & 2.92 (7) & 72.3 (22)$^{\dagger}$ & 0.20 (6) & 0.19 (9) \\
HS51 & 5 & 0.004 (1) & 0.035 (1) & 0.58 (1) & 64.4 (24)$^{\dagger}$ & 0.13 (3) & 0.18 (2) \\
HS53 & 5 & 0.014 (6) & 0.10 (7) & 2.62 (6) & 28.1 (8) & 0.19 (4) & 0.19 (8) \\
QPTEST & 2 & 0.013 (8) & 0.10 (8) & 3.27 (8) & 77.3 (18) & 0.20 (6) & 0.18 (11) \\
QAFIRO & 32 & 0.094 (14) & 0.26 (15) & 5.90 (14) & 46.7 (14) & 0.44 (12) & 0.28 (12) \\
DUAL1 & 85 & 2.06 (12) & 1.63 (12) & 6.80 (13) & 50.5 (15)$^{\dagger}$ & 2.98 (12) & 3.67 (16) \\
DUAL4 & 75 & 1.54 (12) & 1.08 (11) & 6.24 (12) & 721 (200)$^{\dagger}$ & 2.34 (12) & 2.66 (15) \\
QPCBLEND & 83 & 0.44 (17) & 1.30 (19) & 7.63 (17) & 57.7 (17) & 1.27 (14) & 1.02 (19) \\
QPCSTAIR & 467 & 7.00 (24) & 13.7 (23) & 24.4 (24) & 104 (24) & 30.2 (31)$^{\dagger}$ & 10.3 (23) \\
CVXQP1\_S & 100 & 0.40 (9) & 1.28 (10) & 4.49 (9) & 95.0 (27) & 2.45 (14) & 1.53 (12) \\
CVXQP2\_S & 100 & 0.38 (10) & 0.99 (11) & 4.81 (10) & 109 (22) & 2.13 (15) & 0.95 (12) \\
CVXQP3\_S & 100 & 0.58 (11) & 1.01 (12) & 5.65 (11) & 61.6 (18) & 2.37 (17) & 2.67 (12) \\
MOSARQP1 & 2500 & 7.19 (10) & 9.79 (10) & 11.2 (10) & 53.8 (12) & 29.3 (15) & 9.92 (16) \\
PRIMAL1 & 325 & 2.65 (10) & 4.03 (10) & 8.47 (10) & 37.9 (10) & 5.82 (10) & 4.38 (19) \\
LASER & 1002 & 2.99 (11) & 2.53 (9) & 9.62 (11) & 57.1 (14) & 18.7 (15) & 17.6 (22) \\
CONT-050 & 2597 & 36.1 (9) & 42.2 (9) & 20.4 (9) & 58.0 (9) & 48.0 (7) & 17.4 (10) \\
CVXQP1\_M & 1000 & 80.2 (10) & 39.5 (12) & 23.7 (10) & 311 (56) & 217 (28) & 56.8 (14) \\
\bottomrule
\end{tabular}
\end{table}

%% file: tables/tab_soc_cap.tex
\begin{table}[t]
\centering
\caption{SOC-cap $2{\times}2$ on the $n{=}1000$, $k{=}50$ SOCP
(seed~42).
CPU floor-only is the production CPU map; CPU+cap applies the
relative cap
$\mubar_{\mathrm{soc}}\le 10^{-4}\max(\|c_1\|_2^2/4,s_0^2)$
on the same QDLDL backend.
GPU floor turns that cap off on the production GPU path;
GPU cap is the production GPU map.
One-shot cells are a single relative $\mu$ jump; sequential cells
are the $T{=}20$ geometric-mean ratio \eqref{eq:Riter} at
$10^{-2}$.
CPU stays at $R<1$ with or without the cap; GPU stays at $R>1$
with or without it, and turning the cap off makes the GPU
sequential cell worse ($R=1.34$ versus $1.20$).}
\label{tab:soc-cap}
\small
\begin{tabular}{@{}llccc@{}}
\toprule
Protocol & Map & cold & smooth & $R$ \\
\midrule
one-shot $10^{-4}$ & CPU floor & 21 & 13 & $0.62$ \\
one-shot $10^{-4}$ & CPU+cap & 21 & 13 & $0.62$ \\
one-shot $10^{-4}$ & GPU floor & 20 & 21 & $1.05$ \\
one-shot $10^{-4}$ & GPU cap & 20 & 15 & $0.75$ \\
one-shot $10^{-2}$ & CPU floor & 20 & 16 & $0.80$ \\
one-shot $10^{-2}$ & CPU+cap & 20 & 17 & $0.85$ \\
one-shot $10^{-2}$ & GPU floor & 19 & 24 & $1.26$ \\
one-shot $10^{-2}$ & GPU cap & 19 & 22 & $1.16$ \\
sequential $10^{-2}$ & CPU floor & $19.95$ & $15.3$ & $0.77$ \\
sequential $10^{-2}$ & CPU+cap & $19.95$ & $16.95$ & $0.85$ \\
sequential $10^{-2}$ & GPU floor & $18.0$ & $24.0$ & $1.34$ \\
sequential $10^{-2}$ & GPU cap & $18.0$ & $21.65$ & $1.20$ \\
\bottomrule
\end{tabular}
\end{table}